\documentclass[12pt,leqno]{article}
\usepackage{amssymb}
\usepackage{verbatim}
\usepackage{amscd}
\usepackage{amsmath}
\usepackage{tabularx}
\usepackage[dvips]{graphicx}
\usepackage[arrow,matrix,tips,frame,color,line,poly]{xy}
\usepackage{theorem}
\newtheorem{theorem}[equation]{Theorem}
\newtheorem{corollary}[equation]{Corollary}
\newtheorem{definition}[equation]{Definition}
\newtheorem{lemma}[equation]{Lemma}

\newtheorem{proposition}[equation]{Proposition}

{\theorembodyfont{\rmfamily}

\newtheorem{remark}[equation]{Remark}

\newtheorem{example}[equation]{Example}
}
\newcommand{\qed}{\hfill $\square$ \medskip}

\newenvironment{proof}[1][Proof]{\noindent\textbf{#1.} }{\
\rule{0.5em}{0.5em}\medskip}

\renewenvironment{proof}[1][Proof]{\noindent\textbf{#1.} }{\
\qed}

\renewcommand{\exp}{\mathrm{exp}}

\renewcommand{\O}{\mathcal O}
\newcommand{\Ost}{\mathcal{O^{\text{st}}}}

\newcommand{\Orb}{\text{Orb}}
\newcommand{\Orbst}{\text{Orb}^{st}}
\newcommand{\st}{\text{st}}

\newcommand{\caS}{\mathcal{S}}

\newcommand{\caT}{\mathcal{T}}

\newcommand{\Ztwo}{\Z/2\Z}
\newcommand{\wt}{\widetilde}

\newcommand{\bs}{\backslash}

\newcommand{\F}{\mathbb F}
\newcommand{\R}{\mathbb R}

\newcommand{\Gad}{G_{\ad}}

\newcommand{\C}{\mathbb C}
\newcommand{\Z}{\mathbb Z}

\renewcommand{\F}{\mathbb F}

\newcommand{\tH}{\wt{H}}
\newcommand{\tS}{\wt{S}}

\newcommand{\tG}{\widetilde G}
\newcommand{\tg}{\tilde g}

\newcommand{\tM}{\widetilde M}
\newcommand{\ch}[1]{#1^\vee}

\renewcommand{\t}{\mathfrak t}
\renewcommand{\a}{\mathfrak a}

\newcommand{\liem}{\mathfrak m}

\newcommand{\SL}{\mathrm{SL}}

\newcommand{\h}{\mathfrak h}
\renewcommand{\k}{\mathfrak k}

\newcommand{\tsl}{\widetilde{\SL}(2)}

\newcommand\inv{^{-1}}

\newcommand{\g}{\mathfrak g}
\renewcommand{\j}{\mathfrak j}

\newcommand{\Ind}{\mathrm{Ind}}

\newcommand{\Aut}{\text{Aut}}
\newcommand{\diag}{\text{diag}}

\renewcommand{\int}{\text{int}}
\newcommand{\Hom}{\text{Hom}}
\newcommand{\Tr}{\text{Tr}}
\newcommand{\Res}{\text{Res}}
\newcommand{\ad}{\text{ad}}

\newcommand{\Ad}{\text{Ad}}

\newcommand{\Cent}{\text{Cent}}
\newcommand{\Norm}{\text{Norm}}

\newcommand{\Ker}{\text{Ker}}

\renewcommand{\int}{\text{int}}
\newcommand{\sgn}{\text{sgn}}

\renewcommand{\sec}[1]{\section{#1}
\renewcommand{\theequation}{\thesection.\arabic{equation}}
  \setcounter{equation}{0}}

\newcommand{\Lift}{\mathrm{Lift}}

\newcommand{\zetacx}{\zeta_{\text{cx}}}
\newcommand{\Wsharp}{W_{\#}}
\newcommand{\cd}{CD}
\newcommand{\cdg}{CD_{g}}
\def\t{\mathfrak t}
\def\h{\mathfrak h}
\def\k{\mathfrak k}
\def\p{\mathfrak p}

\def\a{\mathfrak a}
\def\g{\mathfrak g}

\def\z{\mathfrak z}

\newcommand\oG{\overline G}
\newcommand\oJ{\overline J}

\newcommand\oH{\overline H}

\newcommand\oT{\overline T}
\newcommand\oM{\overline M}

\newcommand\oA{\overline A}
\newcommand\op{\overline p}

\def\tsl2{\widetilde{\rm{SL}(2)}}

\renewcommand{\tilde}{\widetilde}

\begin{document}

\title{Lifting of Characters for Nonlinear Simply Laced Groups}
\author{Jeffrey Adams and Rebecca Herb}
\maketitle

{\renewcommand{\thefootnote}{}
\footnote{The first author was supported in part by  National Science
Foundation Grant \#DMS-0554278.}}

\section*{Abstract}

One aspect of the Langlands program for linear groups
is  lifting of characters, which relates virtual representations on
a group $G$ with those on an endoscopic group for $G$. The goal of
this paper is to extend this theory to nonlinear two-fold covers of
real groups in the simply laced case.
Suppose $\tG$ is a two-fold cover of a real reductive group $G$.
The main result is that
there is an operation, denoted $\Lift_G^{\tG}$, 
taking a stable virtual character of $G$ to $0$ or a virtual genuine
character of $\tG$, and $\Lift_G^{\tG}(\Theta_\pi)$ 
may be explicitly computed if $\pi$ is a stable sum of standard modules. 

\sec{Introduction}
\label{s:introduction}

The Langlands program is concerned with representation theory and
automorphic forms of algebraic (linear) groups.  Suppose $\tG$ is a
nonlinear group, for example the metaplectic group $Mp(2n)$, the
nontrivial twofold cover of $Sp(2n)$. It is well known that such
groups play an important role in automorphic forms, although they
don't fit into the formalism of the Langlands program.  The goal of
this paper and of \cite{dualityonerootlength} is to extend some of the
results of the Langlands program to certain nonlinear groups.

Let  $G$ be the real points of a connected reductive algebraic
group. Characters of admissible representations of $G$ have been
studied extensively, and the theory is reasonably complete.  In
particular if $\pi$ is a discrete series representation the
character  $\Theta_\pi$ of $\pi$ is known explicitly \cite{herb_ds}. 
The formula for $\Theta_\pi$
implicitly comes from the theory of transfer and endoscopy, originating
in work of Shelstad and Langlands on the trace formula and automorphic
forms. We would like to extend the theory of transfer of characters to
nonlinear groups. 

Suppose $\tG$ is a nonlinear two-fold cover of $G$, i.e. $\tG$ can
not be realized as a subgroup of $GL(n,\C)$. Examples include the
metaplectic group $Mp(2n,\R)$, the unique connected two-fold cover of
$Sp(2n,\R)$, and the twofold cover of $GL(n,\R)$ of \cite{huang_gln}, \cite{ah}.
There is substantial evidence that character  theory for genuine
representations of $\tG$ (those which do not factor to $G$) may be
reduced to that of a linear group.

The case of $GL(n)$ over any local field of characteristic zero
has been studied by Kazhdan,
Patterson and Flicker (\cite{flicker_gl2}, \cite{kf_metaplectic},
\cite{kp_ihes}). The general philosophy is spelled out in 
\cite{kp_towards}. Some discussion of how to extend these ideas to
general groups is given in \cite{adams:jussieu}, and this is carried out for
$G=Sp(2n,\R)$ in \cite{metalift}. 

Here is an outline of the approach. For now suppose $G$ is the
$\F$-points of a connected, reductive algebraic group defined over a
local field of characteristic $0$, and $\tG$ is a nonlinear two-fold
cover of $G$. Identify the kernel of the covering map
$p:\tG\rightarrow G$ with $\pm1$.
We would like to relate genuine characters of $\tG$ 
to characters of $G$. By character we mean
the character of a representation, viewed as a function on
the regular semisimple elements.

The theory of endoscopy for $G$ relates characters of $G$ to
characters of {\it endoscopic groups} of smaller semisimple rank. A
key part of the theory is a relationship between semisimple conjugacy
classes in $G$ and $H$.

Our theory is modelled on this, with $\tG$ in place of $G$, and $G$
playing the role of $H$.
We define a relation on conjugacy
classes as follows. For $g\in G$ define $\psi(g)=s(g)^2$ where
$s(g)\in\tG$ satisfies $p(\tg)=g$. Then $\psi(g)$ is independent of the
choice of $s(g)$, and this induces a map from conjugacy classes of $G$
to those of $\tG$. This idea goes back to  \cite{flicker_gl2}.

Now suppose $\pi$ is an admissible representation of $G$, with
character $\Theta_\pi$ viewed 
as a function on the  (strongly)
regular semisimple elements $G'$ of $G$ (see Section
\ref{s:admissible}). 
For $\tg\in \tG'$ define
\begin{equation}
\label{e:naive}
\psi^*(\Theta_\pi)(\tg)=\sum_{\{h\in G\,|\,\psi(h)=\tg\}}\Theta_\pi(h).
\end{equation}
This is a conjugation invariant function on $\tG'$. One can
ask whether it is the character of a genuine representation $\wt\pi$
of $\tG$.
More generally it may be a virtual character, i.e. the character of
a finite sum of irreducible representations with (possibly negative)
integral coefficients.

Although this definition is too naive for several reasons, we describe
a basic property of characters which suggest it is at least a good
first approximation. 
Suppose $H$ is a Cartan subgroup of $G$, with inverse image  $\tH$ in
$\tG$. Then $\tH$ is not necessarily abelian, and we let
$Z(\tH)$ be its center. The fact that $\tH$ is not abelian is a major
reason why the  representation theory of $\tG$ at least appears to be more
complicated than that of $G$. However, at least as far as
character theory goes, the situation is actually very simple:  the
character $\Theta_{\wt\pi}$ of an admissible genuine representation of
$\tG$ satisfies $\Theta_{\wt\pi}(\tg)=0$ for $\tg\not\in Z(\tH)$
(Lemma \ref{l:relevant}).
It is easy to see that $\psi(H)\subset Z(\tH)$, of finite index,
so at least from this point of view \eqref{e:naive} is similar to the
character of a genuine representation of $\tG$.

An obvious shortcoming of\eqref{e:naive} is that
it is not necessarily the case that
$\psi^*(\Theta_\pi)(-\tg)=-\psi^*(\Theta_\pi)(\tg)$,
an obvious
requirement if $\psi^*(\Theta_\pi)$ is to be the character of a
genuine representation. 
It would be better to define $\phi (h) = h^2$, 
sum over $h$ satisfying
$\phi (h)=p(\tg)$, and modify the definition by some genuine function $\mu$ of
$\tG$:
\begin{equation}
\phi^*(\Theta_\pi)(\tg)=\sum_{\{h\in G\,|\, \phi (h)=p(\tg)\}}\mu(\tg)\Theta_\pi(h).
\end{equation}
We then need to choose $\mu$ appropriately.

There are also some Weyl denominators to take into account. 
Suppose $\Phi^+$ is a set of
positive roots of $H$ in $G$. For $h\in H$ one version of the  Weyl denominator is
$|D(h)|^{\frac12}=|\prod_{\alpha\in\Phi^+}(1-\alpha^{-1}(h))||e^\rho(h)|$
(cf.~\eqref{e:Deltas}). This is
independent of $\Phi^+$ and is well defined.
The character of any representation has a particular form when
muliplied by $|D(h)|^{\frac12}$; for the left hand side of
\eqref{e:naive} to have a chance of having this form we should
multiply each term on the right by the quotient of Weyl denominators
$|D(h)|^{\frac12}/|D(\tg)|^{\frac12}$.

Putting these two considerations  together we look for a function $\Delta$ on $G'\times
\tG'$ satisfying the following conditions:
\begin{equation}
\label{e:tf}
\begin{aligned}
\Delta(h,\tg)=0&\text{ unless }p(\tg)=\phi(h)\\
|\Delta(h,\tg)|&=|D(h)|^{\frac12}/|D(\tg)|^{\frac12}\\
\Delta(xhx\inv,\tilde x\tg\tilde x\inv)&=
\Delta(h,\tg)\quad (\tilde x\in \tG, x=p(\tg))\\
\Delta(h,-\tg)&=-\Delta(h,\tg).
\end{aligned}
\end{equation}
Given such a function we define
\begin{equation}
\label{e:better}
\Lift_G^{\tG}(\Theta_\pi)(\tg)=\sum_{\{h\in G\,|\, \phi (h)=p(\tg)\}}\Delta(h,\tg)\Theta_\pi(h).
\end{equation}
Now $\Lift_G^{\tG}(\Theta_\pi)$ is
a conjugation invariant, genuine function on $\tG$, and
is a reasonable candidate for the character of a virtual representation.
The latter condition
amounts to further conditions on $\Delta$.
If this holds, one can {\it define} $\Lift_G^{\tG}(\pi)$ to be the
virtual representation of $\tG$ whose character is equal to
$\Lift_G^{\tG}(\Theta_\pi)$:
$\Theta_{\Lift_G^{\tG}(\pi)}=\Lift_G^{\tG}(\Theta_\pi)$.

For a version of this for covering groups of $GL(n)$  
see \cite{kp_ihes}.

Formula \eqref{e:tf} is analogous with the main character identity in
the theory of endoscopy.  For example see \cite{shelstad_annalen}. In that setting
$\Delta$ is a {\it transfer factor}, and by analogy we use the same
terminology here.  As in the theory of endoscopy, correctly defining
transfer factors is a difficult part of the theory.

We now discuss another less obvious consideration which arises. Recall
$\Theta_\pi$ is conjugation invariant:
$\Theta_\pi(g)=\Theta_\pi(g')$ if $g,g'\in G'$ and $g'=xgx\inv$
for some $x\in G$.
Following Langlands and Shelstad we say $\pi$ and
$\Theta_\pi$ are {\it stable} if the following stronger condition
holds:
$\Theta_\pi(g)=\Theta_\pi(g')$
if $g,g'\in G'$ and $g'=xgx\inv$ for some $x\in G(\overline\F)$.

For $GL(n)$ every conjugation invariant function is automatically stable,
so stability plays no role in 
\cite{kf_metaplectic}, \cite{flicker_gl2} and \cite{kp_ihes}. 

For guidance we consider the case of $G=Sp(2n,\R)$, $\tG=Mp(2n,\R)$,
as discussed in \cite{metalift}.
Let $G'=SO(n+1,n)$. In this case the map $\phi$ above is replaced by a
bijection between (strongly) regular  semisimple conjugacy classes in
$G$ and $G'$.
Suppose $\pi$ is a {\it stable} representation of
$G'$. 
The main result of \cite{metalift} is that, for appropriate definition
of the transfer factor $\Delta$, if we define
\begin{equation}
\label{e:deflift}
\Lift_{G'}^{\tG}(\Theta_\pi)(\tg)=\Delta(h,\tg)\Theta_\pi(h)\quad(\phi(h)=p(\tg))
\end{equation}
then $\Lift_{G'}^{\tG}(\Theta_\pi)$ is the character of a genuine
virtual representation of $\tG$.

This suggests we should relate characters of $\tG$ with {\it stable} characters  
of real forms of the {\it dual group} $\ch G(\C)$. 
We expect the general theory to have this form.
In this paper we
restrict ourselves to simply laced groups.
This avoids a number of technical complications, and 
the distinction
between $G(\C)$ and $\ch G(\C)$ is less critical.
We may now state a special case of the main result.

\begin{theorem}[Theorem \ref{t:liftstandard} and Corollary
\ref{c:main}]
\label{t:intromain}
Suppose $G(\C)$ is a connected, reductive, simply laced complex group,
with real points $G$.
We assume the derived group of $G(\C)$ is acceptable. 
Suppose $\tG$ is a  admissible cover of $G$ (cf.~Section \ref{s:notation} and
Definition \ref{defadm}).
Then we can define the transfer factor$\Delta(h,\tg)$, satisfying
\eqref{e:tf}, such that for all stable admissible representations $\pi$ of
$G$, 
\begin{equation}
\Lift_G^{\tG}(\Theta_\pi)(\tg)=\sum_{\{h\in G\,|\, \phi (h)=p(\tg)\}}\Delta(h,\tg)\Theta_\pi(h) 
\end{equation}
is the character of a genuine virtual representation $\wt\pi$ of
$\tG$, or $0$.
We say
 $\wt\pi$ is the {\it
  lift} of $\pi$, and write $\wt\pi=\Lift_G^{\tG}(\pi)$.

If $\pi$ is a stable sum of standard modules we compute
$\Lift_G^{\tG}(\pi)$ explicitly.
\end{theorem}

We note that unlike the case of $Mp(2n,\R)$, the role of stability is
very subtle in the simply laced case. It appears in the proof of
Theorem \ref{inveig}; see Remark \ref{r:stable}. 

Theorem \ref{t:intromain} is formally similar to transfer in the setting of endoscopic
groups. 
For example see 
 \cite[Lemma 4.2.4]{shelstad_annalen}. More precisely
$\Lift_G^{\tG}$ is analogous to the simplest case of endoscopy:
transfer from the quasisplit form $G_{qs}$ of $G$ to $G$
\cite{shelstad_inner}.
As we will see below it is often possible to
obtain a single irreducible representation of $\tG$ as a lift, and
there is no natural notion of stability for $\tG$. This suggests that
in the simply laced case this is the only notion of lifting
which is needed. We note that there {\it is} a useful notion of stability for 
$Mp(2n,R)$ \cite{metalift}, and David Renard has defined a family of
``endoscopic groups'' for $Mp(2n,\R)$ and proved lifting results for
them \cite{renard_endoscopy}. We believe this is the only situation in which this is either
necessary or possible.
We plan to return to the two root length case in another paper.

\begin{example}
\label{ex:sl2}
Let $G=SL(2,\R)$ and let $\tG$ be the unique non-trivial two folder
cover of $G$.
Let $B$ be a Borel subgroup of $G$ and write $B=AN$ with $A\simeq
\R^\times$.
For $\delta=\pm1$ and $\nu\in\C$ define 
a character $\chi$ of $\R^\times$ by
$\chi(-1)=\delta$ and 
$\chi(x)=x^\nu$ for $x\in\R^+$.
Let $\pi(\delta,\nu)$ be the corresponding principal series
representation of $G$, i.e. $\Ind_{AN}^G(\chi\otimes 1)$ (normalized
induction). 

Now let $\wt A$ be the inverse image of $A$ in $\tG$. 
Then $\wt A\simeq \R^\times\cup
i\R^\times$. 
We may identify
the inverse image of $N$ with $N$,
and let $\wt B=\wt AN$. For $\epsilon\in\pm1$ and $\gamma\in\C$ define
$\wt\chi(i)=\epsilon i$ and $\wt\chi(x)=x^\gamma$ for $x\in\R^+$.
Let $\wt\pi(\epsilon,\gamma)$ be the genuine principal series
representation  $\Ind_{\wt B}^{\tG}(\wt\chi\otimes1)$ of $\tG$.

Let $\pi=\pi(\delta,\nu)$. 
It is easy to see that $\Lift_G^{\tG}(\Theta_\pi)=0$ if $\delta=-1$. 
If $\delta=1$ then by an easy calculation (see Section \ref{s:ps}),
using the transfer factor of Section \ref{s:special} we see
\begin{equation}
\label{e:sl2}
\Lift_G^{\tG}(\pi (1, \nu))=\wt\pi(1,\nu/2)\oplus\wt\pi(-1,\nu/2).
\end{equation}
Note that the image of $\phi$ restricted to $A$ is $A^0$, 
and the character of this Lift is $0$ on $p\inv(-A^0)=i\R^\times$.

This example illustrates several features. 
In this example computing $\Lift_G^{\tG}(\pi)$
essentially reduces to computing $\Lift_A^{\wt A}(\chi)$. This is a
special case of the general situation.

Note that $\wt\pi(\pm1,\nu/2)$ have different central characters (the
center of $\tG$ is $\wt{Z(G)}\simeq\Z/4\Z$), so
$\Lift_G^{\tG}(\pi(1,\nu))$ does not have a central character. Also
note that this character is supported on $\wt{A^0}$, which is
a proper subgroup of $Z(\wt A)$. 
\end{example}

A hint that Theorem \ref{t:intromain} is not the best result possible is seen
by considering the preceding example from another point of view. View
$SL(2,\R)$ as $Sp(2,\R)$, and apply \cite{metalift}. In this case
$G'=SO(2,1)$. There are two principal series representations $\pi(\pm1,\nu)$
of $SO(2,1)$ analogous to those described for $SL(2,\R)$. 
Writing ${\mathcal Lift}$ for the lifting of  \cite{metalift}  we see 
\begin{equation}
\label{e:liftmp2}
{\mathcal Lift}_{SO(2,1)}^{Mp(2,\R)}(\pi(\pm1,\nu))=\wt\pi(\pm1,\nu).
\end{equation}
See \cite[Proposition 15.10]{metalift}. Thus we may obtain each
principal series
representation $\wt\pi(\pm1,\nu)$ of $\wt{SL}(2,\R)$, rather than just
their sum as in \eqref{e:sl2}.

Note that $SO(2,1)$ is isomorphic to $PSL(2,\R)$, the real points of
$PSL(2,\C)$ (also denoted $PGL(2,\R)$). This is the real form of the
adjoint group, and  suggests it should be possible to generalize lifting of
\eqref{e:deflift} to allow $G$ to be replaced with a real form of a
quotient of $G(\C)$ (without changing $\tG$). We revisit the previous
example from this point of view.

\begin{example}
\label{ex:pgl2}
Let $G(\C)=SL(2,\C)$ and $$\oG(\C)=PSL(2,\C)=SL(2,\C)/\pm I.$$ Let
$G=SL(2,\R)$ and let $\oG=PSL(2,\R)$ be the real points of $\oG(\C)$. 
Recall $\oG\simeq PGL(2,\R)\simeq SO(2,1)$. 

We define an orbit correspondence between $\oG$ and $G$ as
follows. For $g\in \oG$ choose an inverse image $s(g)$ of $g$ in
$G(\C)$, and let $\phi(g)=s(g)^2$. This is independent of the choice
of $s(g)$, and $\phi(g)\in SL(2,\R)$.

Let $A(\C)=\{ \diag(z,1/z) \} \subset G(\C)$, and let $\overline A(\C)$ be
its image in $\oG (\C)$. Then $A(\C)$ and $\overline A(\C)$ are defined
over $\R$, and let $A=A(\R), \overline A=\overline A(\R)$.
Both $A$ and $\overline A$ are isomorphic to
$\R^\times$. Note that while the map $A(\C)\rightarrow \overline A(\C)$ is
surjective, the restriction $A(\R)\rightarrow \overline A(\R)$ is
not.

Write $\overline{\diag}(z,\frac1z)$ for the image of $\diag(z,\frac1z)$
in $\overline A(\C)$. 
Let $g=\overline{\diag}(x,\frac1x)\subset\overline A$ with
$x\in\R^\times$. Then $\phi(g)=\diag(x^2,\frac1{x^2})\in A^0$. 
However suppose $y\in\R^\times$ and let
$g=\overline{\diag}(iy,\frac1{iy})\in\overline A$. Then
$\phi(g)=\diag(-y^2,-\frac1{y^2})$. Therefore (unlike in Example
\ref{ex:sl2}) $\phi$ maps $\overline A$ {\it onto} $A$.

This suggests that if we develop a similar lifting theory from $\oG$
to $\tG$, then the lift of a principal series representation will have
support on all of $\wt A$. In fact this is the case: we can define
$\Lift_{\oG}^{\tG}$. We recover \eqref{e:liftmp2} from this point of
view (the difference between $\nu$ and $\nu/2$ on the right hand side
is an issue of normalization):
\begin{equation}
\Lift_{\oG}^{\tG}(\overline \pi(\pm 1,\nu))=\wt\pi(\pm 1,\nu/2).
\end{equation}
\end{example}

\medskip

Motivated by this example, we generalize \eqref{e:deflift} as follows.
Let $G(\C)$ be our given complex group, with real points $G$ and
nonlinear cover $\tG$.  Suppose $C\subset Z(G)$ is a two group and
let $\oG(\C)=G(\C)/C$, with real points $\oG$. For $g\in\oG$ define
$\phi(g)=s(g)^2$ where $s(g)$ is an inverse image of $g$ in $G(\C)$. 
The fact that $C$ is a two-group implies $\phi(g)\in G$, and this is
independent of the choice of $s(g)$. 
It is not necessarily the case that $\phi$ induces a map on conjugacy
classes; however it does define a map on {\it stable} conjugacy
classes, which is consistent with our application. 

We need to make several technical assumptions on $C$ (Definition
\ref{assumptions}); $C=1$ is allowed. Under these assumptions we can
define the transfer factor $\Delta(h,\tg)$ on $\oG'\times\tG$ and
define the lift  $\Lift_{\oG}^{\tG}(\pi)$ of a stable
representation of $\oG$ as in \eqref{e:deflift}. See Definition
\ref{d:lift}.

Now suppose $H$ is a Cartan subgroup of $G$, with inverse image $\tH$ in
$\tG$ and corresponding subgroup $\oH\subset \oG$.  Our assumptions on
$C$ imply $\phi(\oH)\subset p(Z(\tH))$, which is necessary to have a
meaningful theory.  As in the preceding examples, if this image is
large then $\Lift_{\oG}^{\tG}(\pi)$ will have fewer terms in its sum.
This is desirable; it would be nice to have only one term in the sum.
Taking $C$ larger makes this image larger, so we would like to take
$C$ as large as possible. An optimal choice would be to choose $C$ so
that $\phi(\oH)=p(Z(\tH))$ for all Cartan subgroups, in which case
each Lift would consist of a single term. This is not possible in
general. See Example \ref{notsogood}.

\medskip

The proofs follow the same general outline as those of
\cite{metalift}. A hard part of the theory is the definition of
transfer factors. Once transfer factors and lifting are defined we first show that the lift
of a stable invariant eigendistribution is an invariant
eigendistribution. This requires checking that lifting respects the
Hirai matching conditions, and it is here that stability plays a
crucial role. Once this is done it is fairly easy to compute the lift
of a stable sum of discrete series representations. It is also
straightforward to prove that lifting commutes with parabolic
induction. 
This enables 
us to compute the lift of
any stable sum of standard modules. 
In principal, giving the Kazhdan-Lusztig-Vogan polynomials for $\oG$ and
$\tG$, we may
then compute the Lift of any stable virtual representation of $\oG$.

\medskip

Here is an outline of the contents of the paper.

In Section \ref{s:notation} we make some basic definitions and establish some
notation. Admissible triples $(\tG,G,\oG)$ are defined in Section
\ref{s:admissible}. Such a set consists of a nonlinear group $\tG$, a linear group
$G$, and the real form $\oG$ of a quotient of $G(\C)$, satisfying
certain assumptions.
Section \ref{s:admissible} also contains a discussion of some basic structural facts
about Cartan subgroups, and defines the orbit correspondence 
$\phi$. Section \ref{s:cartan} recalls
some standard facts about Cartan subgroups and Cayley 
transforms, and generalizes them to nonlinear groups.
Cayley transforms are an important tool in the theory. 

The basic, and most important case, that of the real points of a
semisimple, simply connected complex group, is discussed in Section
\ref{s:special}. The transfer factors are canonical in this case, 
have a simple form, and this case provides guidance for the general
theory.  We recommmend the reader restrict  to this case the
first time through. 

Sections \ref{s:characters} and \ref{s:transfer} are technically the
most difficult. Section \ref{s:characters}
defines
certain characters of Cartan subgroups and their application to
transfer factors. Transfer factors are defined in Section \ref{s:transfer}. 
Some constants associated to Cartan subgroups are defined and studied
in Section \ref{s:constants}.

With transfer factors in place we define lifting and study its basic
properties in Section \ref{s:lifting}. The case of tori, which is both a good
example and an important special case, is discussed in Section
\ref{s:tori}. Minimal principal series of split groups are covered
in Section \ref{s:ps}, and discrete series on the compact Cartan in
Section \ref{s:exds}.
Section \ref{s:hirai} summarizes some material about invariant
eigendistributions which we will need.
Some details of an extension of Hirai's results which we need are in
the Appendix (Section \ref{app}); this is work of the second author. 

In Section \ref{s:liftinginveig} we use the results of Section \ref{s:hirai} to 
prove that the lift of a stable invariant eigendistribution
in an invariant eigendistribution.

We prove that lifting commutes with parabolic induction in 
Section \ref{s:levi}.
{\it Modified  character data} appropriate to our setting is defined in
Section \ref{s:characterdata}.
This is a (mild) modification of character data and
the Langlands classification due to Vogan \cite{green}. Formal lifting
of (modified) character data is defined in Section \ref{s:liftingofdata}.
This  essentially comes down to lifting applied to a Cartan subgroup.

We compute the lifting of stable discrete series representations in
Section \ref{s:discreteseries}, and of general standard modules in
Section \ref{s:standard}.

These results are closely related to {\it duality} of representations,
also known as Vogan duality, introduced in \cite{ic4}. Duality for
nonlinear covers of simply laced groups is discussed in a paper by
Peter Trapa and the first author \cite{dualityonerootlength}. These
two results grew up together, and we thank Peter Trapa for many
helpful discussions.

\sec{Some Notation}
\label{s:notation}

A simple root system is said to be {\it simply laced} if all roots
have the same length, and an arbitrary root system is simply laced if this holds for
each simple factor. More succinctly a root system is simply laced if 
whenever $\alpha,\beta$ are non-proportional roots then
$\langle\alpha,\ch\beta\rangle=0,\pm1$. We adopt the convention that in this
case all roots are long.

We say a root system is {\em oddly laced} if whenever $\alpha,\beta$ are
non-proportional roots then $\langle\alpha,\ch\beta\rangle=0$ or is {\it odd}.
Thus {\it oddly laced} is shorthand for {\it each simple factor is simply laced
or of type $G_2$}.
We also adopt the convention that in type $G_2$ all roots are
long.
The reason for these conventions is Lemma \ref{adams1}.
The main results in this paper hold for oddly
laced groups, although (with the general case in mind) we will only
make this assumption when necessary.

Suppose $G(\C)$ is a connected, reductive complex Lie group.
Let $\Gad(\C)$ be the adjoint group.
If $G(\C) $ is defined over $\R$ let $G$ be its real points, and let $\Gad$
be the real points of $\Gad(\C)$.
Write $\int(g)$ for the action of $\Gad(\C)$ on $G(\C)$, or
$G_{ad}$ on $G$.

Let $G_d(\C)$ be the derived group of $G(\C)$, and let $G_d$ be the
derived group of $G$. Note that by \eqref{e:derived} $G_d$ is the
identity component of the real points of $G_d(\C)$.

We denote real Lie algebras by Gothic letters $\h,\g,\t,\dots$, and their
complexifications by $\h(\C),\g(\C),\t(\C),\dots$.
We write $\sigma$ for the action of the non-trivial element of the
Galois group on $\g(\C)$ and $G(\C)$, so $\g=\g(\C)^\sigma$ and
$G(\R)=G(\C)^\sigma$. 

Fix a Cartan involution $\theta $ of $G(\C )$, i.e. $K =
G^{\theta}$ is a maximal compact subgroup of $G$.
Unless otherwise noted all Cartan subgroups of $G$ or $G(\C)$ are
assumed to be $\theta$-stable.
Let $H$ be a ($\theta$-stable) Cartan subgroup of $G$, with
complexification  $H(\C)$.
Write $\h=\t\oplus\a$ as usual, and $H=TA$ with $T=H\cap K$ and
$A=\exp(\a)$. Let $\Phi=\Phi(G,H)$ be the root system of $H(\C)$ in
$G(\C)$. Let $W(G,H)=\Norm_G(H)/H$; this is the {\it real} Weyl group.
It is a subgroup of the absolute Weyl group $W=W(G(\C),H(\C))$
which is isomorphic to the Weyl
group of $\Phi$. The Cartan involution $\theta$ acts on $W$, and
$W(G,H)\subset W^\theta$.

Note $\sigma(\alpha)=-\theta(\alpha)$ for all $\alpha\in\Phi$.
Roots are classified as real, imaginary,
complex, or compact as in \cite{green}. 
Write $\Phi=\Phi_r\cup\Phi_i\cup\Phi_{cx}$ accordingly;
note that $\Phi_r$ and $\Phi_i$  are root systems.
We also have the decomposition of $\Phi_i=\Phi_{i,c}\cup \Phi_{i,n}$
into compact and noncompact roots; $\Phi_{i,c}$ is a root system.
If $\Phi^+$ is a set of positive roots, write $\Phi^+_r=\Phi^+\cap
\Phi_r$, and $\Phi^+_i,\Phi^+_{cx}$ similarly.
Let $\rho=\rho(\Phi^+)=\frac12\sum_{\alpha\in\Phi^+}\alpha$ as
usual, and define $\rho_r,\rho_i$ and $\rho_{cx}$ similarly, so
$\rho=\rho_r+\rho_i+\rho_{cx}$.

We say $G(\C)$ is {\it acceptable} if $\rho$ exponentiates to a
character of $H(\C)$.

An important role is played by certain sets of positive roots:

\begin{definition}
\label{d:special}
A set $\Phi^+$ of positive roots is said to be special if
$\sigma(\Phi^+_{cx})=\Phi^+_{cx}$, or equivalently
$\sigma(\alpha)>0$ for all positive non-imaginary roots.
\end{definition}

Let 
\begin{equation}
\label{defzah}
\Gamma(H)=\exp(i\a)\cap H=\{\exp(iX)\,|\, X\in \a, \exp(2iX)=1\}.
\end{equation}

An important role is played by a certain character of $\Gamma(H)$. 
Let $S$ be a set of complex roots such that
the set of all complex roots is $\{\pm\alpha,\pm\sigma\alpha\,|\,
\alpha\in S\}$. Define 
\begin{equation}
\label{e:zetaxc}
\zetacx(G,H)(h)=\prod_Se^\alpha(h)\quad(h\in \Gamma(H)).
\end{equation}
If $G,H$ are understood we let $\zetacx=\zetacx(G,H)$.
It is elementary to see that $\zetacx$ is independent of the choice of
$S$, factors to $G_{ad}$, 
and for all $h\in \Gamma(H)$ satisfies
\begin{subequations}
\renewcommand{\theequation}{\theparentequation)(\alph{equation}}  
\begin{align}
\zetacx(h)&=\zetacx(wh)\quad (w\in W(G,H))
\label{e:zetacx1}\\
\zetacx(h)&=e^{\rho_{cx}}(h)=e^{\rho-\rho_r}(h)
\label{e:zetacx2}
\end{align}
\end{subequations}
for any special set of positive roots.   
For (b) note that if $\Phi^+$ is a special set of positive roots and $S \subset \Phi^+_{cx}$ is as above, then
for $X\in \mathfrak a(\C)$
\begin{equation}
\rho_{cx}(X)=(\rho-\rho_r)(X)=\sum_{\alpha\in S}\alpha (X),
\end{equation}
and it follows that
$e^{\rho_{cx}}(\exp(X))=e^{\rho-\rho_r}(\exp(X))=e^{\rho_{cx}(X)}$ is a well-defined character of $A(\C)$.

For $\alpha\in\Phi_r$ let
$m_\alpha=\ch\alpha(-1)=\exp(\pi i\ch\alpha)$ and define

\begin{equation}
\label{e:GammaR}
\Gamma_r(H)=\langle m_\alpha\,|\,
\alpha\in\Phi_r\rangle\subset \Gamma(H)\cap G_d^0.
\end{equation}
It is well known that 
\begin{equation}
\label{e:Gamma(H)}
H=\Gamma(H)H^0,\quad H\cap G^0=\Gamma_r(H)H^0.
\end{equation}
It is also well known that if $H_s$ is a maximally
split Cartan subgroup of $G$ then 
\begin{subequations}
\renewcommand{\theequation}{\theparentequation)(\alph{equation}}  
\begin{equation}
\label{e:HsG0}
G=H_sG^0,
\end{equation}
and this implies
\begin{equation}
\label{e:derived}
G_d=G_d^0.
\end{equation}
\end{subequations}

\sec{Admissible Triples}
\label{s:admissible}

Fix $G$ as in Section \ref{s:notation} and suppose 
$p:\tG\rightarrow G$ is a two-fold cover.
We identify the kernel of $p$ with $\pm1$.

If $H$ is a subgroup of $G$ we let $\tH=p\inv(H)$. 
Let $Z(\tH)$ be the center of $\tH$ and let
$Z_0(H)=p(Z(\tH))$. It is immediate that $Z(\tH)=p\inv(Z_0(H))$,
so to
describe $Z(\tH)$ it is enough to describe $Z_0(H)$. 
In particular $Z_0(G)\subset Z(G)$ plays an important role; it is
immediate that $Z_0(G)=Z(G)$ if $G$ is connected (cf.~\eqref{e:zconnected}).

Suppose $H$ is a Cartan subgroup of $G$. 
Typically $\tH$ is not abelian and
$Z_0(H)$ plays an important role. It is easy to see
\begin{equation}
H^0\subset Z_0(H)\subset H.
\end{equation}

Fix a real or  imaginary root $\alpha$. 
Associated to $\alpha$ is the root subgroup $M_\alpha$, which is
locally isomorphic to $SL(2,\R)$ or $SU(2)$. 
As in \cite[Definition 3.2]{dualityonerootlength} we say $\alpha$ is {\it
metaplectic} if $p\inv(M_\alpha)$ is a nonlinear group.
If $M_\alpha$ is compact it has no such cover, nor does 
$SL(2,\R)/\pm 1$, so if $\alpha$ is metaplectic then
$M_\alpha\simeq SL(2,\R)$ and  $\alpha$ is either real or noncompact
imaginary.

Let $\wt m_\alpha$ be an
inverse image of $m_\alpha$ in $\tG$. It is easy to see that 
$\wt m_\alpha$ has order $4$ if $\alpha$ is metaplectic, and $1$ or $2$
otherwise.
For the next Lemma see 
\cite{groupswithcovers} or \cite[Lemma 3.3]{dualityonerootlength}.
 
\begin{lemma}
\label{adams1}
Assume $G(\C )$ is simple and simply connected and that $\tG$ is
nonlinear. Fix a Cartan subgroup $H$ of $G$.
Then a real or imaginary root $\alpha$ of $H(\C)$ in $G(\C)$ is metaplectic
if and only if $\alpha$ is long and is real or noncompact imaginary.
Furthermore $G$ admits a nonlinear
cover if and only if there is a Cartan subgroup $H$ with a long real
or long noncompact imaginary root.
If this condition holds the nonlinear two-fold cover is
unique up to isomorphism.

It is enough to check this condition on a fundamental or 
maximally split Cartan subgroup.
If $G$ is oddly laced   it is enough to check this
condition on any Cartan subgroup.
\end{lemma}

\begin{definition}
\label{defadm} 
We say that $p:\tG \rightarrow G$ is an
admissible two-fold cover if for every  $\theta$-stable Cartan subgroup $H$,
every long real or long noncompact imaginary root is metaplectic.

Equivalently $\tG$ is admissible if and only if $\wt{G_i}$ is
nonlinear for every simple factor $G_i$ of $G$ which admits such a cover.
\end{definition}

See \cite[Definition 3.4]{dualityonerootlength}).

\begin{example}
\label{ex:U11_a}
Let $G=U(1,1)$. 
This has three non-trivial two-fold covers, described by their
restriction to $T\simeq U(1)$, the diagonal 
maximal compact subgroup.
The cover of $T$ is a connected torus $\wt T$, and
is best described by its character lattice. Write the character
lattice $X^*(T)$ as $\Z^2$ in the usual coordinates. The three covers
are given by $X^*(\wt
T)=\frac12\Z\oplus\Z$, $\Z\oplus\frac12\Z$, or
$\Z^2\cup(\Z+\frac12)^2$. 

The derived group is $SU(1,1)\simeq SL(2,\R)$, which has a unique
non-trivial two-fold cover. The first two covers of $U(1,1)$ restrict
non-trivially to $SU(1,1)$, and are therefore admissible. The third
one is trivial when restricted to $SU(1,1)$ (it is the $\sqrt{\det}$
cover), and is not admissible.
\end{example}

\medskip

As discussed in the Introduction we are going to lift characters from
a real form of a quotient of $G(\C)$. We need to impose some
conditions on this quotient. This will take up the remainder of this section.
 
Suppose $C$ is a finite subgroup of $Z(G)$.
Then $\oG(\C)=G(\C)/C$ is defined over $\R$, and let
$\oG$ be its real points. Write $\op :G( \C ) \rightarrow
\oG( \C)$ for the projection map.
Let $\Orb(G)$ be the conjugacy classes of $G$.
If $g,g'\in G$ are conjugate by $G(\C)$ we say $g,g'$ are stably
conjugate, and let $\Orbst(G)$ be the set of stable conjugacy classes.
(This is a naive definition, which agrees with the usual one for
strongly regular semisimple elements.)
Similar
notation applies to other groups.

For $g\in G$ let $\O(G,g)$ be the conjugacy class of $g$, and
$\O^{st}(G,g)=\{xgx\inv\,|\, x\in G(\C), xgx\inv\in G\}$ the stable
conjugacy class. Similar notation applies to other groups.

\begin{definition}
\label{d:phi}
Assume $C$ is a two-group. 
For $h\in \oG$ let $\phi(h)=s(h)^2$ where $s:\oG(\C)\rightarrow G(\C)$ is any section.
\end{definition}

\begin{lemma}
\label{l:phi}
The map $\phi$ is well defined, and satisfies:
\begin{enumerate}
\item  $\phi(\oG)\subset G$ and $\phi(\oG^0)\subset G^0$,
\item $\phi$ induces a map $\Orb(\oG^0)\rightarrow \Orb(G^0)$,
\item $\phi$ induces a map $\Orbst(\oG)\rightarrow \Orbst(G)$.
\end{enumerate}
\end{lemma}

\begin{proof}
Since $C$ is a two-group it is immediate that $\phi(h)$ is independent
of the choice of $s$.

Suppose $h\in\oG, g\in G(\C)$, and $\op(g)=h$. Then
$\op(\sigma(g))=\sigma(\op(g))=\op(g)$, so $\sigma(g)=zg$ for some $z\in
C$. Since $C$ is a two-group $\sigma(g)^2=g^2$, so $\phi(h) = g^2\in G$.
Furthermore since $\op:G^0\rightarrow \oG^0$ is surjective the second
assertion in (1) is clear.

Define
\begin{equation}
\label{e:phi}
\phi(\O(\oG^0,g))=\O(G^0,\phi(g))\quad (g\in\oG^0).
\end{equation}
If $x\in\oG^0$ choose $y\in G^0$ with $\op(y)=x$. 
Then $\phi(xgx\inv)=y\phi(g)y\inv$ so this is well defined, proving (2).
Similarly define 
\begin{equation}
\phi(\Ost(\oG,g))=\Ost(G,\phi(g))\quad (g\in \oG).
\end{equation}
Suppose $x\in\oG(\C)$ and $xgx\inv\in\oG$.
Choose $y\in G(\C)$ with $\op(y)=x$. Then
$\phi(xgx\inv)=y\phi(g)y\inv$, so again $\phi$ is well-defined.
\end{proof}

It is worth noting that $\phi$ does {\it not} define a map
$\Orb(\oG)\rightarrow \Orb(G)$. Suppose we try to define
$\phi(\O(\oG,g))=\O(G,\phi(g))$ as in \eqref{e:phi}.
For this to be well defined we need to know that for $x\in\oG$,
$\phi(xgx\inv)$ is $G$-conjugate to $\phi(g)$. This is true if $x\in
\op(G)$, but not in general. (By (3) $\phi(g)$ is {\it stably} conjugate to
$\phi(xgx\inv)$).

\begin{example}
\label{e:sl2orbit}
Let $G=SL(2,\R)$ and $\oG=PSL(2,\R)\simeq SO(2,1)$  (cf.~Example
\ref{ex:pgl2}).
Let $t(\theta)=\diag(t'(\theta),1)\in SO(2,1)$ (with the appropriate
form) where $t'(\theta)=\begin{pmatrix}
\cos(\theta)&\sin(\theta)\\
-\sin(\theta)&\cos(\theta)
\end{pmatrix}\in SL(2,\R)$.
Then just as in Example \ref{ex:pgl2}
$\phi(t(\theta))=t'(\theta)$.
However while $t(\theta)$ is conjugate to $t(-\theta)$, $t'(\theta)$
is not conjugate to $t'(-\theta)$ (for generic $\theta$). 
Therefore 
the map
$\phi(\O(SO(2,1),t(\theta))=\O(SL(2,\R),t'(\theta))$ is not well defined.
\end{example}

We assume throughout this paper that $C$ is a two-group (and therefore
finite).

Now fix a Cartan subgroup $H$ of $G$, and let $\oH(\C)=H(\C)/C\subset
\oG(\C)$, with real points $\oH$.  We define $\phi:\oH(\C)\rightarrow
H(\C)$ by the same formula as in Definition \ref{d:phi}.
Equivalently, if we write $\exp:\h(\C)\rightarrow H(\C)$,
$\overline{\exp}:\h(\C)\rightarrow \oH(\C)$, then for $X\in \h(\C)$,
$\phi(\overline{\exp}(X))=\exp(2X)$.
It is immediate that
\begin{equation}
\label{e:2alpha}
\alpha(\phi(h))=\alpha(h)^2\quad(\text{for all }\alpha\in\Phi).
\end{equation}

\begin{lemma}
\label{bs4} 
The homomorphism $\phi :\oH( \C) \rightarrow
H(\C ) $ has the following properties.
\begin{enumerate}
\item $\phi(wh)=w\phi(h)$ for all $w\in W(\Phi),h\in \oH(\C)$,
\item $\phi (\oH^0) = H^0$  ,
\item $\phi (h)=1$ for all $h \in \Gamma_r(\oH)$,
\item  $\phi (\oH \cap \oG ^0) = H^0$ ,
\item $ \phi (\Gamma (\oH )) = \Gamma ( H ) \cap C$,
\item $\phi (\oH ) = (\Gamma (H ) \cap C)H^0$.
\end{enumerate}
\end{lemma}

\begin{proof} The first two are clear from the definition.  For 
(3) note that
$\phi(\overline\exp(\pi i\ch\alpha))=\exp(2\pi i\ch\alpha)=1$ 
for all $\alpha \in \Phi _r$.  Then (4) follows since $\oH
\cap \oG ^0 = \Gamma_r(\oH)\oH^0$ (cf. \ref{e:Gamma(H)}).

Let  $h=\overline \exp (iX) \in \Gamma (\oH )$ where $X \in \a $ with
$\overline \exp (2iX ) = 1$.  Now $\phi(h) = \exp (2iX) \in \exp (i\a
)$.  But $ \op (\phi (h) ) = \overline \exp (2iX ) = 1$ so that $\phi
(h) \in C \cap \exp (i\a ) \subset \Gamma (H ) \cap C$.  Conversely,
let $z _G = \exp (iY) \in C \cap \Gamma (H )$ where $Y \in \a $ with
$\exp (2iY)=1$.  Since $z_G \in C$, $\op (z_G) = \overline \exp
(iY)=1$. Let $z = \overline \exp (iY/2)$. Then $z\in \Gamma (\oH ) $ and
$\phi (z) = z_G$. Finally (6) follows since $\oH = \Gamma (\oH ) \oH ^0$.
\end{proof}

It is a standard fact that the character of an irreducible genuine
representation of $\tG$, considered as a function on the regular
semisimple elements, vanishes off of $Z(\tH)$ (Lemma \ref{l:relevant}).
Our character
identities involve the image of $\phi$ in $H$, so we would like to
know that $\phi(\oH)\subset Z_0(H)$, with image as large as possible. 

Recall $C\subset Z(G)$. We first show that we may as well assume
$C\subset Z_0(G)$.

Let $H_s$ be a maximally split Cartan subgroup of $G$, and assume
$\phi(\oH_s)\subset Z_0(H_s)$.
Let $C_s=\Gamma(H_s)\cap C$ so 
that by Lemma \ref{bs4} 
$C_sH_s^0=\phi ( \oH_s)\subset Z_0(H_s)$.
Therefore
$\wt{C_s}$ centralizes $\wt{H_s}$, and 
by \eqref{e:zgg0} $\wt{C_s}$ centralizes $\wt{G^0}$. Therefore $\wt{C_s}$ centralizes
$\wt{H_s}\wt{G^0}$ which equals $\tG$ by \eqref{e:HsG0}.
Therefore $C_s \subset Z_0(G)$.

Now for a general Cartan subgroup  $H=TA$,
there is a maximally split Cartan $H_s=T_sA_s$ so 
that $A\subset A_s$. Then $\Gamma (H) \subset \Gamma (H_s) $, so
that $\phi (\oH) = (C \cap \Gamma (H)) H^0 = (C_s \cap \Gamma (H)
) H^0$.  That is, if we use $C_s$ in place of $C$ we have the same
images $\phi(\oH)$ for every Cartan subgroup $H$.  Thus we may
as well assume that $C\subset Z_0(G)$.

\begin{lemma}
\label{centerh} 
Assume that $C \subset Z_0(G)$.  Then for
every Cartan subgroup $\phi (\Gamma (\oH))
\subset Z_0(G)$ and $\phi ( \oH) \subset Z_0(H)$.
Furthermore $\phi (Z( \oG)) \subset Z_0(G)$. 
\end{lemma}

\begin{proof} The first part follows from the above discussion.  Let
$H_s$ be a maximally split Cartan subgroup with roots $\Phi $.  Let $z \in Z(\oG)
\subset \oH_s, \tilde z \in Z(\tH_s)$ such that $p(\tilde z) = \phi
(z)$.  Then for every $\alpha \in \Phi$, 
$$
\alpha(\tilde z) =\alpha(\phi (z)) = \alpha(z)^2 = 1
$$ 
by \eqref{e:2alpha}.
Thus $\tilde z \in\Cent_{\tG}(\tG^0) \cap Z(\tH_s) =
Z(\tG)$ as above.  
\end{proof}

To reiterate, so far we have assumed $C$ is a (finite) subgroup  
of $ Z_0(G)$, and $z^2=1$ for all $z\in C$.
For the definition of transfer factors we need to impose a further technical
condition on $C$. 

Suppose $\tilde\chi$ is a genuine character of
$Z(\tG)$.  Then $\tilde\chi^2$ factors to a character of $Z_0(G)$. 
Suppose $z\in Z_0(G)$ has order $2$ and $p(\tilde z)=z$.
Then $\tilde z$ has order $2$ or $4$ and $\tilde\chi^2(z)=1$ or
$-1$ accordingly, independent of $\tilde\chi$.     Assume that
$G_d( \C)$ is acceptable.

Let
\begin{equation}
\label{e:zeta2}
\zeta_2(z)=\tilde\chi^2(z)e^\rho(z)\quad z\in Z_0(G) \cap G_d, z^2=1.
\end{equation}

This is well defined since $G_d$ is acceptable, and 
is independent of the choices of $\wt\chi$, $H$ and $\Phi^+$.

\begin{definition}
\label{assumptions} 
An admissible triple is a set $(\tG,G, \oG)$ where:
\begin{enumerate}
\item $G$ is the set of real points of
a connected reductive complex Lie group $G(\C)$ such that the derived
group $G_d( \C)$ is acceptable, with oddly laced root system
i.e. each simple
factor is simply laced or $G_2$  
(cf.~Section \ref{s:notation}).
\item  $p:\tG \rightarrow G$ is an admissible two-fold
cover (Definition \ref{defadm}).
\item $\oG $ is the set of real points of $\oG (\C)= G(\C)/C$ where $C$ is a finite
subgroup of $Z(G)$ satisfying the following conditions:
\begin{enumerate}
\item  $c^2=1$ for all $c \in C$,
\item  $C \subset Z_0(G) = p(Z(\tG))$,
\item $\zeta_2(z)=1$ for all $z\in C\cap G_d$ (cf. \ref{e:zeta2}).
\end{enumerate}
\end{enumerate}
\end{definition}

For example if $\tG$ and $G$ satisfy (1) and (2) then $(\tG,G,G)$ is
an admissible triple.

More generally, suppose $\tG$ and $G$ satisfy conditions (1) and (2).
Then $(\tG,G,\oG)$ is an admissible triple for $\oG(\C)=G(\C)/C$ for 
any subgroup $C$ satisfying conditions (3)(a-c).
We may take $C=1$; by Lemma \ref{bs4}(6) taking $C$ bigger makes $\phi(\oH)\subset Z_0(H)$ 
bigger, and in general, we get
the best lifting results when we take $C$ as large as possible.
However, even if we take $C$ as large as possible for $G$, it may not be
the maximal choice for Levi subgroups of $G$. Thus we do not specify
$C$ beyond the requirements of Definition \ref{assumptions}.

Assume that $(\tG,G, \oG)$ is an admissible triple.
Let $H=TA$ be a Cartan subgroup of $G$. 
Let $M=\Cent_G(A)$ be a cuspidal Levi subgroup. Then
$\tM\subset \tG$ and $\overline M=\Cent_{\oG}(\overline A)\subset \oG$ are also cuspidal
Levi subgroups.
We will show in Section \ref{s:levi} that
$(\tM ,M ,\overline M)$ is also an admissible triple.  Thus
our constructions are 
compatible with parabolic induction.

The following lemma, which is an immediate
consequence of \eqref{e:zetacx(h)2} 
 is useful for producing admissible triples.

\begin{lemma}\label{zeta2}    
Let $H$ be a maximally split Cartan subgroup of $G$.  Then
$\zeta _2(c) = 1$ for all $c \in Z_0(G) \cap \Gamma _r(H)$.
\end{lemma}

\begin{lemma}
\label{notsogood}     
Suppose that $G(\C )$ is simple and simply connected, and let $H$ be
a maximally split Cartan subgroup of $G$.  Then we can choose $C$ so
that $\phi (\oH ) = Z_0(H)$ except when $G = SU(n,n)$ where $n$ is
even, $G = Spin(p,q)$ where $p$ and $q$ are even with $p > q \geq 2$
or when $G=Spin^*(2n)$ where $n\equiv 0$ mod 4.
\end{lemma}

\begin{proof}     

Clearly if $Z_0(H)$ is connected we can take $C = \{ 1 \}$.  We will prove in (4.4) and Proposition \ref{p:ZH}
that $Z_0(H) = Z_0(G)H^0 = Z(G)H^0$ since $G$ is connected.  But since
$H = \Gamma _r(H)H^0$ where  $\Gamma _r(H)$ is a two
group, it is easy to see that $Z(G)H^0 = Z_e(G)H^0$ where $Z_e(G)$ is
the group of elements in $Z(G)$ with even order.  Thus $Z_0(H)$ is
connected when $Z_e(G) = \{ 1 \}$.  This will be the case for all real
forms when $\Phi $ is of type $A_{2n}, E_6, E_8$, or $G_2$.  Thus we
may as well assume that $\Phi $ is type $A_{2n-1}, D_n$, or $E_7$.  In
these cases $Z_e(G) = Z(G)$.  We may also assume that $G$ is a real
form such that $Z_0(H)= Z (G)H^0$ is not connected.  In particular, we
can assume that $G$ is not compact or complex.  Suppose that $Z(G)
\subset \Gamma _r(H)$.  Then $Z(G)$ is a two-group, and by Lemma
\ref{zeta2} we can take $C = Z(G)$, giving us $ \phi (\oH) = Z(G)H^0 =
Z_0(H)$ as desired.  This will always be the case if $G$ is the split
real form.

Suppose that $\Phi$ is of type $A_{2n-1}$. Then $G$ is split or $H$ is connected
unless $G=SU(n,n)$.
In this case $Z_0(H) =
H$ has two components.  Now $Z(G) = <z>$ is cyclic of order $2n$ and
$z_0 = z^n \in \Gamma _r(H)$ is the unique
element of order two.  Thus the biggest $C$ we can take is $C = \{ 1,
z_0 \}$.  When $n$ is odd, $z_0 \not \in H^0$ so that $\phi (\oH ) =
CH^0 = H$, but when $n$ is even, $z_0 \in H^0$ so that $\phi (\oH ) =
CH^0 = H^0$ is a proper subgroup of $Z_0(H)$.
    
Suppose that $\Phi$ is of type $D_n$.  Let $G=Spin(p,q)$ with $p\ge q$. Then $G$ is
split or $H$ is connected unless $p>q\ge 2$.
If $p,q$ are odd, then $Z(G)= Z_2 \subset \Gamma
_r(H)$.  If $p,q$ are even, then $Z(G) = \{ 1, z_0, z_1, z_0z_1 \}$
where $z_0 \in \Gamma _r(H)$, but $z_1 \not \in \Gamma (H)$.  Thus the
biggest $C$ we can take is $C = \{ 1, z_0 \}$.  But $z_1 \not \in
CH^0$ so that $CH^0$ is a proper subgroup of $Z_0(H)$.

Let $G=Spin^*(2n)$. If $n$ is odd $H$ is connected, so assume $n=2k$.
Then $Z(G) = \{ 1, z_0, z_1, z_0z_1 \} $ as
above where $z_1 \in \Gamma _r(H)$, but $z_0 \not \in \Gamma (H)$ so
the biggest $C$ we can take is $C = \{ 1, z_1 \}$.  If $k$ is even,
then $z_0 \not \in CH^0$ while $z _0 \in CH^0$ when $k$ is odd.

Finally, if $\Phi = E_7$, then either $G$ is split or $Z_0(H)$ is connected.
\end{proof}

\begin{example}     

Let $H=TA$ be a maximally split Cartan subgroup of $G$.  By
Lemma \ref{notsogood}  it is not always possible to pick
$C$ satisfying the conditions of Definition \ref{assumptions} such
that $\phi (\oH ) = Z_0(H)$.  Let $M=\Cent_G(A)$.  
Here we show it is
possible to pick a finite central subgroup $C_M$ of $M$ satisfying the
conditions of Definition \ref{assumptions} for $M$ such that $\phi
(\oH) = Z_0(H)$ where $\phi $ is defined using $\oM (\C ) = M(\C
)/C_M$.

  Recall $H=\Gamma (H)H^0$, so $\tH =\wt{\Gamma(H)}\tH^0$ and 
$Z(\tH)= Z(\wt{\Gamma (H)})\tH^0$.  
Then $Z_0(H)=Z_0(\Gamma(H))H^0$, and since $\Gamma(H)$ is a two-group
there exists $C_M\subset Z_0(\Gamma(H))$ such that $C_M\cap H^0=\{1\}$ and
$Z_0(H)=C_MH^0$.
We have to show $C_M$ satisfies the conditions of Definition
\ref{assumptions}(3). Condition (a) is obvious,
and since $M=HM^0=\Gamma(H)H^0M^0$ it follows easily that
$Z_0(\Gamma(H))\subset Z_0(M)$, which is (b). Finally $C_M\cap
M_d\subset C_M\cap H^0=\{1\}$, giving (c).
\end{example}

\sec{Cartan Subgroups and Cayley Transforms}
\label{s:cartan}

We  need some structural facts about Cartan subgroups.     We first define commutators on $G$ with respect to an admissible cover $\tG$.

Assume for the moment that $G_d(\C)$ is simply connected and 
$\tau$ is an automorphism of $G$.
Then $\tau$ stabilizes
$G_d=G_d^0$ (cf. \eqref{e:derived}) and by Lemma \ref{adams1} 
it lifts uniquely to an automorphism
$\tilde\tau$ of the (unique) admissible cover  $\wt{G_d}$ of $G_d$.
Suppose $g\in G_d$ and $\tau(g)=g$.
Choose an inverse image $\tilde g$ of $g$ and set

\begin{equation}
\label{e:taug}
\{\tau,g\}=\tilde\tau(\tilde g)\tilde g\inv.
\end{equation}
Then $p(\{\tau,g\})=1$ so $\{\tau,g\}=\pm1$. It is independent of the choice of $\tilde
g$.

We now drop the assumption that $G_d(\C)$ is simply connected, but
suppose that $\tG$ is an  admissible cover of $G$.
Suppose 
$g,h$ are commuting elements of $G$.
Choose inverse images
$\tg,\wt h$ of $g,h$ in $\tG$ and 
define
\begin{equation}
\label{e:comm}
\{g,h\}=\tg\wt h\tg\inv\wt h\inv.
\end{equation}
Since $p(\tg\wt h\tg\inv\wt h\inv)=ghg\inv h\inv=1$,
$\{g,h\}=\pm1$, 
independent of the choices.
We say $\{g,h\}$ is the commutator of $g,h$
(with respect to the cover $\tG$).
It is continuous in each factor so
\begin{equation}
\label{e:zgg0}
\{Z(G),G^0\}=1.
\end{equation}
Note that this implies
\begin{equation}
\label{e:zconnected}
Z(\tG)=\wt{Z(G)}\quad\text{if $G$ is connected}.
\end{equation}
In general we only have $Z(\tG)\subset \wt{Z(G)}$,
for example if $G$ is a torus and $\tG$ is nonabelian.

Suppose
$\alpha\in \Phi(G,H)$ is a long real root.
Let $M_\alpha$ be as in Section \ref{s:admissible}.
By Lemma \ref{adams1} $\alpha$ is metaplectic, and by the discussion
at the beginning of Section \ref{s:admissible}, $M_\alpha\simeq
SL(2,\R)$.
Let  $\tau$ be an automorphism of $G_d$ satisfying
$\tau(H\cap G_d)=H\cap G_d$ and 
$\tau(\alpha)=\alpha$. 
Then $\tau(M_\alpha)=M_\alpha$, so for $g\in M_\alpha$ define $\{\tau,g\}$ by \eqref{e:taug} applied to
$M_\alpha$ .
Let $T_\alpha\simeq S^1$ be a $\tau$-stable compact Cartan subgroup of $M_\alpha$.
Then for all $t\in T_\alpha$, $\tau(t)=t$ or $t\inv$.
Recall $m_\alpha=\ch\alpha(-1)\in Z(M_\alpha)$.

For the next result we do not need to assume $G$ is oddly laced.

\begin{proposition}
\begin{subequations}
\renewcommand{\theequation}{\theparentequation)(\alph{equation}}  
Suppose $\alpha$ is a long real root. Then
\begin{equation}
\label{e:tau1}
\{\tau,m_\alpha\}=
\begin{cases}
1&\tau(t)=t\quad \text{ for all }t\in T_\alpha\\  
-1&\tau(t)=t\inv\quad \text{ for all }t\in T_\alpha\\  
\end{cases}
\end{equation}
For all $h\in H$,
\begin{equation}
\label{e:tau2}
\{h,m_\alpha\}=\sgn(\alpha(h)).
\end{equation}
If $\beta\in\Phi_r(G,H)$ then
\begin{equation}
\label{e:tau3}
\{m_\alpha,m_\beta\}=(-1)^{\langle\alpha,\ch\beta\rangle}.
\end{equation}
\end{subequations}
\end{proposition}

\begin{proof}
We can compute $\tau(m_\alpha)$ by working in $T_\alpha$. If
$\tau$ acts trivially on $T_\alpha$ then the same holds for the action
of $\tilde\tau$ on $\wt{T_\alpha}$ and $\{\tau,m_\alpha\}=1$.
Suppose $\tau(t)=t\inv$ for all $t\in T_\alpha$.
Then $\wt\tau(\wt t)=\wt t\inv$ for all $\wt t\in \wt T_\alpha$.
Let $\wt m_\alpha$ be an inverse image of $m_\alpha$ in $\wt
T_\alpha$. Then $\{\tau,m_\alpha\}=\wt m_\alpha^{-2}$.  Since $\alpha$
is metaplectic $\wt m_\alpha$ has order $4$ (see Section
\ref{s:admissible}) and  $\wt
m_\alpha^{-2}=-1$.
This gives (a).

It is a standard fact (essentially a calculation in $GL(2,\R)$) 
that for all $h\in H,t\in T_\alpha$,
$hth\inv=t^{\sgn(\alpha(h)}$.
This proves (b), and (c) follows
from this and the identity
$\alpha(m_\beta)=(-1)^{\langle \alpha,\ch\beta\rangle}$.
\end{proof}

The following result is of fundamental importance, and does not hold 
in the two root length case (for example for $Sp(2n,\R)$).

\begin{proposition}
\label{p:ZH}
Assume $G$ is oddly laced.
Then
\begin{equation}
Z(\tH)=Z(\tG)\wt{H^0}
\end{equation}
\end{proposition}

We first prove

\begin{lemma}
\label{l:ztg}
In the setting of the Proposition we have
\begin{equation}
\label{e:ztg1}
Z(\tH)\subset \wt{Z(G)}\wt{H^0}.
\end{equation}
\end{lemma}

\begin{proof}
It is enough to show if $h\in H$ and
$\{h,m_\alpha\}=1$ for all  real roots $\alpha$ then $h\in Z(G)H^0$.
(This is where the oddly laced condition appears: otherwise we only
have this identity for the long real roots.)
By \eqref{e:tau2} it is enough to show
\begin{equation}
\alpha(h)>0\text{ for all }\alpha\in\Phi_r(G,H)\text{ implies }h\in Z(G)H^0
\end{equation}
This is a straightforward calculation using roots and weights. Choose
a basis of the root lattice of the form 
$$
\alpha_1,\dots, \alpha_m,\alpha_{m+1},\dots,\alpha_{m+n},
\alpha_{m+n+1},\dots
\alpha_{m+n+2r}
$$
where $\alpha_i$ is real for $1\le i\le m$, imaginary for $m+1\le i\le
m+n$, and $\theta\alpha_{m+n+2i-1}=-\alpha_{m+n+2i}$ for all $1\le i\le r$.
It is a basic fact about lattices that such a
basis exists.
Let $\ch\lambda_1,\dots,\ch\lambda_{m+n+2r}$ be the dual basis of coweights.

For each $i$ choose $x_i\in\C$ so that $e^{x_i}=\alpha_i(h)$.
If $i\le m$ we may assume $x_i\in \R$ since $\alpha_i(h)>0$.
We may also assume $x_{m+n+2i}=\overline x_{m+n+2i-1}$ for all $1\le
i\le r$.
Let $X=\sum_i x_i\ch\lambda_i$, $h_1=\exp(X)$.
It follows easily that 
$X\in \h$, $h_1\in H^0$, and $\alpha(h)=\alpha(h_1)$ for all
roots $\alpha$. 
Let $z=hh_1\inv$; this is contained in $G$, and also $Z(G(\C))$, since 
$\alpha(z)=1$ for all roots $\alpha$. Therefore $z\in Z(G)$, and 
$h=zh_1\in Z(G)H^0$.
\end{proof}

\begin{proof}[Proof of the Proposition]
The statement is equivalent to
\begin{equation}
\label{e:Z0H}
Z_0(H)=Z_0(G)H^0.
\end{equation}
It is enough to show $Z_0(H)\subset Z_0(G)H^0$ (the reverse inclusion
is obvious). 
We first prove this assuming $H$ is
a maximally split Cartan subgroup.
Suppose $h\in Z_0(H)$.
By the Lemma $h=zy$ with $z\in Z(G), y\in H^0$.
Then $y\in Z_0(H)$, so $z=hy\inv\in Z_0(H)$. 
It is enough to show $z\in Z_0(G)$. We have $\{z,g\}=1$ for all $g\in H$,
and also for $g\in G^0$ by 
\eqref{e:zgg0}.
Since $H$ is maximally split $G=HG^0$, so $z\in Z_0(G)$. 

The general case will be proved  after Lemma
\ref{l:tHaHb}, using Cayley transforms.
\end{proof}

Suppose $H$ is a Cartan subgroup and $\alpha$ is a real or noncompact
imaginary root.  We define the Cayley transform $H_\alpha$ of
$H$ with respect to $\alpha$ as in 
\cite[\S 11.15]{overview}; also see \cite[Proposition 8.3.4 and
8.3.8]{green}.
  Then $H_\alpha$ has a noncompact imaginary or real root
$\beta$, and the Cayley transform of $H_\alpha$ with respect to
$\beta$ is $H$.

In order to emphasize the symmetry of the situation we change notation
and let $H_\alpha=T_\alpha A_\alpha$ be a Cartan subgroup with a real
root $\alpha$. Then we let $H_\beta=T_\beta A_\beta$ be its Cayley
transform, with noncompact imaginary root $\beta$. 

Define $Z_\alpha\in\mathfrak p$ and $Z_\beta\in \k$
as in \cite[Proposition 8.3.4 and 8.3.8]{green}, and let
$B_\alpha=\exp(\R Z_\alpha)\simeq \R^+\subset H_\alpha$, and 
$B_\beta=\exp(\R Z_\beta)\simeq S^1\subset H_\beta$. 
It follows easily from  \cite[8.3.4 and 8.3.13]{green} that we have
\begin{subequations}
\renewcommand{\theequation}{\theparentequation)(\alph{equation}}  
\label{e:HaHb}
\begin{align}
\label{e:HaHb1}
(H_\alpha\cap H_\beta)B_\alpha&\subset
H_\alpha\quad\text{index $1$ or $2$}\\
\label{e:HaHb2}
(H_\alpha\cap H_\beta)B_\beta&=H_\beta.
\end{align}
\end{subequations}
The root $\alpha$ takes positive real values on $(H_\alpha\cap
H_\beta)B_\alpha$.
We say $\alpha$ is type $I$ if  inclusion (a) is
an equality.
Otherwise $\alpha$ is of type $II$, in which case $H_\alpha$ is
generated by the left hand side and an element $t$ satisfying
$\alpha(t)=-1$. 

We say $\beta$ is type $I$ if $s_\beta\not\in W(G,H_\beta)$, and type
$II$ otherwise. Then $\alpha,\beta$ are both of type $I$ or both of type $II$.

There is an element $g\in G_{ad}(\C)$ so that if $c=\Ad(g)$,
$c^*=\Ad^*(g)$ then
\begin{subequations}
\label{e:c}
\renewcommand{\theequation}{\theparentequation)(\alph{equation}}  
\begin{align}
&c\text{ centralizes } \h_\alpha(\C)\cap \h_\beta(\C)\\
&c(\h_\beta(\C))=\h_\alpha(\C)\\
&c^*(\alpha)=\beta.
\end{align}
Note that any $c$ satisfying (a) and (b) satisfies
$c^*(\alpha)=\pm\beta$. 
\end{subequations}

\begin{definition}
\label{d:cayley}
Suppose $\chi_\alpha$ is a character of $H_\alpha$, and
$\chi_\beta$ is a character of $H_\beta$.
We say $\chi_\alpha$ is a Cayley transform of  $\chi_\beta$, and 
vice-versa, if
\begin{subequations}
\renewcommand{\theequation}{\theparentequation)(\alph{equation}}  
\begin{align}
\chi_\alpha(h)&=\chi_\beta(h)\quad (h\in H_\alpha\cap H_\beta)
\label{e:defcayley1}\\
\Ad^*(c)(d\chi_\alpha)&=d\chi_\beta.
\end{align}
\end{subequations}
where $c\in G_{ad}(\C)$ satisfies \eqref{e:c}(a-c).
\end{definition}

Here is a convenient alternative characterization of Cayley
transforms, which does not refer to the element $c$.
The proof is elementary.

\begin{lemma}
\label{l:noc}
In the setting of the Definition,
$\chi_\alpha$ is a Cayley transform of  $\chi_\beta$, and vice versa, 
if and only if  \eqref{e:defcayley1} holds, and 
\begin{equation}
\langle d\chi_\alpha,\ch\alpha\rangle= \langle d\chi_\beta,\ch\beta\rangle.
\end{equation}
\end{lemma}

\begin{remark}
In the definition of Cayley transform we are allowed to replace
$\alpha$ and $\beta$ with their negatives. This does not affect the
Cartan subgroups. Suppose $\chi_\beta$ is the Cayley transform of
$\chi_\alpha$ with respect to $\{\alpha,\beta\}$. Then the Cayley
transform of $\chi_\alpha$ with respect to $\{\alpha,-\beta\}$ is
$s_\beta(\chi_\beta)$, as is clear from the preceding Lemma. Similar
commments apply to $\alpha$.
\end{remark}

\begin{lemma}
\label{l:cayley}
Fix $\alpha$ and $\beta$.
\hfil\break
\noindent (1) Fix a character $\chi_\alpha$ of $H_\alpha$. 
There is a Cayley transform $\chi_\beta$ of $\chi_\alpha$ if and only if 
\begin{equation}
\langle d\chi_\alpha,\ch\alpha\rangle\in\Z
\end{equation}
and
\begin{equation}
\label{e:chi}
\chi_\alpha(m_\alpha)=(-1)^{\langle d\chi_\alpha,\ch\alpha\rangle}.
\end{equation}
Assume these hold. Then the
 Cayley transform of $\chi_\alpha $ is given by 
\eqref{e:defcayley1} and 
\begin{equation}
\chi_\beta(\exp(xZ_\beta))=e^{ i\langle d\chi,\ch\alpha\rangle x}.
\end{equation}

\noindent (2)  Fix a character $\chi_\beta$ of $H_\beta$. 
There are one or two choices of Cayley transform $\chi _\alpha $ 
of $\chi_\beta $, given as follows.
Define $\chi_\alpha$
restricted to $(H_\alpha\cap H_\beta)B_\alpha$ 
by 
\eqref{e:defcayley1}, and 
\begin{equation}
\chi_\alpha(\exp(xZ_\alpha))
=
e^{\langle d\chi_\beta,\ch\beta\rangle x}.
\end{equation}
If $\beta$ is of type $I$ this defines the  character $\chi_\alpha$ of $H_\alpha$.
If $\beta$ is of type II define $\chi_\alpha^{\pm}$ to be the  two
extensions of $\chi_\alpha$ to $H_\alpha$.
\end{lemma}

\begin{proof}[Sketch of proof]
This is elementary from the identities
\cite[8.3.4 and 8.3.13]{green}.
See  \cite[Lemma 8.3.7 and 8.3.15]{green},
where 
the setting is a regular character 
and the construction differs from this one by a
$\rho$-shift. For this reason the
proof of the Lemma is in fact much easier
than those in \cite{green}.
\end{proof}

We now consider Cayley transforms for $\tG$ in the oddly laced case.
The analogue of \eqref{e:HaHb} is:

\begin{lemma}
\label{l:tHaHb}
Assume $G$ is oddly laced and  $\tG$ is an admissible cover of
$G$. Then
\begin{subequations}
\renewcommand{\theequation}{\theparentequation)(\alph{equation}}  
\begin{align}
\label{e:zha}
(Z(\tH_\alpha)\cap Z(\tH_\beta))\wt{B_\beta}&=Z(\tH_\beta)\\
\label{e:zhb}
(Z(\tH_\alpha)\cap Z(\tH_\beta))\wt{B_\alpha}&=Z(\tH_\alpha)
\end{align}
\end{subequations}
\end{lemma}
  
\begin{proof}
We first prove (a).
We have to show
\begin{equation}
(Z_0(H_\alpha)\cap Z_0(H_\beta))B_\beta=
Z_0(H_\beta).
\end{equation}
Since $B_\beta$ is connected,
$B_\beta \subset H_{\beta }^0\subset Z_0(H_\beta)$, and the inclusion $\subset$ is clear.
For the opposite inclusion suppose $g\in Z_0(H_\beta)$. 
By \eqref{e:HaHb} write $g=hb$ with $h\in H_\alpha\cap H_\beta,
b\in B_\beta$. Then $h=gb\inv\in Z_0(H_\beta)$ since both $g$ and $b$
are in $Z_0(H_\beta)$. It is enough to show $h\in Z_0(H_\alpha)$.

Note that
$\{h,x\}=1$ for $x\in H_\alpha\cap H_\beta$ (since $h\in
Z_0(H_\beta)$).
Since $B_\alpha$ is connected $B_\alpha \subset H_\alpha ^0\subset
Z_0(H_\alpha)$, and therefore $\{h,x\}=1$ for $x\in B_\alpha$.
Therefore $\{h,x\}=1$ for $x\in (H_\alpha\cap H_\beta)B_\alpha$. 
If $\alpha$ is type $I$ the right hand side is $H_\alpha$ and we are
done.  Otherwise choose $t$ satisfying $\alpha(t)=-1$, so 
that $H_\alpha$ is generated by $(H_\alpha\cap H_\beta)B_\alpha$ and
$t$.

It is enough to show $\{h,t\}=1$. 
If this is not the case replace $h$ with $hm_\alpha$ and $b$ with
$bm_\alpha$ (note that $m_\alpha\in H_\alpha\cap H_\beta$ and $B_\beta$).
By \eqref{e:tau2} $\{m_\alpha,t\}=-1$, so $\{hm_\alpha,t\}=1$. 

For (b) the inclusion $\subset$ is immediate since $B_\alpha$ is
connected.
On the other hand suppose $g\in H_\alpha$ but $g\not\in (H_\alpha\cap
H_\beta)B_\alpha$. Then $\alpha(g)<0$, and by 
\eqref{e:tau2} $\{m_\alpha,g\}=-1$, so $g\not\in Z_0(H_\alpha)$.
This proves the reverse inclusion.
\end{proof}

We may now complete the proof of Proposition \ref{p:ZH}.

\begin{proof}
We have already shown this for the most split Cartan subgroup. 
By repeated use of the Cayley transform it is enough to show that in
the previous setting
\begin{equation}
Z(\tH_\alpha)=Z(\tG)\tH_\alpha^0
\Rightarrow
Z(\tH_\beta)=Z(\tG)\tH_\beta^0
\end{equation}
By the Lemma we have
\begin{equation}
\begin{aligned}
Z(\tH_\beta)&=(Z(\tH_\alpha)\cap Z(\tH_\beta))\wt{B_\beta}\\
&=
(Z(\tG)\tH_\alpha^0\cap Z(\tH_\beta))\wt{B_\beta}\\
&=
Z(\tG)(\tH_\alpha^0\cap Z(\tH_\beta))\wt{B_\beta}\\
&=
Z(\tG)\tH_\beta^0
\end{aligned}
\end{equation}
For the last equality we have used that
$(H^0_\alpha\cap H_\beta)B_\beta\subset H_\beta^0$,
which implies
$(\tH_\alpha^0\cap Z(\tH_\beta))\wt{B_\beta}\subset \tH_\beta^0$.
\end{proof}

We  now define Cayley transforms for genuine characters 
as in Definition \ref{d:cayley}, using $Z(\tH_\alpha), Z(\tH_\beta)$,
and Lemma  \ref{l:tHaHb} in place of \eqref{e:HaHb}.

\begin{definition}
\label{d:cayley2}
Assume $G$ is oddly laced and  $\tG$ is an admissible cover of $G$.
Fix $\alpha,\beta,\wt H_\alpha$ and $\wt H_\beta$.
Suppose $\wt\chi_\alpha$ is a genuine character of $Z(\tH_\alpha)$, and
$\wt\chi_\beta$ is a genuine character of $Z(\tH_\beta)$.
We say $\wt\chi_\alpha$ is a Cayley transform of  $\wt\chi_\beta$, and 
vice-versa, if
\begin{subequations}
\renewcommand{\theequation}{\theparentequation)(\alph{equation}}  
\begin{align}
\wt\chi_\alpha(h)&=\wt\chi_\beta(h)\quad (h\in Z(\tH_\alpha)\cap Z(\tH_\beta))
\label{e:defcayley2}\\
\Ad^*(c)(d\wt\chi_\alpha)&=d\wt\chi_\beta.
\end{align}
\end{subequations}
where $c\in G_{ad}(\C)$ satisfies \eqref{e:c}(a-c).
\end{definition}  

Unlike in the linear case this defines $\wt\chi_\beta$ uniquely, 
since 
\eqref{e:zha} is an equality, rather than
the containment of \ref{e:HaHb}(a)).

The analogue of Lemma \ref{l:noc} follows easily from Proposition  \ref{p:ZH}.

\begin{lemma}
In the setting of the Definition $\wt\chi_\alpha$ and
$\wt\chi_\beta$ are each others Cayley transforms if and only if
\begin{subequations}
\renewcommand{\theequation}{\theparentequation)(\alph{equation}}  
\begin{align}
\label{e:cayleynewA}
d\wt\chi_\alpha(X)&=d\wt\chi_\beta(X)\quad (X\in\h_\alpha\cap\h_\beta),\\
\label{e:cayleynewa}
\wt\chi_\alpha(z)&=\wt\chi_\beta(z)\quad(z\in Z(\tG))\\
\label{e:cayleynewb}
\langle d\wt\chi_\alpha,\ch\alpha\rangle&= \langle d\wt\chi_\beta,\ch\beta\rangle.
\end{align}
\end{subequations}
\end{lemma}

The analogue of Lemma \ref{l:cayley} is:

\begin{lemma}
\label{l:cayleynonlinear}
Suppose we are in the setting of Definition \ref{d:cayley2}.
\hfil\break
\noindent (1) Suppose $\wt\chi_\alpha$ is a genuine character of
$Z(\tH_\alpha)$.    If $\alpha$ is type II, then
the Cayley transform $\wt\chi_\beta$  of $\wt\chi_\alpha$ exists if
and only if 
\begin{subequations}
\renewcommand{\theequation}{\theparentequation)(\alph{equation}}  
\begin{equation}
\label{e:dtildechi}
\langle d\wt\chi_\alpha,\ch\alpha\rangle\in \Z+\frac12.
\end{equation}
Assume this holds.
Then $\wt\chi_\beta$ is given by \eqref{e:defcayley2} and 
\begin{equation}
\label{e:chitilde1}
\wt\chi_\beta(\tilde{\exp}(xZ_\beta))=e^{\langle d\wt\chi_\alpha,\ch\alpha\rangle ix}.
\end{equation}
If $\alpha$ is type $I$
let $\wt m_\beta=\wt{\exp}(\pi Z_\beta)\in Z(\tH_\alpha)\cap
Z(\tH_\beta)$. 
Note that $\wt\chi_\beta(\wt m_\beta)=\pm i$.
Then
the Cayley transform $\wt\chi_\beta$  of $\wt\chi_\alpha$ exists if and only if 
(\ref{e:dtildechi}) and
\begin{equation}
\label{e:chitilde2}
e^{  i\pi\langle d\wt\chi_\alpha,\ch\alpha\rangle}=\wt\chi(\wt m_\beta).
\end{equation}
In this case it is given by \eqref{e:defcayley2} and \eqref{e:chitilde1}.

\noindent (2) Suppose $\wt\chi_\beta$ is a genuine character of
$Z(\tH_\beta)$. 
Then there is a unique Cayley transform $\wt\chi_\alpha$ of
$\wt\chi_\beta$ given by 
\eqref{e:defcayley2} and 
\begin{equation}
\label{e:chitilde3}
\wt\chi_\alpha(\wt{\exp}(xZ_\alpha))=
e^{\langle d\wt\chi_\beta,\ch\beta\rangle x.}
\end{equation}
\end{subequations}
\end{lemma}

\begin{proof}
For (1) it is immediate from the definition that $\wt\chi_\beta$ is defined by
\eqref{e:defcayley2} and   \eqref{e:chitilde1}.
If $\alpha$ is type II then $Z(\tH_\alpha)\cap Z(\tH_\beta)\cap
\wt{B_\beta}=1$, and there is no further condition. If $\alpha$ is
type I then $\wt m_\beta\in Z(\tH_\alpha)\cap
Z(\tH_\beta)\cap\wt{B_\beta}$, which gives condition
\eqref{e:chitilde2}.
Case (2) is similar, and easier.
\end{proof}

\sec{Transfer Factors: Special Case}
\label{s:special}
Let $(\tG,G,\oG)$ be an admissible triple as in
Definition \ref{defadm}.
The basic idea of transfer factors is simple when $G=\oG$ is connected and semisimple
(for example take $G(\C)=\oG(\C)$  semisimple and simply
connected).
We discuss this case, before turning to the general case in the next
two sections.
We start with some general notation before specializing to the case at hand.

Suppose $\h=\t\oplus\a$ is a Cartan subalgebra of $\g$,
with corresponding Cartan subgroups $\tH,H$ and $\oH$ of $\tG,G$ and $\oG$,
respectively. Let $\tH',H',\oH'$ be the regular elements of these groups.

Let $\Phi ^+$ be a set of positive roots of $H$ in $G$.
For $h\in H'$ define
\begin{subequations}
\renewcommand{\theequation}{\theparentequation)(\alph{equation}}  
\label{e:Deltas}
\begin{align}
\Delta^0(h,\Phi^+)&=\prod_{\alpha\in\Phi^+}(1-e^{-\alpha}(h))\\
\epsilon_r(h,\Phi^+)&=sign\prod_{\alpha\in\Phi_r^+}(1-e^{-\alpha}(h))\\
\Delta^1(h,\Phi^+)&=\epsilon_r(h,\Phi^+)\Delta^0(h,\Phi^+)
\end{align}

Let  $D(h)$ be the coefficient of $t^{\text{rank G}}$ in
$\text{Det}(t+1-\Ad(h))$, and let $n$ be the number of positive roots.
It is easy to see that
\begin{align}
\label{e:easy1}
(-1)^n D(h)=\Delta^0(h)^2e^{2\rho}(h)\\
\label{e:easy2}
|D(h)|^{\frac12}=|\Delta^0(h)||e^\rho(h)|
\end{align}
where the final term is shorthand for $|e^{2\rho}(h)|^{\frac12}$.
It is evident these expressions are independent of the choice of $\Phi^+$.

If $G(\C)$ is acceptable define
\begin{equation}
\label{e:easy3}
\Delta(h,\Phi^+)=e^\rho(h)\Delta^0(\Phi^+,h),
\end{equation}
in which case
\begin{equation}
\label{e:easy4}
|D(h)|^{\frac12}=|\Delta(h,\Phi^+)|,
\end{equation}
independent of the choice of $\Phi^+$.

The same definitions apply to $\tG$ (by pulling back from $G$) and to $\oG$.
\end{subequations}

\begin{definition}
\label{d:XG}
Let
\begin{subequations}
\renewcommand{\theequation}{\theparentequation)(\alph{equation}}  
\begin{equation}
\label{e:XG}
X(\oG,\tG)=\{(h,\tg)\in \oG\times \tG\,|\, \phi(h)=p(\tg)\}.
\end{equation}
and 
\begin{equation}
\label{e:XH}
X(\oH,\tH)=X(\oG,\tG)\cap \oH\times\tH.
\end{equation}
For $\tg\in \tG$ let 
\begin{equation}
\label{e:XoGtg}
X(\oG,\tg)=\{h\in \oG\,|\, \phi(h)=p(\tg)\}
=
\{h\in \oG\,|\, (h,\tg)\in X(\oG,\tG)\}.
\end{equation}
and define $X(\oH,\tg)$ similarly.

Let $X'(\oG,\tG)$, $X'(\oH,\tH)$ be the subsets consisting of regular
semisimple elements.
\end{subequations}
With $G$ playing the role  of $\oG$ we have $\phi(h)=h^2$ and
\begin{equation}
X(H,\tH)=\{(h,\tg)\in H\times \tH\,|\, p(\tg)=h^2\}.
\end{equation}
For $(h,\tg)\in X(H,\tH)$, choose $\tilde h$ satisfying $p(\tilde h)=h$ and define
\begin{equation}
\label{e:tau}
\tau(h,\tilde g)=\tilde h^2\tg\inv.
\end{equation}
This is clearly independent of the choice of $\tilde h$, and 
$\tau(h,\tg)=\pm1$ since  $p(\tilde h^2\tg\inv)=h^2p(\tg\inv)=1$.
\end{definition}

Now suppose $G=\overline G$ is semisimple and connected. 
Recall (Definition \ref{assumptions}) we always assume $G_d(\C)$ is
acceptable. 
Since $G=G_d$,
$\Delta$ is defined, and for  $(h,\tg)\in X(H,\tH)$
it is natural to consider
\begin{equation}
\frac{\Delta(h,\Phi^+)}{\Delta(\tg,\Phi^+)}\Gamma_0(h,\tg)
\end{equation}
for some factor $\Gamma_0(h,\tg)$.
See the Introduction.
Note that the quotient is independent of the choice of $\Phi^+$. 
We want $\Gamma_0$ to be a genuine function of $\tg$,
i.e. $\Gamma_0(h,\pm\tg)=\pm\Gamma_0(h,\tg)$. 
We also want $\Gamma_0(h,\tg)$ to have small finite
order for all $h,\tg$.

Considerations of harmonic analysis make it natural to include the
$\epsilon_r$ terms (cf.~\ref{e:Deltas}(c)) and instead define
\begin{equation}
\Delta(h,\tg)=
\frac{\Delta(h,\Phi^+)\epsilon_r(h,\Phi^+)}{\Delta(g,\Phi^+)\epsilon_r(g,\Phi^+)}
\Gamma_0(h,\tg)\quad
\end{equation}
with $\Gamma_0(h,\tg)$ to be determined. Now the quotient depends on the
choice of $\Phi^+_r$.

Fix $(h,\tg)\in X(H,\tH)$.
The simplest genuine function on $X(H,\tH)$ is $\tau$ \eqref{e:tau}. 
Motivated by the case of $GL(n)$ (see the Introduction) it is natural to use this on $H^0$:
\begin{equation}
\Gamma_0(h,\tg)=\tau(h,\tg)=\tilde h^2\tg\inv=\pm1\quad (h\in
H^0).
\end{equation}

Here is another way to think of this term. Let $\tilde\chi$ be a
genuine character of $\widetilde{H^0}$. Then $\tilde\chi^2$ factors to a
character of $H^0$: for $h\in H^0$ define
$\tilde\chi^2(h)=\tilde\chi(\tilde h^2)$ where $p(\tilde h)=h$. 
Then it is easy to see that
\begin{equation}
\Gamma_0(h,\tg)=\frac{\tilde \chi(\tg)}{\tilde\chi^2(h)}\quad (h\in H^0).
\end{equation}

This suggests a way to define $\Gamma_0(h,\tg)$ for all $h$; 
choose a character $\chi_0$ of $H$ restricting to $\tilde\chi^2$ on $H^0$,
and define:
\begin{equation}
\label{e:definehtg}
\Gamma_0(h,\tg)=\frac{\tilde \chi(\tg)}{\chi_0(h)}\quad \text{for all }(h,\tg)\in X(H,\tH).
\end{equation}

Recall since $G$ is connected by \eqref{e:Gamma(H)} $H=\Gamma_r(H)H^0$. 
A basic property of $\tilde\chi$ is the following 
(cf.~\ref{e:zetacx(h)2}):
\begin{equation}
\label{e:tildechi^2}
\tilde\chi^2(t)=e^{\rho_r}(t)\quad t\in\Gamma_r(H)\cap H^0
\end{equation}
where $\rho_r=\rho(\Phi^+_r)$.
Note that $e^{\rho_r}$ only depends on the positive real roots, and
its restriction to $H^0$ is independent of all choices, since any two
differ by a sum of real roots, which is necessarily trivial on $\Gamma_r(H) \cap H^0$.

Therefore we can define $\chi_0$ by
\begin{equation}
\label{e:chi0}
\chi_0(h)=
\begin{cases}
\tilde\chi^2(h)&h\in H^0\\
e^{\rho_r}(h)&h\in \Gamma_r(H).  
\end{cases}
\end{equation}
Note that $\chi_0$ depends on the choice of $\Phi^+$,
although its restriction to $\Gamma_r(H)\cap H^0$ does not. 
Defining $\Gamma_0(h,\tg)$ by \eqref{e:definehtg}, we see:

\begin{lemma}
Fix a  set of positive roots $\Phi^+$, and let
$\rho_r=\rho(\Phi^+_r)$. 
There is a unique function $\Gamma_0$ on $X(H,\tH)$ satisfying
\begin{equation}
\label{e:Gammahspecial}
\begin{aligned}
\Gamma_0(h,\tg)&=\tau(h,\tg)\quad (h\in H^0)\\
\Gamma_0(th,\tg)&=e^{-\rho_r}(t)\Gamma_0(h,\tg)\quad (t\in \Gamma_r(H), h\in H).
\end{aligned}
\end{equation}
\end{lemma}

\begin{example}
\label{ex:sl2_special}
Let $G=SL(2,\R)$ and let $H$ be the diagonal Cartan subgroup,
which we identify with $\R^\times$.
If $p(\tilde x)=x$ and $\epsilon=\pm1$ then $(x,\tilde x^2\epsilon)\in X(H,\tH)$ and
\begin{equation}
\Gamma_0(x,\tilde x^2\epsilon)=\sgn(x)\epsilon,
\end{equation}
independent of  $\Phi^+$.
\end{example}

While this is a natural extension of $\Gamma_0(h,\tg)$ from $H^0$ to
$H$, the compelling evidence that it is the correct definition 
amounts to the matching conditions involving different Cartan
subgroups. See Lemma \ref{l:HHalpha}.

\begin{definition}
\label{d:transferspecial}
Assume $G=\oG$ is connected and semisimple.
Fix a set of positive roots $\Phi^+$. For $(h,\tg)\in X(H,\tH)$
define the transfer factor 
\begin{equation}
\label{e:transferspecial}
\Delta(h,\tg)=
\frac{\Delta(h,\Phi^+)\epsilon_r(h,\Phi^+)}{\Delta(g,\Phi^+)\epsilon_r(g,\Phi^+)}.
\Gamma_0(h,\tg)
\end{equation}
\end{definition}

While $\Gamma_0(h,\tg)$ depends on a choice of $\Phi^+$ (actually only
$\Phi^+_r$), the transfer factor itself is independent of this choice:

\begin{lemma}
$\Delta(h,\tg)$ is independent of the choice of
$\Phi^+$.
\end{lemma}

\begin{proof}
Write $h=th_0$ with $t\in \Gamma_r(H), h_0\in H^0$.
By \eqref{e:Gammahspecial}
\begin{equation}
\label{e:formulaforDelta}
\Delta(h,\tg)=
\frac{\Delta(h,\Phi^+)}{\Delta(g,\Phi^+)}
\frac{\epsilon_r(h,\Phi^+)}{\epsilon_r(\tg,\Phi^+)}e^{\rho_r(\Phi^+)}(t)
\tau(h_0,\tilde g)
\end{equation}
Let $\Phi^+_r(h_0)=\{\alpha\in\Phi_r\,|\, e^\alpha(h_0)>1\}$. 
By \eqref{e:epr12} this gives
\begin{equation}
\Delta(h,\tg)=
\frac{\Delta(h,\Phi^+)}{\Delta(\tg,\Phi^+)}
e^{\rho(\Phi^+_r(h_0))}(t)
\tau(h_0,\tg)
\end{equation}
which is obviously independent of $\Phi^+$.
\end{proof}

\begin{lemma}
For all $(h,\tg)\in X(G,\tG)$ we have
\begin{equation}
\begin{aligned}
\Delta(h,\tg)=c(h,\tg)\frac{|D(h)|^{\frac12}}{|D(\tg)|^{\frac12}}
\end{aligned}
\end{equation}
and
\begin{equation}
\Delta(h,\tg)=c'(h,\tg)\frac{\Delta(h,\Phi^+)}{\Delta(\tg,\Phi^+)}
\end{equation}
where $c(h,\tg)^4=c'(h,\tg)^2=1$.

Suppose $\tilde x\in\tG$ and let $x=p(\tilde x)$. Then
\begin{equation}
\Delta(xhx\inv,\tilde x\tilde g\tilde x\inv)=\Delta(h,\tilde g).
\end{equation}

\end{lemma}
This is immediate.

In the next section we will drop the assumption that $G$ is acceptable,
in which case 
$\Delta(h,\Phi^+)$  and $\Delta(\tg,\Phi^+)$ are  not necessarily well
defined. With this in mind we rewrite 
\eqref{e:transferspecial} as follows:
\begin{equation}
\begin{aligned}
\Delta(h,\tg)
&=\frac
{\Delta(h,\Phi^+)\epsilon_r(h,\Phi^+)}
{\Delta(\tg,\Phi^+)\epsilon_r(\tg,\Phi^+)}
\frac{\tilde\chi(\tg)}{\chi_0(h)}\\
&=\frac
{\Delta^0(\Phi^+,h)\epsilon_r(h,\Phi^+)}
{\Delta^0(\Phi^+,\tg)\epsilon_r(\tg,\Phi^+)}
\frac{\tilde\chi(\tg)e^\rho(h)}{\chi_0(h)e^\rho(\tg)}\\
&=\frac
{\Delta^1(\Phi^+,h)}
{\Delta^1(\Phi^+,\tg)}
\frac{\tilde\chi(\tg)}{(\chi_0 e^\rho)(h)}\\
\end{aligned}
\end{equation}
since $p(\tg)=h^2$.
Let $\chi=\chi_0e^\rho$ and assume $\Phi^+$ is special
(Definition~\ref{d:special}).
Then
by \eqref{e:zetacx2} and \eqref{e:chi0}
for $h\in \Gamma_r(H)$ we have
\begin{equation}
\chi(h)=e^{\rho_r}(h)e^\rho(h)=\zetacx(h).
\end{equation}
Thus $\chi$ is defined by

\begin{subequations}
\renewcommand{\theequation}{\theparentequation)(\alph{equation}}  
\label{e:generaltransfer}
\begin{equation}
\chi(h)=
\begin{cases}
(\tilde\chi^2e^\rho)(h)&h\in H^0\\
\zetacx(h)&h\in \Gamma_r(H).
\end{cases}
\end{equation}
Define 
\begin{equation}
\Gamma(h,\tg)=\tilde\chi(\tg)/\chi(h), 
\end{equation}
and then (still assuming $G=\oG$ is connected and semisimple)
\begin{equation}
\Delta(h,\tg)=\frac{\Delta^1(\Phi^+,h)}{\Delta^1(\Phi^+,\tg)}\Gamma(h,\tg).
\end{equation}
\end{subequations}

This has the advantage that the quotient is defined in general and
this is the starting point of the definition of transfer factors.
There are several issues in extending it to the general case.  First
of all $e^\rho$ is only defined on $H_d^0$ (recall $G_d$ is
acceptable), so \eqref{e:generaltransfer} only defines $\chi$, and
hence $\Gamma$ and $\Delta$, on $\Gamma_r(H)H_d^0=H\cap G_d^0$
(cf.~\eqref{e:Gamma(H)}). This may be a proper subgroup of $H$.
Furthermore if $\oG\ne G$ then some additional modification is
required.  We take up these issues in the next two sections.

\sec{Characters of Cartan Subgroups.}
\label{s:characters}

Let $(\tG, G, \oG)$ be an admissible triple. In the previous section
we defined the factor $\Gamma(h,\tg)$ \eqref{e:generaltransfer}(b)
if $G=\oG$ is connected and semisimple.
In this section we consider the general case, which we use in
the subsequent section to define transfer factors in general.

Suppose $H$ is a  Cartan subgroup of $G$ and let 
$\Phi^+$ be a set of positive roots.
Suppose $\wt\chi$ is a genuine character of $Z(\tH)$. 
Since $G_d(\C)$ is acceptable $\rho$
exponentiates to a character of $(H\cap G_d)^0$, and $\wt\chi^2e^\rho$
is defined on $(H\cap G_d)^0$. 
Let $\oH$ be the corresponding Cartan subgroup of $\oG$ and set
$\oH_d=\oH\cap \oG_d$.
We now use property (3c) of  Definition \ref{assumptions}, which says that 
$\wt\chi^2e^\rho$ factors to $\oH_d^0$.
The canonical character $\zetacx$ of $\Gamma(\oH)$ was defined in Section
\ref{s:notation}. 

We use \eqref{e:generaltransfer} as our starting point.

\begin{definition}
\label{d:Gamma}
Let $H$ be a  Cartan subgroup of $G$. Choose
\begin{enumerate}
\item a special set of positive roots $\Phi^+$ (Definition \ref{d:special})
\item a genuine character $\tilde\chi$ of $Z(\tH)$,
\item a character $\chi$ of $\oH$.
\end{enumerate}
Assume these satisfy:
\begin{subequations}
\label{e:cond}
\renewcommand{\theequation}{\theparentequation)(\alph{equation}}  
\begin{align}
\chi(h)&=(\tilde\chi^2e^{\rho})(h), \quad(h\in \oH_d^0)
\label{e:cond1}\\
\chi(h)&=\zetacx(h)\quad(h\in \Gamma_r(\oH)).
\label{e:cond2}
\end{align}
\end{subequations}

Fix $(\Phi^+,\tilde\chi,\chi)$  satisfying these conditions.
Suppose  $(h,\tg)\in X(\oH,\tH)$ (cf. \eqref{e:XH}).
Then we define
\begin{equation}
\label{e:Gamma}
\Gamma(\tilde\chi,\chi)(h,\tg)=\frac{\tilde\chi(\tg)}{\chi(h)}.
\end{equation}
Let
\begin{subequations}
\label{e:ST}
\renewcommand{\theequation}{\theparentequation)(\alph{equation}}  
\begin{align}
\caS(H,\Phi^+)&=\{(\tilde\chi,\chi)\}\text{ satisfying \eqref{e:cond}(a) and (b)} \} ,\\
\caT(H,\Phi^+)&=\{\Gamma(\tilde\chi,\chi)\,|\,(\tilde\chi,\chi)\in\caS(H,\Phi^+)\}\\
\caS(H)&=\cup_{\Phi^+}\caS(H,\Phi^+)\\
\caT(H)&=\cup_{\Phi^+}\caT(H,\Phi^+).
\end{align}
\end{subequations}
(The unions are over special sets of positive roots).
\end{definition}

\begin{example}
\label{e:sl2_a}
Let $G=\oG=SL(2,\R)$ and let $H$ be the diagonal Cartan subgroup,
which we identify with $\R^\times$. Let $\Phi^+=\{\alpha\}$ where
$\alpha(x)=x^2$. If $p(\tilde x)=x$  and $\epsilon=\pm1$ then a short calculation gives
\begin{equation}
\Gamma(x,\tilde x^2\epsilon)=x\inv\epsilon.
\end{equation}
This only depends on $\Phi^+$; the other choice of $\Phi^+$ gives
$x\epsilon$.
Compare Example \ref{ex:sl2_special}.
\end{example}

We begin with some elementary properties. Fix $H$, $(\tilde\chi,\chi)$ as
above and let $\Gamma=\Gamma(\tilde\chi,\chi)$.    First of all we
obtain a
character $\lambda$ of $Z_0(H)$.

\begin{lemma}
\label{l:mu}
Let $\lambda=(\tilde\chi^2/\chi)$, viewed as a  character of $Z_0(H)$,
i.e. $\lambda(h)=\tilde\chi^2(\tilde h)/\chi(\overline p(h))$ where $p(\tilde h)=h$.
Then
\begin{subequations}
\renewcommand{\theequation}{\theparentequation)(\alph{equation}}  
\begin{equation}
\label{e:rho}
\lambda(h)=e^{-\rho}(h)\quad(h\in Z_0(H)\cap G_d).
\end{equation}
Also let
\begin{equation}
\mu=-d\lambda-\rho=d\chi-2d\wt\chi-\rho.
\end{equation}
Then $\mu|_{\h_d}=0$ and we identify it with an element of $\z(\C)^*$;
as such it is the negative of the differential of $\lambda |_{Z(G)^0}$.

Furthermore for all $(h,\tg)\in X(\oH,\tH)$
\begin{equation}
\label{e:Gamma^2}
\Gamma(h,\tg)^2=\lambda(\phi(h)).
\end{equation}
\end{subequations}
In particular if $h\in \oG_d$ then $\Gamma(h,\tg)^2=e^{-2\rho}(h)$.
\end{lemma}

\begin{proof}
If $h\in H_d^0$ then by
\ref{e:cond}(b)  $\chi(\overline p(h))=(\wt\chi^2e^\rho)(h)$.
The same identity holds for $h\in \Gamma_r(H)\cap Z_0(H)$ by 
\eqref{e:zetacx(h)2}.
In both cases we conclude $\lambda(h)=\wt\chi^2(h)/\chi(\overline
p(h))=e^{-\rho}(h)$. This holds on $\Gamma_r(H)H_d^0\cap Z_0(H)$,
which equals $Z_0(H)\cap G_d$ by 
\eqref{e:Gamma(H)}. This gives (a).

For (b), 
if $(h,\tg)\in X(\oH,\tH)$ then
$\Gamma(h,\tg)^2=\tilde\chi^2(\tg)/\chi(h^2)$.
Since $p(\tilde g)=\phi(h)$ and $\overline p(\phi(h))=h^2$ this equals
$\tilde\chi^2(\phi(h))/\chi(\phi(h))=\lambda(\phi(h))$. The final equality
follows from (a) and (b) and the fact that $\overline p(\phi(h))=h^2$.
\end{proof}

Recall \eqref{e:tau} for $(h,\tg)\in X(\oH,\tH)$, $\tau(h,\tg)=\wt
h^2\tg\inv=\pm1$.  
\begin{lemma}
\label{e:Gammaformulas}
Suppose $(h,\tg)\in X(\oH,\tH)$.

\noindent (1) Assume $h\in \overline p(Z_0(H))$, and choose $y\in Z_0(H)$
satisfying $\overline p(y)=h$. Then
with $\lambda$ as in \eqref{e:mu},
\begin{equation}
\Gamma(h,\tg)=\lambda(y)\tau(y,\tg).
\end{equation}

\noindent (2) If $h\in\oG_d$  write $h=th_0$ with $t\in\Gamma_r(\oH)$,
$h_0\in\oH_d^0$, and choose $y_0\in H^0$ satisfying $\overline
p(y_0)=h_0$. 
 Then 
\begin{equation}
\Gamma(h,\tg)=\zeta_{cx}(t)e^{-\rho}(y_0)\tau(y_0,\tilde g),
\end{equation}
\label{e:GammaonGd}
and this is independent of all choices except that of $\Phi^+$.

\noindent (3)
In general write $h=th_0$ with $t\in
\Gamma(\oH)$ and $h_0\in \oH^0$ \eqref{e:Gamma(H)}. Write
$\tg=\tilde\gamma\tg_0$ where $p(\tilde\gamma)=\phi(t)$ and
$p(\tg_0)=\phi(h_0)$. Finally choose $y_0\in H^0$ satisfying
$\overline p(y_0)=h_0$. Then
\begin{equation}
\label{e:Gammahtg}
\Gamma(h,\tg)=\frac{\tilde\chi(\tilde\gamma)}{\chi(t)}\lambda(y_0)\tau(y_0,\tg_0).
\end{equation}
Note that $\tilde\chi(\tilde\gamma)^4=1$ and $\chi(t)^2=1$.

\end{lemma}

\begin{proof}
For (1) choose $\tilde y\in p\inv(y)$ and write 
$\tilde g=\tilde y^2\tau(y,\tilde g)$ (cf. \ref{e:tau}). Then
$\Gamma(h,\tg)=\tilde\chi(\tilde
y^2)\tau(y,\tg)/\chi(h)=(\tilde\chi^2(y)/\chi(\overline
p(y)))\tau(y,\tg)=\lambda(y)\tau(y,\tg)$. 
Part (2) follows from this, \eqref{e:cond2} and \eqref{e:rho}.
Part (3) is similar, noting that
$\Gamma(h,\tg)=\Gamma(t,\tilde\gamma)\Gamma(h_0,\tg_0)$.
\end{proof}

\begin{example}
\label{ex:gln}
Let $G=\oG=GL(n,\R)$ and suppose $H$ is the diagonal (split) Cartan
subgroup. For simplicity we assume $n$ is even, in which case
$Z_0(H)=H^0$.
Note that $\Gamma(H)=\{\diag(\epsilon_1,\dots,\epsilon_n)\}$
with $\epsilon_i=\pm1$, and $\Gamma_r(H)$ is the subset of $\Gamma(H)$
of elements of determinant $1$.
If $h\in H$ write $h=th_0$ with $t\in\Gamma(H)$ and $h\in
H^0$. 
Let $\Phi^+$ be any (necessarily special) set of positive roots.

The cover $Z(\tH)\rightarrow Z_0(H)$ splits; let $\tilde\chi$ be the
unique, genuine quadratic character of $Z(\tH)$.
Note that $\rho$ exponentiates to $H^0$, and we can define
$\chi(th_0)=e^{\rho}(h_0)$.
Then $(\tilde\chi,\chi)\in\caS(H,\Phi^+)$ and by \eqref{e:Gammahtg} we have
\begin{equation}
\label{e:gln}
\Gamma(h,\tg)=e^{-\rho}(h_0)\tau(h_0,\tg).
\end{equation}

Given $\tilde\chi$, the only other allowed choice of $\chi$  and $\Gamma$
is obtained by multiplying by
$\nu (\det(h))$ for $\nu $ a character of $\R^\times$ (see Lemma
\ref{l:simplytransitive}).

\end{example}

\begin{lemma}
\label{l:elementary}
Fix a special set of positive roots $\Phi^+$.
\hfil\break
\noindent (1) 
Suppose $\nu$ is a character of $Z_0(H)$ and $(\tilde\chi,\chi)\in \caS(H,\Phi^+)$.
Then $(\tilde\chi(\nu \circ p),\chi(\nu \circ\phi))\in \caS(H,\Phi^+)$ and 
\begin{equation}
\Gamma(\tilde\chi(\nu\circ p),\chi(\nu\circ\phi))=\Gamma(\tilde\chi,\chi).
\end{equation}

\noindent (2)
Suppose $\Phi^+_1$ is another special set of positive roots. 
There exists a character $\tau$  of $\oH$  such that for all $(\tilde\chi,\chi)\in \caS(H,\Phi^+)$,
$(\tilde\chi,\chi\tau)\in \caS(H,\Phi^+_1)$, and 
\begin{equation}
\Gamma(\tilde\chi,\chi)(\tg,h)
=
\Gamma(\tilde\chi,\chi\tau)(\tg,h)
\end{equation}
for all $(h,\tg)\in X(\oH,\tH)$ satisfying $h\in \Gamma(\oH)Z(\oG)$.
\end{lemma}

\begin{proof}
The first part is straightforward. For the second define
$\tau(h)=e^{\rho(\Phi _1^+)-\rho(\Phi^+)}(h)$ for $h\in \oH_d^0$, and
$\tau(h)=1$ for $h\in \Gamma(\oH)Z(\oG)$. It is easy to see that the
fact that both $\Phi^+$ and $\Phi^+_1$ are special implies 
$e^{\rho(\Phi^+_1)-\rho(\Phi^+ )}(h)=1$ for all $h\in
\Gamma(\oH)Z(\oG)\cap \oH^0_d$, so $\tau$ is well defined. It follows
easily that it has the desired properties.  
\end{proof}

There is a natural action of the characters of $\oG$ on $\caT(H,\Phi^+)$:
if $\psi$ is a one-dimensional representation of $\oG$ define
\begin{equation}
\psi\cdot\Gamma(\tilde\chi,\chi)
=\Gamma(\tilde\chi,\chi\psi|_{\oH}).
\end{equation}

\begin{lemma}
\label{l:simplytransitive}
The space $\caT(H,\Phi^+)$ is non-empty, and 
the action of the group of  characters of $\oG$ on
$\caT(H,\Phi^+)$ is transitive.
If $H$ is maximally split it is simply transitive.
\end{lemma}

\begin{proof}
Fix any genuine character $\tilde\chi$ of $Z(\tH)$. 
For existence it is enough to show that for any special set of
positive roots,
\begin{equation}
\label{e:zetacx(h)}
\zetacx(h)=e^\rho\tilde\chi^2(h)\quad h\in \Gamma_r(\oH)\cap \oH_d^0.
\end{equation}
Note that $p\inv(\oH_d^0)\subset p\inv(\oH^0)=H^0C\subset H^0Z_0(G)=Z_0(H)$ by
\ref{assumptions}(3b)) and  Proposition \ref{p:ZH}.
Pulling back to $G$ it is therefore  enough to show
\begin{equation}
\label{e:zetacx(h)2}
\zeta_{cx}(h)=e^\rho\tilde\chi^2(h)\quad h\in\Gamma_r(H)\cap Z_0(H).
\end{equation}

\begin{lemma}
\label{l:commutes}
$
\wt\chi^2(h)=e^{\rho_r}(h)
$
for all $h\in Z_0(H)\cap \Gamma_r(H)$.
\end{lemma}

\begin{proof}
Fix $h\in Z_0(H)\cap \Gamma_r(H)$.
By \eqref{e:GammaR} write
$h=\prod_{i=1}^n m_i$
where $\{\alpha_1,\dots, \alpha_n\}$ is a subset of the simple roots 
for $\Phi_r^+$, and $m_i=m_{\alpha_i}$.
Then $e^{\rho_r}(h)=(-1)^n$.

For each $i$ choose an inverse image $\wt m_i$ of $m_i$ in $\tH$.
By Lemma \ref{adams1} and the discussion preceding it $\wt m_i^2=-1$
for all $i$. 
Assume for the moment that $\wt m_i,\wt m_j$ commute for all $i,j$.
Then
\begin{equation}
\begin{aligned}
\wt\chi^2(h)&=\wt\chi([\wt m_1\dots\wt m_n]^2)\\
&=\wt\chi(\wt m_1^2\dots \wt m_n^2)\\
&=\wt\chi((-1)^n)=(-1)^n
\end{aligned}
\end{equation}
since $\wt\chi$ is genuine, proving the result.

Thus it is enough to show the $\wt m_i$ commute.
Using the commutator notation of Section \ref{s:cartan},
the fact that $\{m_i,h\}=1$  implies 
by \eqref{e:tau3} that, for each $i$,
$\alpha_i$ is
not orthogonal to an even number of other $\alpha_{j}$.
Consider the subdiagram of the Dynkin diagram
of $\Phi_r$
consisting of the $\alpha_i$ ($1\le i\le n$): every node is
adjacent to an even number of other nodes. This implies the nodes are
all orthogonal, and  the $\wt m_i$ commute.
\end{proof}

Therefore
\begin{equation}
\tilde\chi^2e^\rho(h)
=
e^{\rho_r}(h)e^{\rho_r+\rho_{cx}+\rho_i}(h)\\
=
e^{\rho_{cx}}(h)
\end{equation}
since $e^{2\rho_r}(h)=e^{\rho_i}(h)=1$. 
This equals $\zeta_{cx}(h)$ by \eqref{e:zetacx2}.

Now suppose $\Gamma_1,\Gamma_2\in\caT(H,\Phi^+)$.
If $h\in \oH$ then $\psi_0(h)=\Gamma_1(h,\tg)/\Gamma_2(h,\tg)$ is
independent of the choice of $\tg$ so that $(h,\tg)\in X(\oH,\tH)$,
and is a character of $\oH$. By 
\eqref{e:GammaonGd} $\psi_0(h)=1$ for $h\in \oG_d$.

We have a surjection
\begin{equation}
\label{e:ogogd}
\Hom(\oG/\oG_d,\C)\rightarrow \Hom(\oH/\oH_d,\C)
\end{equation}
dual to the injection
$\oH/\oH_d\hookrightarrow \oG/\oG_d$,
which is an isomorphism if $H$ is maximally split.
If we choose any preimage $\psi$ of $\psi_0$ in \eqref{e:ogogd} then
it is easy to see that
$\Gamma_1=\psi\Gamma_2$.
It follows that characters of $\oG$ act transitively on
$\caT(H,\Phi^+)$, and this action is simply transitive if $H$ is
maximally split.
\end{proof}

We now need to choose elements of $\caT(H,\Phi^+)$ consistently
for all $H$. We do this by reducing to the most split Cartan
subgroup. For this we need a Lemma.

\begin{lemma}
\label{l:winv}
Suppose $(\tilde\chi,\chi)\in \caS(H,\Phi^+)$. 
Then
\begin{equation}
\chi(h)=\chi(wh)\quad(w\in W(\oG,\oH),h\in \Gamma(\oH)Z(\oG)).
\end{equation}
\end{lemma}

\begin{proof}
We need to show $\chi(g(wg)\inv)=1$ for $g\in \Gamma(\oH)$.
This is obvious if $w\in W(\Phi_i)$. If $w\in W(\Phi_r)$ then 
$g(wg)\inv\in \Gamma_r(\oH)$. Then
$\chi(g(wg)\inv)=\zetacx(g(wg)\inv)=1$ by \eqref{e:cond2}
and \eqref{e:zetacx1}.

By \cite[Proposition 3.12]{ic4} $W(\oG(\C),\oH(\C))^\theta$, which contains
$W(\oG,\oH)$, is generated by $W(\Phi_i),W(\Phi_r)$ and  elements of
the form $s_\alpha s_{\theta\alpha}$  where $\langle\alpha,\theta\ch\alpha\rangle=0$ and
$s_\alpha s_{\theta\alpha}\Phi _r^+ = \Phi _r^+$.
So it is enough to show
$\chi(g(wg)\inv)=1$ for $w$ of this form.

Write $g=\overline\exp{\pi iX}$ with $X\in \a$. 
Since $g^2=1$, $\alpha(X)\in \Z$.
Then 
\begin{equation}
g(wg)\inv=\overline\exp{\pi i\alpha(X)(\ch\alpha+\theta\ch\alpha)}\in \Gamma(\oH)\cap \oH_d^0
\end{equation}
Let $\rho=\rho(\Phi^+)$.
Then
\begin{equation}
\begin{aligned}
\chi((wg)g\inv)
&=
(\tilde\chi^2e^{\rho})(g(wg)\inv)\quad \text{(by }\eqref{e:cond1})\\
&=
(-1)^{\alpha(X)\langle 2d\tilde\chi,\ch\alpha+\theta\ch\alpha\rangle}e^{\rho-w\inv\rho}(g)
\end{aligned}
\end{equation}
By \cite[Lemma 6.11]{dualityonerootlength} 
$\langle d\tilde\chi,\ch\alpha+\theta\ch\alpha\rangle\in\Z$.
The fact that $\Phi^+$ is special and $w\Phi _r^+ = \Phi _r^+$ implies
$\rho-w\inv\rho$ is a sum of imaginary roots  and 
terms $\beta-\theta\beta$ with $\beta$ complex.
Therefore $e^{w\inv\rho-\rho}(g)=1$.
\end{proof}

Fix a maximally split Cartan subgroup $H_s=T_sA_s$ of $G$ and
$(\tilde\chi_s,\chi_s)\in \caS(H_s)$. 

Suppose $H$ is a Cartan subgroup, and write
$\oH=\oT\oA$ as usual.
There exists $x\in\oG$ such that
$x\oA x\inv\subset \oA_s$,
and therefore $x\Gamma(\oH)x\inv \subset \Gamma(\oH_s)$. 
Suppose $(h,\tg)\in X(\oH,\tH)$ where $h\in \Gamma(\oH)Z(\oG)$.
By Lemmas \ref{bs4}(5)
and \ref{centerh}
$\phi(h)\in Z_0(G)$
so $\tg\in Z(\tG)\subset Z(\tilde{H_s})$.
Therefore $(xhx\inv,\tg)\in X(\oH_s,\tH_s)$. 
Define
\begin{equation}
\label{e:Gammas}
\Gamma_s(h,\tg)=\Gamma(\tilde\chi_s,\chi_s)(xhx\inv,\tg).
\end{equation}

By the Lemma this is independent of the choice of $x$.

\begin{proposition}
\label{p:Gammamain}
Fix a maximally split Cartan subgroup $H_s$ and
$(\tilde\chi_s,\chi_s)\in \caS(H_s)$.
Suppose $H$ is a Cartan subgroup, 
$\tilde\chi$ is a genuine character of $Z(\tH)$, and $\Phi^+$ is a
special set
of positive roots for $H$. Then there exists a unique character $\chi$ of $\oH$
such that  $(\tilde\chi,\chi)\in \caS(H,\Phi^+)$ and 
\begin{equation}
\label{e:Gammamain}
\Gamma(\tilde\chi,\chi)(h,\tg)=\Gamma_s(h,\tg)
\end{equation}
for all  $(h, \tg) \in X(\oH,\tH)$ with $h\in\Gamma(\oH)Z(\oG)$.
\end{proposition}

Explicitly \eqref{e:Gammamain} is
\begin{equation}
\label{e:tildechichi}
\frac{\wt\chi(\tg)}{\chi(h)}=\frac{\wt\chi_s(\tg)}{\chi_s(xhx\inv)}
\end{equation}
or alternatively
\begin{equation}
\label{e:compatible}
\chi(h)=(\wt\chi/\wt\chi_s)(\phi(h))\chi_s(xhx\inv).
\end{equation}

Before giving the proof of the Proposition we give the main definition
of this section. 

\begin{definition}
\label{d:compatible}
Fix a split Cartan subgroup $H_s$ and $(\wt \chi_s,\chi_s)\in\caS(H_s)$.
Suppose $H$ is a Cartan subgroup and $\Phi^+$ is a special set of positive
roots.

\smallskip
\noindent (1) Let $\caS(H,\Phi^+,\wt\chi_s,\chi_s)$ be the set of pairs
$(\wt\chi,\chi)\in \caS(H,\Phi^+)$  satisfying \eqref{e:Gammamain}. 

\smallskip
\noindent (2) Suppose $\wt\chi$ is  any genuine character of $\tH$. 
Choose $\chi$ so that $(\wt\chi,\chi)\in
\caS(H,\Phi^+,\wt\chi_s,\chi_s)$ and let

\begin{subequations}
\renewcommand{\theequation}{\theparentequation)(\alph{equation}}  
\label{e:compatible2}
\begin{equation}
\Gamma(H,\Phi^+)=\Gamma(\wt \chi,\chi)\in\caT(H,\Phi^+).
\end{equation}
Also define
\begin{equation}
\label{e:muHPhi+}
\lambda(H,\Phi^+)=\tilde\chi^2/\chi\in \widehat{Z_0(H)}
\end{equation}
and
\begin{equation}
\label{e:mu}
\mu=-d\lambda-\rho=d\chi-2d\wt\chi-\rho\in\z(\C)^*
\end{equation}
\end{subequations}
(cf.~Lemma \ref{l:mu}). 
If it is necessary to indicate the dependence on $(\wt\chi_s,\chi_s)$
we write $\Gamma(\wt\chi_s,\chi_s,H,\Phi^+)$, 
$\lambda(\wt\chi_s,\chi_s,H,\Phi^+)$ and $\mu(\wt\chi_s,\chi_s)$.
\end{definition}
It follows easily from 
Lemma \ref{l:elementary} that $\Gamma(H,\Phi^+)$ and $\lambda(H,\Phi^+)$
are independent of the choice of $\tilde\chi$. It is also easy to see
that $\mu=d\chi_s -2d\wt\chi_s$ restricted to $\z(\C)^*$, and is
therefore independent of the choice of $(H,\Phi^+)$.

We  now prove Proposition \ref{p:Gammamain}.     
Using Cayley transforms 
(Section \ref{s:cartan}) we will reduce to the following Lemma.
Suppose $H_\alpha$ is a Cartan subgroup, with real root $\alpha$,
and $H_\beta$ is the Cayley transform of $H_\alpha$, with noncompact
imaginary root $\beta$.   See the discussion preceding \eqref{e:HaHb}.
Fix $c\in\Aut(\g)$ satisfying \eqref{e:c}(a-c).

\begin{lemma}
\label{l:HHalpha}
Suppose $\Phi^+_\alpha$ is a special set of positive roots for $H_\alpha$
and $(\tilde\chi_\alpha,\chi_\alpha)\in \caS(H_\alpha,\Phi^+_\alpha)$. 
Suppose $\Phi^+_\beta$ is a special set of positive roots for
$H_\beta$, and $\tilde\chi_\beta$ is a genuine character of $Z(\tH_\beta)$.
Then there exists a unique character $\chi_\beta$ of $H_\beta$ 
such that  $(\tilde\chi_\beta,\chi_\beta)\in \caS(H_\beta,\Phi^+_\beta)$ and 
\begin{equation}
\label{e:Gammaequals}
\Gamma(\tilde\chi_\alpha,\chi_\alpha)(h,\tg)=
\Gamma(\tilde\chi_\beta,\chi_\beta)(h,\tg)
\end{equation}
for all $(h,\tg)\in X(\oH_\beta,\tH_\beta)$ satisfying  $h\in
\Gamma(\oH_\beta)Z(\oG)$. 

In addition assume $c^*(\Phi^+_\alpha)=\Phi_\beta^+$. Then 
\eqref{e:Gammaequals} holds for all $h\in \oH_\alpha\cap \oH_\beta$. 
\end{lemma}

\begin{proof}
To avoid runaway notation let $\oH=\oH_\beta$ for a moment. 
The character $\chi_\beta$ is determined on $\Gamma(\oH)Z(\oG)$  by
\eqref{e:Gammaequals} and on $\oH_d^0$ by \eqref{e:cond}(a).
Since $\oH^0\subset Z(\oG)\oH_d^0$, by
\eqref{e:Gamma(H)} $\oH=\Gamma(\oH)Z(\oG)\oH_d^0$, and
uniqueness is immediate.

For existence, by Lemma \ref{l:elementary} it is enough to prove this for a single
choice of $\Phi^+_\alpha$ and $\Phi^+_\beta$. 
It is not hard to see we can choose  $\Phi^+_\alpha$ special so that 
$\Phi^+_\beta=c^*(\Phi^+_\alpha)$ is also special. This implies
$\alpha$ is simple for $\Phi^+_r$ (cf.~Lemma \ref{sstable}). 

Fix $\tilde\chi_\beta$. By \ref{l:elementary}(1)  we may
choose $\tilde\chi _\alpha$ arbitrarily, so by Lemma
\ref{l:cayleynonlinear} choose $\tilde\chi_\alpha$ to be the Cayley
transform of $\tilde\chi_\beta$.  
Since $\tilde\chi_\beta$ is genuine and $\beta $ is a
noncompact imaginary root,  $\langle
d\tilde\chi_\beta,\ch\beta \rangle\in \Z+\frac12$ by \cite[Lemma 6.11]{dualityonerootlength}.  
Therefore $\langle
d\tilde\chi_\alpha,\ch\alpha \rangle = \langle
d\tilde\chi_\beta,\ch\beta \rangle\in \Z+\frac12$.  Let $\chi _\alpha
$ be any character of $\oH_\alpha $ with
$(\tilde\chi_\alpha,\chi_\alpha)\in \caS(H_\alpha,\Phi^+_\alpha)$.
Then
\begin{equation} 
\chi_\alpha (m_\alpha ) = \zeta_{cx} (m_\alpha )
 = (-1)^{\langle \rho - \rho _r,\ch\alpha \rangle} = -(-1)^{\langle \rho,\ch\alpha \rangle}
\end{equation}
since $\alpha$ is simple for $\Phi^+_r$.
Furthermore
\begin{equation} \langle d\chi_\alpha ,\ch\alpha\rangle 
=\langle 2d\tilde\chi_\alpha +\rho ,\ch\alpha \rangle  \in 2\Z +1 + \langle  \rho ,\ch\alpha \rangle.\end{equation}
Thus $\chi _\alpha$ satisfies the conditions of
Lemma  \ref{l:cayley}.
Let $\chi_\beta$ be the (unique) Cayley transform of $\chi_\alpha$.

It is enough to show $(\wt\chi_\beta,\chi_\beta)\in
\caS(H_\beta,\Phi^+_\beta)$, for then it is obvious from the
definition of Cayley transforms that \eqref{e:Gammaequals} holds for all $h\in \oH_\alpha\cap \oH_\beta$. 

We have
\begin{equation}
\begin{aligned}
d\chi_\beta&=\Ad^*(c  )(d\chi_\alpha)\\
&=\Ad^*(c)(2d\wt\chi_\alpha+\rho(\Phi^+_\alpha))\\
&=2 d\tilde\chi_\beta+\rho(\Phi^+_\beta)
\end{aligned}
\end{equation}
by the choice of $\Phi^+_\alpha$ and $\Phi^+_\beta$.
This verifies \eqref{e:cond1}.

We now verify \eqref{e:cond}(b). 
Suppose $\gamma$ is a real root of $\oH_\beta$.
Then $m_\gamma\in \oH_\alpha\cap\oH_\beta$, and
$\chi_\beta(m_\gamma)=\chi_\alpha(m_\gamma)=\zeta_{cx}(\oG,\oH_\alpha)(m_\gamma)$,
by \eqref{e:cond}(b) for $(\wt\chi_\alpha,\chi_\alpha)$. We want to
show this equals $\zeta_{cx}(\oG,\oH_\beta)(m_\gamma)$, i.e.
\begin{equation}
\label{e:zetacxalphabeta}
\zeta_{cx}(\oG,\oH_\alpha)(m_\gamma)
=
\zeta_{cx}(\oG,\oH_\beta)(m_\gamma).
\end{equation}

If the simple factor containing $\alpha,\beta$ is of type $G_2$ an
explicit calculation shows that both sides are $1$, so assume $G$ is
simply laced.

The root $\wt\gamma=\Ad^*(c\inv)(\gamma)$ of $\oH_\alpha$ is also
real, and $\langle\wt\gamma,\ch\alpha\rangle=0$.
Suppose $\delta$ is a complex root of
$\oH_\alpha$ with $\langle\delta,\ch{\wt\gamma}\rangle \not = 0$.    Then the corresponding root of $\oH_\beta$ is also complex.     Conversely, suppose $\delta$ is a complex root of
$\oH_\beta$ with $\langle\delta,\ch{\gamma}\rangle \not = 0$.    Then the corresponding root $\wt\delta$ of $\oH_\alpha$ is complex unless $s_\beta\delta=\sigma\delta$, in
 which case it is real.   In this case, let $\delta ' = s_{ \gamma}\delta$.    It is another complex root of
 $\oH_\beta$ with $\wt\delta '$ real.    Since $\langle\delta,\ch{\gamma}\rangle \not = 0$, $ \delta ' \not = \delta$, and since $\gamma $ is real, $ \delta ' \not = \sigma \delta$.
Since $\delta \not = \pm \gamma,  \delta ' \not = -\delta$.   Suppose $\delta '  =  - \sigma\delta = -s_\beta\delta$.   
Then
\begin{equation}
2\delta=\langle\delta,\ch{\gamma}\rangle\gamma+
\langle\delta,\ch\beta\rangle\beta.
\end{equation}
This cannot happen since in the simply laced case, two orthogonal
roots have no other 
roots in their (real) span.

Therefore both $\delta$ and $s_{ \gamma}(\delta)$ contribute to
$\zeta_{cx}(\oG,\oH_\beta)(m_\gamma)$, and their total 
contribution is
$\delta(m_\gamma)(s_{ \gamma}\delta)(m_\gamma)=\delta(m_\gamma)\delta(m_\gamma)=1$.

It follows that the terms in 
$\zeta_{cx}(\oG,\oH_\alpha)$ and 
$\zeta_{cx}(\oG,\oH_\beta)$ are the same, with the exception of those
just discussed, which are $1$.
\end{proof}

\begin{remark}
Note that $\Gamma(\tilde\chi_\alpha,\chi_\alpha)(m_\alpha,1)=\chi_\alpha(m_\alpha)$ and
$\Gamma(\tilde\chi_\beta,\chi_\beta)(m_\alpha,1)=\chi_\beta(m_\alpha)$.
As in the proof of the Lemma
$\chi_\beta(m_\alpha)=-(-1)^{\langle\rho,\ch\alpha\rangle}$. The
equality of the Lemma then implies
$\chi_\alpha(m_\alpha)=\zeta_{cx}(m_\alpha)$. This motivates  condition
\eqref{e:cond2}.

\end{remark}

\begin{proof}[Proof of the Proposition]
This is now straightforward. If $H$ is maximally split there is
nothing to prove.
If $H$ is obtained from $H_s$ by a series of Cayley transforms we
conclude the result by a repeated application of 
Lemma \ref{l:HHalpha}.
It is easy to check that if the conditions of the
Proposition hold for a Cartan subgroup $H$, they hold for every
$G$-conjugate of $H$. 
Up to conjugacy
every Cartan subgroup is conjugate to one obtained by a series of
Cayley transforms from $H_s$, and the result follows.
\end{proof}

For later use we note a consequence of the last part of the proof of
Lemma \ref{l:HHalpha}.

\begin{lemma}
\label{l:cayleyH1H2}
Suppose $H_1=T_1A_1$ and $H_2=T_2A_2$ are Cartan subgroups with
$A_1\subset A_2$.
Suppose $\alpha\in \Phi_r(G,H_1)$. 

\noindent (1) 
$$
\zetacx(G,H_1)(m_{\alpha})=\zetacx(G,H_2)(m_{\alpha}).
$$

\noindent (2) Let $M_1=\Cent_G(A_1)$ and suppose $\alpha\in \Phi_r(M_1,H_2)$.
Then
$$
\zetacx(G,H_2)(m_{\alpha})=\zetacx(M_1,H_2)(m_{\alpha}).
$$
\end{lemma}

\begin{proof}
The first equality follows from a repeated application of 
\eqref{e:zetacxalphabeta}.
The second is similar.
Choose a set $S$ of complex roots such that 
\begin{equation}
\{\beta\in \Phi_{cx}(G,H_2)\backslash\Phi_{cx}(M_1,H_2)
\,|\, \langle\beta,\ch\alpha\rangle\ne 0\}
=
\{\pm\beta,\pm\sigma\beta\,|\,\beta\in S\}.
\end{equation}
If $\beta\in S$ then $s_\alpha\beta$ is also complex,
is not contained in $\Phi(M_1,H_2)$, and is not equal to 
$\pm\beta,\pm\sigma\beta$
Then $S$ is can be written as a union over  pairs $\{\beta,s_\alpha\beta\}$
and the result follows as in  the proof of 
\eqref{e:zetacxalphabeta}.
\end{proof}

The dependence of $\Gamma(H,\Phi^+)$ on $\Phi^+$ 
follows easily from Lemma \ref{l:elementary}(2) and its proof.

\begin{lemma}
\label{l:GammaPhi}
Suppose $\Phi^+_1,\Phi^+_2$ are special sets of positive roots for  $H$.
Let $(h, \tg) \in X(\oH, \tH)$, and write 
$h=\gamma h_0$ with $\gamma\in \Gamma(\oH)$ and $h_0\in\oH^0$. Then
\begin{equation}
\Gamma(H,\Phi^+_1)(h,\tg)
=
\Gamma(H,\Phi^+_2)(h,\tg)
e^{\rho(\Phi^+_2)-\rho(\Phi^+_1)}(h_0).
\end{equation}
\end{lemma}

The factors $\Gamma(H,\Phi^+)$ are conjugation invariant in the following
sense.

\begin{lemma}
\label{l:Gammaconj}
Suppose $\Phi^+$ is a special set of positive roots for $H$, and
$(h,\tg)\in X(\oH,\tH)$. Fix $\wt x\in \tG$ and let $x=\overline
p(p(\wt x))\in \oG$. Let $\oH_1=x\oH x\inv$ and
$\Phi^+_1=\Ad(x)\Phi^+$. Then
\begin{equation}
\Gamma(H_1,\Phi^+_1)(xhx\inv,\wt x\tg\wt x\inv)=\Gamma(H,\Phi^+)(h,\tg)
\end{equation}  
\end{lemma}

\begin{proof} 
Choose  $(\wt\chi, \chi )\in \caS(H,\Phi^+,\wt\chi_s,\chi_s)$. 
Let $\wt\chi_1=\wt x\tilde\chi\wt x\inv$ and $\chi_1=x\chi x\inv$. 
It is clearly enough to show
$(\wt\chi_1,\chi_1)\in\caS(H_1,\Phi^+_1,\wt\chi_s,\chi_s)$.
Condition \eqref{e:cond}(b) is immediate, and 
for all $h\in(\oH_1\cap \oG_d)^0$
$$
\chi_1(h)=\chi(x^{-1}tx)=\tilde\chi^2e^{\rho(\Phi^+)}(x^{-1}tx)=
\tilde\chi_1^2e^{\rho(\Phi^+_1)}(h)
$$
so \eqref{e:cond}(a) holds as well.

Assume $\gamma\in\Gamma(\oH _1)Z(\oG)$ and choose $y$ satisfying
$y\Gamma(\oH)y\inv\subset \Gamma(\oH_s)$.
By \eqref{e:compatible}
\begin{equation}
\begin{aligned}
\chi_1(\gamma)=\chi(x\inv \gamma x)&=(\wt\chi/\wt\chi_s)(\phi(x\inv\gamma x))
\chi_s(yx\inv\gamma xy\inv)\\
&=
(\wt\chi_1/\wt\chi_s)(\phi(\gamma))
\chi_s((yx\inv)\gamma (yx\inv)\inv).
\end{aligned}
\end{equation}
This proves \eqref{e:compatible} holds for $(\wt\chi_1,\chi_1)$.
\end{proof}

We summarize the important properties of $\Gamma(H,\Phi^+)$,
reformulating them slightly.

\begin{proposition}
\begin{subequations}
\renewcommand{\theequation}{\theparentequation)(\alph{equation}}  
Fix a split Cartan subgroup, a special set $\Phi^+_s$ of positive
roots,
and $(\wt \chi_s,\chi_s)\in\caS(H_s,\Phi^+_s)$.
Suppose $H$ is any Cartan subgroup, $\Phi^+$ is special, and
$(h,\tg)\in X(\oH,\tH)$.

\noindent (1) 
\begin{equation}
\Gamma(H,\Phi^+)(h,\tg)=e^{-\rho}(y)\tau(y,\tg)\quad (h\in\oH_d^0).
\end{equation}

\noindent (2) Suppose $\wt x\in\tG$, $x=\overline p(p(\wt x))$. Then
\begin{equation}
\Gamma(xHx\inv,\Ad(x)\Phi^+)(xhx\inv,\wt x\tg\wt x\inv)=\Gamma(H,\Phi^+)(h,\tg)
\end{equation}

\noindent (3) Suppose $(i=1,2)$ $(h_i,\tg_i)\in X(\oH,\tH)$. Then
\begin{equation}
\Gamma(H,\Phi^+)(h_1h_2,\tg_1\tg_2)=
\Gamma(H,\Phi^+)(h_1,\tg_1)\Gamma(H,\Phi^+)(h_2,\tg_2).
\end{equation}

\noindent (4) For $i=1,2$ let $H_i$ be a Cartan subgroup and
$\Phi^+_i$ a special set of positive roots for $H_i$.
Assume $(h,\tg)\in X(\oH_i,\tH_i)$ for $i=1,2$.
Choose $c\in G_{ad}(\C)$ satisfying $\Ad(c)\h_2=\h_1$,
$\Ad(c)|_{\h_1\cap \h_2}=1$, and $\Ad^*(c)(\Phi^+_1)=\Phi^+_2$. 
Then
\begin{equation}
\Gamma(H_1,\Phi^+_1)(h,\tg)=
\Gamma(H_2,\Phi^+_2)(h,\tg).
\end{equation}

Conditions (1,2,4) and Condition (3) for $h_1\in \Gamma(\oH)Z(\oG),\tg_1\in Z(\tG)$
uniquely determine the functions $\Gamma(H,\Phi^+)$,
given the restriction of $\Gamma(H_s,\Phi^+_s)$ to 
$(\Gamma(\oH_s)Z(\oG) \times Z(\tG))\cap X(\oH_s,\tH_s)$.
\end{subequations}
\end{proposition}

\sec{Transfer Factors}
\label{s:transfer}

We now define transfer factors for a general admissible triple
$(\tG, G,\oG )$, generalizing Definition \ref{d:transferspecial}. 

\begin{definition}
\label{d:transfer}
A set of lifting data for $(\tG,G,\oG)$ is a pair
$(\wt\chi_s,\chi_s)\in\caS(H_s)$ for some
maximally split Cartan subgroup $H_s$ of $G$.

The transfer factor $\Delta_{\oG}^{\tG}$ associated to
$(\tilde\chi_s,\chi_s)$ is the function on
$X'(\oG,\tG)$ (Definition \ref{d:XG}) defined as follows. 

Suppose $(h,\tg)\in X'(\oG,\tG)$. Let $H=\Cent_G(p(\tg))$; this is a Cartan
subgroup of $G$.
Choose a special set of positive roots $\Phi^+$ for $H$. 
Recall $\Gamma(H,\Phi^+)$ is given by Definition \ref{d:compatible},
and $\Delta^1$ is given in \eqref{e:Deltas}(c). Define the transfer factors

\begin{equation}
\label{e:transfer}
\Delta_{\oG}^{\tG}(h,\tg)=\frac{\Delta^1(\Phi^+,h)}{\Delta^1(\Phi^+,\tg)}
\Gamma(H,\Phi^+)(h,\tg).
\end{equation}
\end{definition}

\begin{proposition}
\label{p:independent}
$\Delta_{\oG}^{\tG}$ is well defined, i.e. 
independent of the choice of $\Phi^+$.
\end{proposition}

\begin{proof} 
We have to show that
$\frac{\Delta^1(\Phi^+,h)}{\Delta^1(\Phi^+,\tg)}$ 
and $\Gamma(H,\Phi^+)$ satisfy inverse
transformation properties with respect to $\Phi^+$.

Suppose $(h,\tg)\in X'(\oH,\tH)$.
By \eqref{e:Gamma(H)} write $h=\gamma h_0$ with $\gamma\in\Gamma(\oH)$
and $h_0\in \oH^0$.
Suppose $\Phi^+_1,\Phi^+_2$ are two choices of special positive
roots. 
\begin{subequations}
\renewcommand{\theequation}{\theparentequation)(\alph{equation}}  
By Lemma \ref{l:GammaPhi}
\begin{equation}
\Gamma(H,\Phi_1^+)(h,\tg)=
\Gamma(H,\Phi_2^+)(h,\tg)
e^{\rho(\Phi_2^+)-\rho(\Phi_1^+)}(h_0).
\end{equation}
By the definition of $\Delta^{\tG}_{\oG}$ it  is enough to prove 
\begin{equation}
\frac{\Delta ^1( \Phi  _1^+,h)}{\Delta ^1( \Phi  _1^+,\tg)} =
\frac{ \Delta ^1(\Phi _2 ^+,h)}{ \Delta ^1(\Phi _2 ^+,\tg)}
e^{\rho(\Phi_1^+)-\rho(\Phi_2 ^+)}(h_0).
\end{equation}  
\end{subequations}
Every term factors to $\oG$, so we can replace $\tg$ with $\overline
p(p(\tg))=\overline p(\phi(h))=h^2$.
This reduces us to showing
\begin{equation}
\label{e:enough}
\frac{\Delta ^1( \Phi  _1^+,h)}{\Delta ^1( \Phi  _1^+,h^2)} =
\frac{ \Delta ^1(\Phi _2 ^+,h)}{ \Delta ^1(\Phi _2 ^+,h^2)}
e^{\rho(\Phi_1^+)-\rho(\Phi_2 ^+)}(h_0).
\end{equation}  

\begin{subequations}
\renewcommand{\theequation}{\theparentequation)(\alph{equation}}  

We first compute for any set of positive roots $\Phi^+$:
\begin{equation}
\begin{aligned}
\frac{\epsilon_r(h,\Phi^+)}{\epsilon_r(h^2,\Phi^+)}
&=
\sgn\frac
{\prod_{\alpha\in\Phi^+_r}(1-e^{-\alpha}(h))}
{\prod_{\alpha\in\Phi^+_r}(1-e^{-\alpha}(h^2))}\\
&=
\sgn
\prod_{\alpha\in\Phi^+_r}(1+e^{-\alpha}(\gamma)e^{-\alpha}(h_0))\\
\end{aligned}
\end{equation}
A given term is positive unless $e^\alpha(\gamma)=-1$ and
$e^\alpha(h_0)<1$. 
Let 
\begin{equation}
\Phi^+_r(h)=\{\alpha\in\Phi_r\,|\, e^\alpha(h_0)>1\}.
\end{equation}
Then
\begin{equation}
\frac{\epsilon_r(h,\Phi^+)}{\epsilon_r(h^2,\Phi^+)}
=
\prod_{\alpha\in\Phi^+_r\cap(-\Phi^+_r(h))}e^\alpha(\gamma)
=e^{\rho_r(\Phi^+)-\rho (\Phi_r^+(h))}(\gamma).
\end{equation}

Letting $\Phi^+_{r,j}=\Phi^+_j\cap\Phi_r$
we have
\begin{equation}
\label{e:epr12}
\frac{\epsilon_r(h,  \Phi _1^+) } {\epsilon_r(h^2,  \Phi _1^+)} = e^{\rho (\Phi
_{r,1}^+) - \rho (\Phi _{r,2}^+) }(\gamma) \frac{\epsilon_r(h, \Phi _2^+)
} {\epsilon_r(h^2, \Phi _2^+)}.
\end{equation}

Now consider the $\Delta^0$ term of \eqref{e:Deltas}(c).
Write $\Phi _1 ^+ = w\Phi _2^+$ for some $w \in W(\Phi )$, so that
\begin{equation}
\begin{aligned}
\frac{\Delta^0(\Phi_1^+,h)}{\Delta^0(\Phi_1^+,h^2)}&=
\frac{\epsilon(w)e^{\rho(\Phi_2^+)-\rho(\Phi_1
^+)}(h)\Delta^0(\Phi_2^+,h)}{\epsilon(w)e^{\rho(\Phi_2^+)-
\rho(\Phi_1^+)}(h^2)\Delta^0(\Phi_2^+,h^2)}\\
&=
e^{\rho(\Phi^+_1)-\rho(\Phi^+_2)}(h)
\frac{\Delta^0(\Phi_2^+,h)}
{\Delta^0(\Phi_2^+,h^2)}
\end{aligned}
\end{equation}
since $\rho(\Phi^+_1)-\rho(\Phi^+_2)$ is a sum of roots.

Now
\begin{equation}
\begin{aligned}
e^{\rho(\Phi_{r,1}^+)-\rho(\Phi_{r,2}^+)}(\gamma)&
e^{\rho(\Phi_1^+)-\rho(\Phi_2^+)}(h)\\
&=e^{\rho(\Phi_{r,1}^+)-\rho(\Phi_{r,2}^+)-\rho(\Phi_2^+)+\rho(\Phi_1^+)}(\gamma)
e^{-\rho(\Phi_2^+)+\rho(\Phi_1^+)}(h_0)\\
&=e^{\rho(\Phi_1^+)-\rho(\Phi_2^+)}(h_0)
\end{aligned}
\end{equation}
since $\gamma$ has order $2$ and (by \eqref{e:zetacx2})
\begin{equation}
e^{\rho(\Phi_1^+)-\rho(\Phi_{r,1}^+)}(\gamma)
=\zeta_{cx}(\oG,\oH)(\gamma)
=e^{\rho(\Phi_2^+)-\rho(\Phi_{r,2}^+)}(\gamma).
\end{equation}
\end{subequations}
Multiplying (c) and (e), and using (f) shows
\eqref{e:enough}, and  completes the proof.
\end{proof}

Here are some elementary properties of transfer factors.

\begin{lemma}
Suppose $(h,\tg)\in X(\oH,\tH)$. Fix a special set of positive roots
$\Phi^+$ of $H$, and let $\rho=\rho(\Phi^+)$, $\lambda=\lambda(H,\Phi^+)$
(cf.~\ref{e:muHPhi+}).  Let $g=p(\tg)$. 
Then
\begin{subequations}
\label{e:D}
\renewcommand{\theequation}{\theparentequation)(\alph{equation}}  
\begin{equation}
\label{e:Da}
\Delta(h,\tg)^2=\frac{D(h)}{D(g)}e^{2\rho}(h)\lambda(\phi(h)).
\end{equation}
Write $h=th_0$ with $t\in \Gamma(\oH)Z(\oG)$ and $h_0\in\oH _d^0$. 
Then 
\begin{equation}
\label{e:Db}
\Delta(h,\tg)^2=\frac{D(h)}{D(g)}e^{2\rho}(t)\lambda(\phi(t)).
\end{equation}
In particular
\begin{equation}
\label{e:Dc}
\Delta(h,\tg)=c(h,\tg)\lambda(\phi(t))^{\frac12}\frac{|D(h)|^{\frac12}}{|D(\tg)|^{\frac12}}
\end{equation}
where $c(h,\tg)^4=1$.
Finally if $\lambda$ is unitary then
\begin{equation}
|\Delta(h,\tg)|=\frac{|D(h)|^{\frac12}}{|D(\tg)|^{\frac12}}.
\end{equation}
\end{subequations}
\end{lemma}

\begin{proof}
The main point is that
\begin{equation}
\left[
\frac{\Delta^0(\Phi^+,h)\epsilon_r(h,\Phi^+)}
{\Delta^0(\Phi^+,\tg)\epsilon_r(\tg,\Phi^+)}
\right]^2
=
\frac{\Delta^0(\Phi^+,h)^2e^{2\rho}(h)}
{\Delta^0(\Phi^+,\tg)^2e^{2\rho}(\tg)}
\frac{e^{2\rho}(\tg)}{e^{2\rho}(h)}
=
\frac{D(h)}{D(\tg)}e^{2\rho}(h),
\end{equation}
the last equality follows from \eqref{e:easy1} and the fact $\overline
p(p(\tg))=h^2$. Then \eqref{e:Da}  and \eqref{e:Db} follow from
\eqref{e:mu},
and (c) and (d) are elementary consequences of this.
\end{proof}

\begin{remark}
Recall $\mu\in\z(\C)^*$ is given in Definition \ref{d:compatible}.
If $\mu=0$ (i.e. $\lambda|_{Z(G)^0}=1$) then in \eqref{e:D}(c) $\lambda(\phi(t))^2=1$,
and $\Delta(h,\tg)$ differs from $|D(h)/D(g)|^{\frac12}$ by a fourth
root of unity. By varying our choice of lifting data
$(\tilde\chi_s,\chi_s)$ we are free to choose the restriction of
$\mu$ to $Z(G)^0$, up to the constraints \eqref{e:rho}. 
For simplicity assume $G=\oG$.
It is easy
to see we can take this restriction to be trivial if and only if
$e^\rho(z)=1$ for all $z\in Z(G)^0\cap H_d^0=Z(G)^0\cap Z(G_d)$.
This is equivalent to: $\rho$ exponentiates to a character of $H^0$. 
By the right hand side of the equality this is independent of the
Cartan, and we say $G$ is  {\it real-admissible} if it holds. 

This condition is empty unless $Z(G)$ contains a compact torus, and it
is weaker than admissibility. For example $GL(2,\C)$ is not
admissible; $GL(2,\R)$ is real-admissible, and $U(1,1)$ is not.

In any
event we can always choose the restriction of $\lambda$ to $Z_0(G)$ to be unitary, so that
\eqref{e:D}(d) holds.
\end{remark}

\begin{example}
\label{ex:U11_b}
Let $G=U(1,1)$ (see Example \ref{ex:U11_a}).
Then $G$ is not real-admissible, since $-I\in Z(G)^0$ and
$e^\rho(-I)=-1$. 
Let $T$ be the compact diagonal Cartan subgroup, and write
$X^*(T)=\Z^2$ in the usual coordinates. Let $\tG$ be an 
admissible cover of $G$. Recall (cf. Example~\ref{ex:U11_a}) $X^*(\wt
T)=\frac12\Z\oplus\Z$ or $\Z\oplus\frac12\Z$. 

Suppose $\wt\chi=(x,y)\in X^*(\wt T)$ and $\chi=(a,b)\in X^*(T)$. 
Then \eqref{e:cond1} holds 
if and only if
$a-2x-\frac12=b-2y+\frac12$, in which case
\begin{equation}
\mu=(a-2x-\frac12,a-2x-\frac12)\in\z(\C)^*.
\end{equation}
Since $a,2x\in\Z$ we cannot choose $\mu$ to be $0$.
  
\end{example}

The transfer factors satisfy the following invariance property with
respect to conjugation by $\tG$.

\begin{lemma}\label{invx} Suppose $(h, \tg)\in X'(\oH, \tH)$.   
Fix $\wt x\in\tG$  and let $x=\overline p(p(\wt x))\in \oG$.
Then
\begin{equation}
\Delta_{\oG}^{\tG}(xhx^{-1},\tilde x\tg\tilde x^{-1})
=\Delta_{\oG}^{\tG}(h,\tg). 
\end{equation}
\end{lemma}

Follows easily from Lemma \ref{l:Gammaconj}.

From Lemma \ref{l:simplytransitive} we see:

\begin{lemma}\label{choicedata} 
The group of characters of $\oG$ acts  simply transitively on the set of
transfer factors.

More precisely suppose for $i=1,2$
$(\wt\chi_{s}^i,\chi _{s}^i)\in\caS(H_s)$, and let
$\Delta_{\oG}^{\tG}(\chi_s^i,\tilde\chi_s^i)$ denote the corresponding transfer factors.
Then there is a character $\psi :
\oG \rightarrow \C ^{\times}$ so that for all $(h,\tg)\in X'(\oG,\tG)$
\begin{equation}
\label{choice12}
\Delta_{\oG}^{\tG}(\chi_s^1,\tilde\chi_s^1)(h,\tg)=\psi(h)\Delta_{\oG}^{\tG}(\chi_s^2,\tilde\chi_s^2)(h,\tg)
\end{equation} 
Conversely, suppose that $\psi : \oG \rightarrow \C
^{\times}$ is a character and $( \tilde \chi_s^1,\chi_s^1)\in\caS(H_s)$.
Then there exists $( \tilde \chi_s^2,\chi_s^2)\in\caS(H_s)$ 
such that (\ref{choice12}) is satisfied.
 \end{lemma}

Suppose $(h,\tg)\in X(\oH,\tH)$ and $e^{\rho}(h)$ is defined. This
happens, for example, if $\oG$ is admissible, or if $\oG_d$ is
admissible and   $h\in\oG_d$.
Then
\begin{equation}
\frac{\Delta^1(\Phi^+,h)\epsilon_r(h,\Phi^+)}
{\Delta^1(\Phi^+,\tg)\epsilon_r(\tg,\Phi^+)}
=
\frac{\Delta(h,\Phi^+)\epsilon_r(h,\Phi^+)}
{\Delta(\tg,\Phi^+)\epsilon_r(\tg,\Phi^+)}e^\rho(h).
\end{equation}
Together with \eqref{e:Gammaformulas} this gives simpler formulas for
the transfer factor in some situations.
For example \eqref{e:formulaforDelta} holds 
for $(h,\tg)$ if $\oG_d$ is acceptable and $h\in\oG_d$.

\begin{example}
\label{ex:gln2}
Let $G=\oG=GL(n,\R)$ and let $H$ be the diagonal split Cartan subgroup
of $G$.
Choose $(\tilde\chi,\chi)$  as in Example \ref{ex:gln}, and use other
notation as in that Example.
Suppose $(h,\tg)\in X'(H,\tH)$ 
and write $h=th_0\in \Gamma(H)H^0$.
An easy calculation gives
\begin{equation}
\Delta^1(\Phi^+,h)e^{\rho}(h_0)=|D(h)|^{\frac12}.
\end{equation}
From this and \eqref{e:gln} we compute
\begin{equation}
\Delta(h,\tg)=\frac{|D(h)|^{\frac12}}{|D(\tg)|^{\frac12}}\tau(h_0,\tg).
\end{equation}
This agrees with the transfer factors of \cite{ah}.
\end{example}

\begin{lemma}
\label{cc} 
Let $H$ be a Cartan subgroup of $G$.
Suppose $(h,\tg)\in X'(\oH,\tH)$.
Also suppose $z\in Z(\oG), \wt z\in Z(\tG)$ and $(z,\wt z)\in X(\oH,\tH)$.
Then $(zh,\wt z\tg)\in X'(\oH,\tH)$ and
\begin{equation}
\Delta _{\oG} ^{\tG} (zh,\wt z \tg) 
=\tilde\chi_s(\wt z)/\chi_s(z)\Delta_{\oG}^{\tG}(\tg,h).
\end{equation}
\end{lemma}

\begin{proof} 
Let $\Phi ^+$ be a special set of positive roots and choose
$(\wt\chi,\chi)\in \caS(H,\Phi^+,\wt\chi_s, \chi _s )$.
Since $z$ and $\wt z$ are central

$$
\Delta^1(\Phi^+,zh)=\Delta^1(\Phi^+,h),\,\Delta^1(\wt \tg,\Phi^+)=
\Delta^1(\Phi^+,\tg).
$$ 
Furthermore
\begin{equation}
\Gamma (H,\Phi ^+)(zh, \tilde z\tg) = 
\Gamma(H,\Phi^+)(z,\wt z)
\Gamma(H,\Phi^+)(h,\tg)
\end{equation}
and $\Gamma(H,\Phi^+)(z,\wt z)=\tilde\chi_s(\wt z)/\chi_s(z)$ by 
\eqref{e:tildechichi}.
Inserting this into the definition of $\Delta^{\tG}_{\oG}$ gives the result.
\end{proof}

\sec{Some Constants}
\label{s:constants}
We need to take care of some constants. Fix an admissible triple
$(\tG,G,\oG)$. Suppose $H$ is a Cartan subgroup of $G$, with
corresponding Cartan subgroups $\tH$ of $\tG$ and $\oH$ of $\oG$.
Let $\phi_{\oH}$ be the
restriction of $\phi$ to $\oH$.    Define

\begin{definition}
\label{d:constants}  

\begin{equation}
\label{e:constants}
\begin{aligned}
c(H)&=|\Ker(\phi_{\oH})||\tH/Z(\tH)|^{-\frac12}|Z(\tH)/p\inv(\phi(\oH))|^{-1}\\
&=|\Ker(\phi_{\oH})||H/Z_0(H)|^{-\frac12}|Z_0(H)/\phi(\oH)|^{-1}
\end{aligned}
\end{equation}
Fix a maximally split Cartan subgroup $H_s$ of $G$, and define
\begin{equation}
c=c(H_s),
\quad C(H)=c(H)/c(H_s).
\end{equation}
If it is necessary to specify $\oG$ we write $c_{\oG}(H)$, $c_{\oG}=c_{\oG}(H_s)$
and $C_{\oG}(H)=c_{\oG}(H)/c_{\oG}(H_s)$. 
\end{definition}

Let $H_2=\{h\in H\,|\, h^2=1\}$ and $H^0_2=H^0\cap
H_2$. Note that $|H^0_2|=2^{\dim(H\cap K)}$.
Let $H_s, H_f$ be maximally split and fundamental Cartan subgroups, respectively.

\begin{proposition}
\label{p:cC}  
\hfil\break
\noindent (1) $c_G(H)=|H^0_2||H/Z_0(H)|^{\frac12}$.

\medskip
\noindent (2) $c_{\oG}(H)=c_G(H)\frac{|\Gamma(H)\cap C|}{|C|}$.

\medskip
\noindent (3)  $c(H)$ and $C(H)$ are  integers, and
are powers of $2$.

\medskip
\noindent (4) $c=c(H_s)\le c(H)\le c(H_f)$.

\medskip
\noindent (5) $1\le C(H)\le C(H_f)$

\end{proposition}

\begin{proof}

By definition we have 
\begin{equation}
\label{e:constants2}
\begin{aligned}
c_{\oG}(H)
&=|\Ker(\phi_{\oH})||H/Z_0(H)|^{-\frac12}|Z_0(H)/\phi(\oH)|^{-1}\\
&=|\Ker(\phi_{\oH})||H/Z_0(H)|^{-\frac12}|H/\phi(\oH)|^{-1}|H/Z_0(H)|\\
&=|\Ker(\phi_{\oH})||H/Z_0(H)|^{\frac12}|H/\phi(\oH)|^{-1}.
\end{aligned}
\end{equation}
If $\oG=G$ then $\Ker(\phi_{\oH})=H_2$,
$\phi(\oH)=H^0$, and
\begin{equation}
c_G(H)
=|H_2||H/Z_0(H)|^{\frac12}|H/H^0|^{-1}.
\end{equation}
It is easy to see  $H\simeq H^0\times H/H^0$, and $H/H^0$ is a two-group,
which implies $|H_2|=|H^0_2||H/H^0|$.
Therefore $c_G(H)=|H^0_2||H/Z_0(H)|^{\frac12}$, which is (1).

\begin{subequations}
\renewcommand{\theequation}{\theparentequation)(\alph{equation}}  
By (1) and \eqref{e:constants2} 
\begin{equation}
\label{e:chh}
c_{\oG}(H)/c_G(H)=
\frac{|\Ker(\phi_{\oH})||Z_0(H)/\phi(\oH)|\inv}
{|H_2||Z_0(H)/H^0|\inv}
=
\frac{|\Ker(\phi_{\oH})||\phi(\oH)/H^0|}
{|H_2|}
\end{equation}

Recall (Lemma \ref{bs4}) $\phi(\oH)=(\Gamma(H)\cap C)H^0$. The
composition
of maps
$\oH\overset{\phi_{\oH}}\rightarrow\phi_{\oH}(\oH)\rightarrow \phi_{\oH}(\oH)/H^0$
gives an exact sequence
\begin{equation}
1\rightarrow\frac{\Ker(\phi_{\oH})}{\Ker(\phi_{\oH})\cap \oH^0}
\rightarrow
\frac{\oH}{\oH^0}
\rightarrow
\frac{\phi_{\oH}(\oH)}{H^0}
\rightarrow 1
\end{equation}

Solving the resulting identity
for $|\phi_{\oH}(\oH)/H^0|$, inserting in \eqref{e:chh} and
cancelling terms gives
\begin{equation}
c_{\oG}(H)/c_G(H)=
|\oH/\oH^0||\Ker(\phi_{\oH})\cap\oH^0|/|H_2|.
\end{equation}
It is easy to see that $\Ker(\phi_{\oH})\cap \oH^0=H^0_2/C\cap H^0$, so
the right hand side equals
\begin{equation}
\label{e:c(d)}
\frac{|\oH/\oH^0|}{|H_2/H^0_2||C\cap H^0|}
\end{equation}

On the other hand there is an exact sequence
\begin{equation}
1\rightarrow
\frac C{C\cap H^0}
\rightarrow
\frac H{H^0}
\rightarrow
\frac{\oH}{\oH^0}
\rightarrow
\Gamma(H)\cap C
\rightarrow 1.
\end{equation}
The  map  to $\Gamma(H)\cap C$ 
is given by $h\rightarrow y\inv\sigma(y)$ where $y\in H(\C)$
satisfies $\op (y)=h$.
This gives
\begin{equation}
\frac{|\oH/\oH^0|}{|H/H^0|}=\frac{|\Gamma(H)\cap C||C\cap H^0|}{|C|}.
\end{equation}
Using the fact that $H/H^0\simeq H_2/H^0_2$ 
and plugging this into \eqref{e:c(d)} gives (2).
\end{subequations}

We will see in Lemma \ref{adams3} that $|H/Z_0(H)|^{\frac12}$ is an
integer.    Further, $H/Z_0(H)$ is a quotient of $H/H^0$, which is a two
group, so by (1) 
$c_G(H)$ is an integer and a power of $2$.
The general case of (3) follows from (1) and (2) if we can show
\begin{equation}
\label{e:Z}
|H^0_2||\Gamma(H)\cap C|/|C|\in \Z.
\end{equation}
To see this note that $C\subset \Gamma(H)H^0_2$.
For $g\in C$ write $g=\gamma h$ accordingly, and 
define $\psi(\gamma h)=h(H^0_2\cap C)$. This is a well defined
homomorphism from $C$ to $H^0_2/H^0_2\cap C$.   It factors to an
inclusion
$C/C\cap \Gamma (H) \hookrightarrow H^0_2/H^0_2\cap C$, and \eqref{e:Z}
follows from this. This proves (3) for $c_{\oG}(H)$.

It is enough to prove (4); assertion (3) for $C_{\oG}(H)$ and (5)
follow easily. For this we use  Cayley transforms. Suppose $H_\alpha,H_\beta$ are as
in Section \ref{s:cartan}.
It is enough to show
\begin{equation}
c_{\oG}(H_\alpha)\le c_{\oG}(H_\beta).
\end{equation}

First assume $\oG=G$. Write $H_\alpha=T_\alpha A_\alpha$,
$H_\beta=T_\beta A_\beta$. Then $\dim(T_\beta)=\dim(T_\alpha)+1$, and
by part (1) of the Lemma it is enough to show
\begin{equation}
\label{e:half}
|H_\alpha/Z_0(H_\alpha)|^{\frac12}\le 2 |H_\beta/Z_0(H_\beta)|^{\frac12}.
\end{equation}

If $\alpha$ is of type $II$ choose $t\in H_\alpha$ satisfying
$\alpha(t)=-1$; otherwise let $t=1$. 
By \eqref{e:HaHb} write $H_\alpha=\langle (H_\alpha\cap
H_\beta)B_\alpha ,t\rangle$. By \eqref{e:zha}  
\begin{equation}
\label{e:halphazhalpha}
\begin{aligned}
\frac{H_\alpha}{Z_0(H_\alpha)}&\simeq
\frac{\langle (H_\alpha\cap H_\beta)B_\alpha,t\rangle}
{[Z_0(H_\alpha)\cap Z_0(H_\beta)]B_\alpha}\\
&\simeq
\frac{\langle H_\alpha\cap H_\beta,t\rangle}
{[Z_0(H_\alpha)\cap Z_0(H_\beta)][\langle H_\alpha \cap  H_\beta ,t \rangle \cap B_\alpha]}\\
&\simeq
\frac{\langle H_\alpha\cap H_\beta,t\rangle}
{Z_0(H_\alpha)\cap Z_0(H_\beta)}
\end{aligned}
\end{equation}
since $\langle H_\alpha \cap  H_\beta ,t \rangle \cap B_\alpha \subset \langle H_\beta ,t \rangle\cap
B_\alpha=1$. 
Similarly
\begin{equation}
\label{e:hbetazhbeta}
\begin{aligned}
\frac{H_\beta}{Z_0(H_\beta)}&\simeq
\frac{(H_\alpha\cap H_\beta)B_\beta}
{[Z_0(H_\alpha)\cap Z_0(H_\beta)]B_\beta}\\
&\simeq
\frac{ H_\alpha\cap H_\beta}
{[Z_0(H_\alpha)\cap Z_0(H_\beta)][H_\alpha\cap H_\beta\cap B_\beta]}\\
&\simeq
\frac{ H_\alpha\cap H_\beta}
{\langle Z_0(H_\alpha)\cap Z_0(H_\beta),m_\alpha\rangle}
\end{aligned}
\end{equation}
We have $m_\alpha\in H_\beta^0\subset Z_0(H_\beta)$. By \eqref{e:tau2} 
$m_\alpha\in Z_0(H_\alpha)$ if and only if $\alpha$ is of type $I$.
By \eqref{e:halphazhalpha} and \eqref{e:hbetazhbeta} we see 
$|H_\alpha/Z_0(H_\alpha)|=\tau|H_\beta/Z_0(H_\beta)|$ where $\tau=1$
(resp. $4$) if $\alpha$ is of type $I$ (resp. II). This proves
\eqref{e:half}.

Now supppose $\oG$ is not equal to $G$.
By (2) of the Lemma 
\begin{equation}
\frac{c_{\oG}(H_\alpha)}{c_{\oG}(H_\beta)}
=
\frac{c_G(H_\alpha)}{c_G(H_\beta)}
\frac{|\Gamma(H_\alpha)\cap C|}{|\Gamma(H_\beta)\cap C|}
\end{equation}
We need to show  the right hand side is $\le 1$.

By the preceding argument the first quotient on the right hand side is
equal to $1$ if $\alpha$ is of type II for $H_\alpha$, or $\frac12$
otherwise.
It is clear that $\Gamma(H_\beta)\subset \Gamma(H_\alpha)$ of index
$1$ or $2$, so
$\Gamma(H_\beta)\cap C$ is of index $1$ or $2$ in $\Gamma(H_\alpha)\cap C$.
We need to show that if $\alpha$ is of type II for $H_\alpha$ then
$\Gamma(H_\alpha)\cap C=\Gamma(H_\beta)\cap C$.

By the assumption on $\alpha$ we can choose $t\in \Gamma(H_\alpha)$
such that $\alpha(t)=-1$.  Since $\alpha(g)=1$ for all $g\in \Gamma(H_\beta)$,
$t\not\in \Gamma(H_\beta)$, so
$\Gamma(H_\alpha)=\langle \Gamma(H_\beta),t\rangle$.
Therefore $\alpha(g)=-1$ for all $g\in \Gamma(H_\beta)t$, so 
$\Gamma(H_\alpha)\cap C=\Gamma(H_\beta)\cap C$.

\end{proof}

In fact we have shown that
\begin{equation}
\frac{c_{\oG}(H_\alpha)}{c_{\oG}(H_\beta)}=
\begin{cases}
1&\alpha\text{ type II for }H_\alpha\\  
1\text{ or }\frac12&\alpha\text{ type I for }H_\alpha
\end{cases}
\end{equation}
(Only $1/2$ occurs in the second case if $C=1$).
This implies that for any Cartan subgroup $H$
\begin{equation}
\label{e:Cestimate}
C_{\oG}(H)\le2^{r(H_s)-r(H)}
\end{equation}
where $r(*)$ denotes split rank.

\begin{example}
If $G=GL(n,\R)$ every real root $\alpha$ is of type
II, so $c_G(H)$ is independent of $H$ and
$C_{G}(H)=1$ for all $H$. See \cite{ah}.

On the other hand if $G=\oG=U(p,q)$ then every real root is of type $I$,
and equality holds in \eqref{e:Cestimate}. 
\end{example}

\sec{Lifting: Definition and Basic Properties}
\label{s:lifting}
Assume that $(\tG, G, \oG )$ is an admissible triple.  
In this section we  define lifting from $\oG$ to $\tG$ and derive
some of its basic properties.
Fix lifting data $(\wt\chi_s,\chi_s)$ 
with corresponding transfer factors
$\Delta_{\oG}^{\tG}$ 
(Definition \ref{d:transfer}).

Recall (Section \ref{s:admissible}) $\O(G,g)$ is the conjugacy
class of $g\in G$, and 
\begin{equation}
\O^{\st}(G,g)=\{xgx\inv\,|\, x\in G(\C)\}\cap G
\end{equation}
is the stable orbit. We also have  $\O(\oG,g)$ and
$\O^{\st}(\oG,g)$, and the map (Lemma \ref{l:phi}(3))
$\phi:\Orbst(\oG)\rightarrow \Orbst(G)$. 

\begin{definition}
\label{d:relevant}
Suppose $\tg$ is a strongly regular semisimple element of $\tG$.
Let $\tH=p\inv(\Cent_{G}(p(\tg)))$. We say $\tg$ and $\O(\tG, \tg)$ are {\it relevant} if
$\tg\in Z(\tH)$.
\end{definition}

The reason for this terminology is:

\begin{lemma}[\cite{adams_japan}, Proposition 2.7]
\label{l:relevant}
Let $\wt\pi$ be a genuine  admissible representation of $\tG$.
Suppose $\tg$ is a strongly regular semisimple element which is not
relevant. Then $\Theta_{\wt\pi}(\tg)=0$.
\end{lemma}
This is elementary, and goes back to  $GL(2)$
\cite{flicker_gl2}. 

Here are a few basic properties about orbits.

\begin{lemma}
\label{l:basicorbits}
Suppose $\tg$ is a semisimple element of $\tG$.
Then
\begin{enumerate}
\item $p(\O(\tG,\tg))=\O(G,p(\tg))$
\item $p\inv(\O(G,p(\tg)))=\O(\tG,\tg)\cup \O(\tG,-\tg)$
\item Suppose $\tg$ is  relevant and strongly regular. Then $\O(\tG,\tg)\ne \O(\tG,-\tg)$.
\end{enumerate}
\end{lemma}

\begin{proof} Part (1) and (2) are routine.
For (3) suppose not.
Then $\wt x\tg\wt x\inv=-\tg$ for some $\wt x\in
\tG$. Let $H=\Cent_{G}(p(\tg))$.
Then $p(\wt x)p(\tg)p(\wt x\inv)=p(\tg)$, so $p(\wt x)\in H$. But then
$\wt x\tg\wt x\inv=\wt g$ since $\tg\in Z(\tH)$ and $\wt x\in \tH$.
\end{proof}

\begin{lemma}
\label{l:fordef}
Suppose $\overline{\O}^{\st}\in\Orbst(\oG)$ and write
$\phi(\overline{\O}^{st})=\Ost(G,g)$ for some $g\in G$.
Assume $g$ is strongly regular.
Then there is a unique $h \in \overline\O^{\st}$ such that $\phi (h) =
g$.
Furthermore  $h \in \oH$ and is also strongly regular.     \end{lemma}

\begin{proof}  By definition there are $h' \in \overline\O^{\st}$ and
$g' \in \O^{\st}(G,g)$ such that $\phi (h') = g'$.
Let $x \in G(\C )$ such that $g' = xg x^{-1}$.
Since $g, g' \in G$ we have $\sigma (x)^{-1}x \in \Cent_{G(\C)}(g)=H(\C)$.    
Let $y = \op (x)\in \oG(\C )$ and $h = y^{-1}h'y$.
Since $g = \phi (h) = s(h)^2$, 
$s(h) \in \Cent_{G(\C)}(g) = H(\C)$.
But $\sigma (y)^{-1}y \in \oH (\C)$  so $\sigma (h) = h$ and
$h \in \oH \cap \overline\O^{\st}$ with $\phi (h) = g$.   

Suppose $r \in \Cent_{\oG (\C)}(h)$, and choose a preimage $s$ of $r$
in $G(\C )$.  Then $g = \phi (h) = \phi (rhr^{-1}) = s\phi (h)s^{-1} =
sgs^{-1}$ so that $s \in H(\C)$.  Thus $r \in \oH (\C)$, so that $h$
is also strongly regular.

Suppose that $h_1 \in \overline\O^{\st}$ with $\phi (h_1) = g$.  Then
there is $u \in \oG (\C)$ with $uhu^{-1} = h_1.$ Choose a preimage $v$
of $u$ in $G(\C )$.  Then $g = \phi (h_1) = \phi (uhu^{-1} ) = v\phi
(h)v^{-1} = vgv^{-1}$.  Thus $v \in H(\C)$.  Now $u \in \oH (\C) $ so
$h_1=h$.
\end{proof}

\begin{remark}  Note that Lemma \ref{l:fordef} implies the following.
Let $\O^{st} \in \Orbst(G)$ and $\overline\O^{\st} \in \Orbst(\oG)$
with $\phi (\overline\O^{\st}) = \O^{st}$.
If $\O^{st}$ is strongly regular and semisimple,
then so is $\overline\O^{\st}$, and the restriction of $\phi $ to
$\overline\O^{\st}$
gives a bijection between $\overline\O^{\st}$ and $\O^{\st}$.
Further, for any $g \in \O^{\st} \cap H$,
\begin{equation}
\label{forlift}
\{\overline\O^{\st}\,|\,
\phi(\overline\O^{\st})=\O^{\st}\}
=
\{\Orbst(\oG,h)\,|\, h \in \oH, \phi(h)=g\}.
\end{equation}
\end{remark}

 \begin{definition}
Suppose $\wt\O\in \Orb(\tG)$ and $\overline\O^{\st}\in \Orbst(\oG)$
are strongly regular and semisimple.
Let $p(\wt\O)^{\st}$ be the stable orbit for $G$ containing $p(\wt\O)$.

Define 
$\Delta^{\tG}_{\oG}(\overline\O^{\st},\wt\O)=0$
unless 
$\phi(\overline\O^{\st})=p(\wt\O)^{\st}$.

Suppose $\phi(\overline\O^{\st})=p(\wt\O)^{\st}$. Choose $\tg\in
\wt\O$.   
By Lemma \ref{l:fordef} there is a unique 
$h\in\overline\O^{\st}$ with
$\phi(h)=p(\tg)$. Define
\begin{equation}
\Delta^{\tG}_{\oG}(\overline\O^{\st},\wt\O)=
\Delta^{\tG}_{\oG}(h,\tg).
\end{equation}
By Lemma \ref{invx} this is independent of the choice 
of $\tg$.
\end{definition}

Recall a virtual character $\pi$ of $\oG$ is said to be {\it stable}
if the function $\theta_\pi$ representing its character is constant on
the stable orbit   $\O^{st}(G,g)$ for  all 
strongly regular semisimple elements $g$. 

\begin{definition}
\label{d:lift}
Suppose $\overline\Theta$ is a stable virtual character of $\oG$.
Then $\overline
\Theta(\overline\O^{\st})$ is defined for any strongly regular
semisimple orbit $\overline\O^{\st}\in\Orbst(\oG)$.
Suppose $\wt\O\in\Orb(\tG)$ is strongly regular and semisimple.
Recall $c=c_{\oG}$ is given by Definition \ref{d:constants}.

Define
\begin{subequations}
\renewcommand{\theequation}{\theparentequation)(\alph{equation}}  
\label{e:liftorbits} 
\begin{equation}
\Lift_{\oG}^{\tG}(\overline\Theta)(\wt\O)=
c\inv
\sum_{\{\overline\O^{\st}\in \Orbst(\oG)\}}
\Delta^{\tG}_{\oG}(\overline\O^{\st},\wt\O)\overline\Theta(\overline\O^{\st}).
\end{equation}
This is a finite sum, over 
\begin{equation}
\{\overline\O^{\st}\,|\, \phi(\overline\O^{\st})=p(\wt\O)^{\st}\}.
\end{equation}
\end{subequations}
\end{definition}

Thus ${\Lift}_{\oG} ^{\tG}\overline\Theta$ 
is a genuine class function defined on the set of strongly regular
semisimple orbits in $\tG$. The following lemma follows easily from
(\ref{forlift}) and the fact 
that $\phi$ 
restricted to $\oH$ is a homomorphism. Recall $X(\oH,\tg)$ is given
by Definition \ref{d:XG}.

\begin{lemma}
Suppose $\tg$ is a strongly regular semisimple element.
Let $H=\Cent_{G}(p(\tg))$ and $\oH=\op(H)$.
Then
\begin{equation}
\label{e:lift1}
\Lift_{\oG}^{\tG}(\overline\Theta)(\tg)
=c\inv\sum_{\{h\in X(\oH,\tg)\}}
\Delta^{\tG}_{\oG}(h,\tg)\overline\Theta(h).
\end{equation}
In particular $\Lift_{\oG}^{\tG}(\overline\Theta)(\tg)=0$ unless
$p(\tg)$ is in the image of $\phi$. Assume this holds, and choose $h$
satisfying $\phi(h)=p(\tg)$. Then
\begin{equation}
\label{e:lift2}
\Lift_{\oG}^{\tG}(\overline\Theta)(\tg)
=c\inv\sum_{\{t\in\oH\,|\, \phi(t)=1\}}
\Delta^{\tG}_{\oG}(th,\tg)\overline\Theta(th).
\end{equation}
\end{lemma}

\begin{remark}
Formula \eqref{e:lift1} can be used to extend
$\Lift_{\oG}^{\tG}(\Theta)$ to a genuine class function on all 
regular (not just strongly regular) semisimple elements.
\end{remark} 

We derive some elementary properties of lifting. Recall the
definitions of \eqref{e:Deltas}.

\begin{lemma}\label{supptt} ${\rm Lift}_{\oG} ^{\tG}(\overline\Theta)$ is
supported on $Z(\tG)\tG_d^0$.
\end{lemma}

\begin{proof} Let $\oH$ be a Cartan subgroup  of $\oG$.  Then $\oH =
Z(\oG)\Gamma (\oH) \oH _d^0$, so that $\phi (\oH ) \subset
Z_0(G)H_d^0$.  Thus $p^{-1}\phi (\oH ) \subset Z(\tG)\tilde
G_d^0$.
\end{proof}

Define $Z_1({\oG} ) = \{ z \in Z({\oG} ): \phi (z) =1 \}$.

\begin{lemma}\label{cl0} 
Let $\overline\Theta $ be a stable character of ${\oG}$ with central
character $\zeta$.  Then $\Lift_{\oG} ^{\tilde G}(\overline \Theta)
\equiv 0$ unless $\zeta(z)=\chi_s(z)$ for all $z\in Z_1({\oG})$.
\end{lemma}

\begin{proof} Fix a Cartan subgroup $H$, $\tg \in \tH'$,
and $z \in Z_1({\oG} )$.  Then
$$\{ h \in \oH \,|\,\phi (h) =  p(\tg)\} = \{ zh\,|\, h   \in \oH  , \phi (h) =  p(\tg)\}.$$
Further, for all $h \in \oH $ such that $\phi (h) = \tilde p(\tilde
g)$,
$$\Delta_{\oG} ^{\tG} (\tg, zh) = \chi _s(z^{-1})\Delta _{\oG} ^{\tG}
(\tg, h)$$ by Lemma \ref{cc} applied to $h, \tg, z$, and
$\tilde z=1$.  We have
$$ 
Lift _{\oG} ^{\tG}(\overline\Theta)(\tg) =c\inv
\sum_{X(\oH,\tg)}\Delta _{\oG} ^{\tG}
(\tg, zh) \overline\Theta (zh) =c\inv\zeta (z)\chi _s(z^{-1}) Lift
_{\oG}^{\tG}(\overline\Theta)(\tg).
$$
\end{proof}

Fix a character $\zeta$ of $Z(\oG)$ satisfying
$\zeta (z) = \chi _s(z)$ for all $z\in Z_1({\oG})$. 
For $\wt z\in \wt S = p^{-1}\phi (Z(\oG))$ define
\begin{equation}
\label{e:zeta}
\wt\zeta(\wt z)=\wt\chi_s(\wt z)(\zeta\chi_s\inv)(z)
\end{equation}
where $z\in Z(\oG)$, $\phi(z)=p(\wt z)$.
This is independent of the choice of $z$, and is a genuine character
of $\wt S$.

\begin{lemma}\label{ccc} Let $\overline\Theta $ be a stable character of ${\oG}
$ with central character $\zeta$ such that $\zeta(z)=\chi _s(z)$ for
all 
$z\in Z_1({\oG} )$.  Then for all $\tg \in \tG, \tilde z \in
\tilde S$,
$$\Lift _{\oG} ^{\tG}(\overline\Theta)(\tilde z \tg) =  \tilde \zeta (\tilde z) Lift _{\oG} ^{\tG}
(\overline\Theta)(\tg).$$
\end{lemma}

\begin{proof} 
Fix a Cartan subgroup $H$ of $G$, $\tg \in \tH'$ and 
$\tilde z \in \tilde S$. Choose $z\in Z ({\oG})$ such that $\phi (z) =
p(\tilde z)$.  
Then 
$$
\{h\in\oH\,|\,(h,\wt z\tg)\in X(\oH,\tH)\} = 
\{zh\,|\,h\in\oH,(h,\tg)\in X(\oH,\tH)\}.
$$ 
Thus by Lemma \ref{cc},
$$
\begin{aligned}
Lift _{\oG}^{\tG}(\overline\Theta)(\tilde z \tg) &=   
c\inv\sum _{X(\oH,\tg)}
\Delta _{\oG }^{\tG}(\tilde z \tg, zh)(\overline\Theta  (zh)\\
&=(\tilde\chi_s(\wt z)/\chi_s(z))\zeta(z)Lift _{\oG}^{\tilde
G}(\overline\Theta)(\tg) = \tilde\zeta (\tilde z) Lift _{\oG}^{\tG}\overline\Theta)(\tg).
\end{aligned}
$$
\end{proof}

\sec{Lifting for Tori}
\label{s:tori}

We give some details of lifting for tori. This illustrates some of the
basic principles, and plays an important role in lifting for general reductive groups.

Let $G(\C)$ be an algebraic torus, with real points $G$.
Let $p:\tG\rightarrow G$ be any two-fold cover; any such cover is
necessarily admissible (Definition \ref{defadm}). Let $C$ be a  subgroup of $p(Z(\tG))$
consisting of elements of order $2$, and let $\oG=G(\C)/C$, with real
points $\oG$.
Then $(\tG,G,\oG)$ is an admissible triple. 

Since $G$ is a torus, $G=H=H_s$ is a maximally split Cartan subgroup.
Let $X_g(Z(\tG))$ denote the genuine characters of $Z(\tG)$ and
let $X(\oG)$ denote the characters of $\oG$.
Since $\g_d=0$ condition \eqref{e:cond} is empty, so
$(\wt\chi,\chi)\subset \caS(G)$ for all $\wt\chi\in X_g(Z(\tG)),\chi\in
X(G)$.

The transfer factor $\Delta_{\oG}^{\tG}$ of
Definition \ref{d:transfer} is then given by
$$
\Delta_{\oG}^{\tG}(h,\tg)=\wt\chi(\tg)\chi(h)\inv\quad((h,\tg)\in X(\oG,\tG)).
$$
Let 
\begin{equation}
\label{e:S}
\tS=\{\tg\in\tG\,|\, p(\tg)\in \phi(\oG)\}=\{\tg\,|\, X(\oG, \tg)\ne\emptyset\}.
\end{equation}
If $\psi\in X(\oG)$ then

\begin{equation}
\label{abelian1} 
Lift _{\oG}^{\tG}(\psi)(\tg)= 
c\inv\wt\chi(\tilde g)\sum _{h\in X(\oG,\tg)} \chi ^{-1}(h)\psi (h),
\end{equation} 
with $c\in\Z$  given by  Definition \ref{d:constants}.
In particular $\Lift_{\oG}^{\tG}(\psi)(\tg)=0$ for $\tg\not\in\tS$.

Suppose $\tg\in\tS$ and fix $h\in X(\oG, \tg)$. Then
\begin{equation}
\label{abelian2} 
Lift _{\oG}^{\tG}\psi (\tilde
g) = c\inv \tilde \chi (\tg) \chi ^{-1}(h)\psi (h) \sum _{ \{ t \in
\oG\,|\, \phi (t) =
 1 \}} \chi ^{-1}(t)\psi (t) .
\end{equation} 

Thus $\Lift
_{\oG}^{\tG}\psi =0$ unless $\chi (t) = \psi (t)$ for all $t \in\Ker(\phi)$.
Assume this holds. For $\tg\in\tS$ choose $h\in X(\oG, \tg)$ and define
\begin{equation}
\label{e:psi0}
\wt\psi_0(\tg)=\wt\chi(\tg)\psi(h)\chi(h)\inv.
\end{equation}
This is independent of the choice of $h\in X(\oG,\tg)$, so  $\wt\psi_0$ is
a well defined  genuine character of $\tS$, and it follows immediately that

\begin{equation}
\label{abelian3} 
Lift _{\oG}^{\tG}(\psi)(\tg)=
\begin{cases}
0&\tg\not\in\tS\\
c\inv|\Ker(\phi)|\tilde \psi _0(\tg)&\tg \in \tS.
\end{cases}
\end{equation}
It follows easily from the induced character formula 
that
\begin{equation}
\Lift_{\oG}^{\tG}(\psi)=c\inv|\Ker(\phi)||\tG/\wt S|^{-1}\,\Ind_{\tS}^{\tG}(\wt\psi_0).
\end{equation}
By induction by stages we have
\begin{equation}
\Ind_{\tS}^{\tG}(\wt\psi_0)=
\Ind_{Z(\tG)}^{\tG}\Ind_{\tS}^{Z(\tG)}(\wt\psi_0)
\end{equation}
In order to identify this 
as a sum of characters of $\tG$, we use the following elementary
result.

\begin{lemma}[\cite{shimura}, Proposition 2.2]
\label{adams3} 
Let $\tilde \psi \in X_g(Z(\tG))$.  Then
there is a unique irreducible genuine representation $\tilde \tau =
\tilde \tau (\tilde \psi) $ of $\tG$ such that $\tilde \tau
|_{Z(\tG)}$ is a multiple of $\tilde \psi $.  The map $\tilde
\psi \rightarrow \tilde \tau
 (\tilde \psi )$ is a bijection between
$X_g(Z(\tG))$ and the set of equivalence classes of irreducible
genuine representations of $\tG$, and
 $$
Ind_{Z(\tG) }^{\tG} \tilde \psi = |\tG/Z(\tG)|^{\frac12} \tilde \tau (\tilde \psi) .
$$
    The dimension of $\wt\tau(\wt\psi)$ is $|\tG/Z(\tG)|^{\frac12}$; in
particular this is an integer.    
\end{lemma}

Let 
\begin{equation}
\label{e:Xg}
X_g(Z(\tG),\wt\psi_0)=\{\wt\psi\in X_g(Z(\tG))\,|\, \wt\psi|_{\tS}=\wt\psi_0\}.
\end{equation}
By Frobenius reciprocity we have
\begin{equation}
\label{e:toruslift}
\Lift_{\oG}^{\tG}(\psi)=c\inv|\Ker(\phi)||\tG/\wt S|^{-1}
|\tG/Z(\tG)|^{\frac 12}\sum_{\wt\psi\in X_g(Z(\tG),\wt\psi_0)}\Theta _{\tilde \tau (\tilde \psi) }.
\end{equation}
By (\ref{e:constants}) the term in front of the sum is equal to $1$
(of course we 
defined $c$ precisely to make this hold).

We summarize these results.  As defined, $\Lift_{\oG}^{\tG}(\psi)$ is
a class function on $\tG$, which by \
\eqref{e:toruslift}) 
the character of a representation. 
We identify $\Lift_{\oG}^{\tG}(\psi)$ with this representation.

\begin{proposition}
\label{p:tori}
Fix $\wt\chi\in X_g(Z(\tG))$ and $\chi\in X(\oG)$, and use them to define
$\Lift_{\oG}^{\tG}$. Fix $\psi\in X(\oG)$.
Then $\Lift_{\oG}^{\tG}(\psi)\ne 0$ if and only if $\chi (t) = \psi (t)$ for all $t
\in\Ker(\phi)$.
Assume this holds. Define $\wt\psi_0$ as in \eqref{e:psi0} and
$X_g(Z(\tG),\wt\psi_0)$ as in \eqref{e:Xg}. Then
\begin{equation}
\Lift_{\oG}^{\tG}(\psi)
=\sum_{\wt\psi\in X_g(Z(\tG),\wt\psi_0)} \tilde \tau (\tilde \psi)  .
\end{equation}
This is an identity of genuine representations of $\tG$; 
the right hand side is the direct sum of 
$|p(Z(\tG))/\phi(\oG)|$ irreducible representations.
The differentials satisfy
\begin{equation}
\label{e:dtpsi}
d\wt\psi=\frac12(d\psi-d\mu)
\end{equation}
where $d\mu=d\chi-2d\wt\chi\in\g(\C)^*$ as in \eqref{e:mu}.

In particular assume $C$ is chosen so that
$\phi(\oG)=p(Z(\tG))$. Then $X_g(Z(\tG),\wt\psi_0)$ consists of the 
single character $\wt\psi_0$, and
\begin{equation}
\label{e:goodc}
\Lift_{\oG}^{\tG}(\psi)= \tilde \tau (\tilde \psi_0) .
\end{equation}
\end{proposition}

If it is necessary to indicate the dependence of $\Lift_{\oG}^{\tG}$
on $\wt\chi$ and $\chi$ we will write
\begin{equation}
\label{e:specify}
\Lift_{\oG}^{\tG}(\wt\chi,\chi,\psi)=
\Lift_{\oG}^{\tG}(\psi).
\end{equation}

\begin{corollary}
\label{gammaexists}      
Fix $(\wt\chi,\chi)$ as above, and let $\tilde\psi$ be a 
genuine character  of $Z(\tG)$. Then there is a unique character $\psi $ of $G$ such that
$\tilde \tau (\tilde \psi)$  occurs in $\Lift_{\oG}^{\tG}(\wt\chi,\chi,\psi)$.
It is given by
\begin{equation} 
\psi (h) = \chi (h) (\tilde \psi \tilde \chi^{-1})(\phi (h))\quad(h\in\oG).
\end{equation} 
\end{corollary}

\begin{example} 
\label{ex:ctori}   
Let $G(\C)$ be an algebraic torus with connected real points $G$ and
let $\tilde G$ be a two-fold cover.  Then $\tilde G = Z(\tG)$ is also
abelian.  We may as well take $\oG=G$. Let $\tilde\chi$ be any genuine
character of $\tG$, and let $\chi=\tilde\chi^2$. 
Suppose $\psi \in X(G)$.  Using Proposition \ref{p:tori}
$\Lift_G^{\tG}(\psi)\not = 0$ if and only if (\ref{e:psi0})  gives a 
well-defined genuine character $\tilde \psi $ of $\tG$.  It is easy to
see that in this case $\tilde \psi $ is the unique genuine character
with $\psi = \tilde \psi ^2$.

Therefore, with this choice of $(\tilde\chi,\chi)$, $Lift _G^{\tG}(\psi)\not = 0$
if and only if $\psi $ has a genuine square root $\tilde \psi $ on
$\tG$, in which case $Lift _G^{\tG}(\psi)= \tilde \psi $.

\end{example}

\sec{Example: Minimal Principal Series of Split Groups}
\label{s:ps}

Now assume 
$(\tG, G, \oG )$ is an
admissible triple where $G$ is split.  
Let $H_s$ be a split
Cartan subgroup of $G$, let $\Phi$ denote the roots of $H_s$,
and let $W=W(\Phi)$.
Fix a set of positive roots $\Phi^+$ and choose 
$(\wt\chi_s,\chi_s)\in \caS (G,H_s,\Phi^+)$
(cf. \ref{e:ST}).
Let $\rho=\rho(\Phi^+)$ and define
$\Delta_{\oG}^{\tG}$ as in Definition \ref{d:transfer} accordingly.

Corresponding to each character $\psi$ of $\oH_s$
is a principal series representation $\pi(\psi)$ of $\oG$. Its 
character
$\Theta_{\pi(\psi)}$ is 
supported on the conjugates of $\oH _s$ and satisfies

\begin{equation}
\label{split1} 
\Theta_{\pi(\psi)}(h)=
|D_{\oG}(h)|^{-\frac12} 
\sum _{w \in W} \psi (wh)\quad(h \in \oH_s').
\end{equation} 
We wish to compute $\Lift_{\oG}^{\tG}(\pi(\psi))$.  Using
(\ref{e:lift1}) it is clear that
$\Lift_{\oG}^{\tG}(\Theta_{\pi(\psi)})(\tg) =0$ unless $\tg$ is
conjugate to an element of $\tH_s$.

Suppose $\tg \in \tH_s'$. Then by \eqref{e:lift1}
\begin{equation}
\label{split2} 
\Lift_{\oG}^{\tG}(\Theta_{\pi(\psi)})(\tg)=
c\inv\sum_{h\in X(\oH_s,\tg)}
\Delta_{\oG}^{\tG}(h,\tg)||D_{\oG}(h)|^{-\frac12}\sum_{w\in W}\psi(wh).
\end{equation} 

Suppose $h\in X(\oH_s,\tg)$. From the definitions we have
\begin{equation}
\label{split3} 
\Delta _{\oG}^{\tG}(h,\tg)|D_{\oG}(h)|^{-\frac12}
=|D_{\tG}(\tilde g)|^{-\frac12}|e^{\rho}(h)|\tilde\chi_s(\tg)\chi_s^{-1}(h).
\end{equation} 
If $w\in W$ then $wh\in X(\oH_s,w\tg)$ and by Lemma \ref{invx} we have
\begin{equation}
\label{split4} 
\Delta_{\oG}^{\tG}(h,\tg)|D_{\oG}(h)|^{-\frac12}=
\Delta_{\oG}^{\tG}(wh,w\tilde g)
|D_{\oG}(wh)|^{-\frac12}.
\end{equation} 
Therefore
$\Lift_{\oG}^{\tG}(\Theta_{\pi(\psi)})(\tg)=$
\begin{equation}
\label{split5} 
c\inv\sum _{w \in W} |D_{\tG}(w\tg)|^{-\frac12} \tilde \chi _s(w\tg)
\sum_{h\in X(\oH_s,\tg)}
|e^{ \rho }(wh)| \chi _s^{-1}(wh)\psi (wh).
\end{equation} 

Since $\wt\chi_s \in X_g(Z(\tH_s))$ and $\chi_s |e^{-\rho}| \in
X(\oH_s)$, we can use the pair $(\wt\chi_s,\chi_s|e^{-\rho}|)$ 
to define
$\Lift_{\oH_s}^{\tH_s}$ as in \S \ref{s:tori}.
By Proposition \ref{p:tori}
and a short calculation, we can write
\begin{equation}
\label{split6}  
\Lift_{\oG}^{\tG}(\Theta_{\pi(\psi)}(\tg)=
 |D_{\tG}(\tg)|^{-\frac12}\sum _{w \in W}
\Lift_{\oH_s}^{\tH_s}(\psi))(w \tg).
\end{equation} 
By Proposition \ref{p:tori} the lift is non-zero if and only if 
\begin{equation}
\label{e:psih}
\psi(h)=|e^{-\rho}(h)|\chi_s(h)\quad(h\in \oH_s, \phi(h)=1).
\end{equation}
If $h=m_\alpha$ then by \eqref{e:cond}(b)  the right hand side
is trivial 
($\zetacx=1$  since there are no complex roots).
Therefore $\Lift_{\oG}^{\tG}=0$ independent of the choice of transfer
factors unless $\psi(m_\alpha)=1$ for all $\alpha\in \Phi$.

Assume \eqref{e:psih} holds.
We apply  Proposition \ref{p:tori} to the summands of
\eqref{split6}.
Let $\wt S=p\inv(\phi(\oH_s))\subset
Z(\tH_s)$ and let $n=|Z(\tH_s)/\wt S|^{\frac12}$.
Define a character of $\tS$:
\begin{equation} 
\wt\psi_0(\tg)=\wt\chi_s(\tg)\psi(h)\chi_s(h)\inv|e^{\rho}(h)|
\end{equation} 
where $\phi(h)=p(\tg)$.
Write
\begin{equation}
\label{e:Xgsplit}
\{\wt\psi\in X_g(Z(\tH_s))\,|\,\wt\psi|_{\wt S}=\wt\psi_0\}=
\{\wt\psi_1,\dots,\wt\psi_n\}.
\end{equation}

By Proposition \ref{p:tori}
\begin{equation}
\Lift_{\oH_s}^{\tH_s}(\psi)
=\sum_{i=1}^n \pi(\wt\tau(\wt \psi_i))
\end{equation}
and plugging this into \eqref{split6} gives
\begin{equation}
\label{split7}  
\Lift_{\oG}^{\tG}(\Theta_{\pi(\psi)})(\tg)=
|D_{\tG}(\tg)|^{-\frac12}\sum_{i=1}^n\sum _{w \in W}
\Theta_{\wt\tau(\wt\psi _i)}(w\tg).
\end{equation} 
It is straightforward to identify the right hand side of
\eqref{split7} with a sum of principal series characters.

Suppose $\wt\tau$ is  an irreducible genuine representation of
$\tH_s$. Associated to $\wt\tau$ is a genuine principal series
representation 
$\pi(\wt\tau)$
of $\tG$; see Section \ref{s:characterdata} for details.
Then
$\Theta_{\pi(\tilde \tau)}$ is supported on the conjugates of
$Z( \tH_s )$, and for  $\tg\in  \tH_s'$ we have
\begin{equation}
\label{split10} 
\Theta_{\pi(\wt\tau)}(\tg) =
|D_{\tG}(\tg)|^{-\frac12}\sum _{w
\in W} {\rm Tr}( \tilde\tau(w\tg) ).
\end{equation} 
Write $\wt\tau=\wt\tau(\wt\psi)$ for $\wt\psi$ a genuine character of
$Z(\tH_s)$ as in Lemma \ref{adams3}, and let
$\pi(\wt\psi)=\pi(\wt\tau(\wt\psi))$.
Comparing \eqref{split7} and  \eqref{split10} we obtain the following
result. Let
$\mu=\mu(\wt\chi_s,\chi_s)=d\chi_s-2d\wt\chi_s-\rho\in\z(\C)^*$
(Definition \ref{d:compatible}).

\begin{proposition}
Fix $(\wt\chi_s,\chi_s)\in \caS(H_s)$ for $(\tG,G,\oG)$,
and use this to define
$\Lift_{\oG}^{\tG}$. Define $\Lift_{\oH_s}^{\tH_s}$ using
$(\wt\chi_s,\chi_s|e^{-\rho}|)\in\caS(H_s)$ for $(\tH_s,H_s,\oH_s)$.
Suppose $\psi\in X(\oH _s)$. Then
$\Lift_{\oG}^{\tG}(\pi(\psi))=0$ unless \eqref{e:psih} holds,
so assume this is the case.
Define $\mu,\wt S,n=|Z(\tH_s)/\wt S|$ and $\{\wt\psi_1,\dots, \wt\psi_n\}$ as above.
Then
\begin{equation}
\label{e:liftps}
\Lift_{\oG}^{\tG}(\pi(\psi))=\sum_{i=1}^n \pi(\wt\psi_i).
\end{equation}
Each term in the sum has 
infinitesimal character $\frac12(d\psi-\mu)$.

Conversely, for  $\tilde \psi \in X_g(Z(\tH _s)) $ define
\begin{equation} \psi (h) = |e^{ -\rho }( h)| \chi _s(h) (\tilde \psi
\tilde \chi _s^{-1})(\phi (h))\quad(h \in \oH).
\end{equation} 
Then  $\Theta_{\pi(\wt\psi)}$ occurs in $\Lift _{\oG}^{\tG}(\pi(\psi))$.
\end{proposition} 

\begin{remark}
\label{r:distinctps}
Note that $\wt{H_s^0}\subset \wt S\subset
Z(\tH_s)=Z(\tG)\wt{H_s^0}$.  This implies
$d\wt\psi_i=\frac12(d\psi-\mu)$ (independent of $i$)
and that the restrictions of the
$\wt\psi_i$ to $Z(\tG)$ are distinct. Therefore the principal series
representations $\pi(\wt\psi_i)$ occuring in \eqref{e:liftps} have
distinct central characters, and are {\it a fortiori} not isomorphic.
\end{remark}

\begin{example}
Suppose $G(\C)$ is simple and simply connected.
In this case lifting is canonical (cf. Section \ref{s:special}) and $\mu=0$.
Also $Z(\tG)=\wt{Z(G)}$ and
$|Z(\tH_s)/\wt{H_s^0}|=|Z(G)|$.
In fact $|Z(G)|=4$ (in type $D_{2n}$), $2$ (in type $A_{2n+1}$,
$D_{2n+1}$ or $E_7$) or $1$ otherwise. See \cite[Table I]{shimura}. 

Suppose $\oG=G$. Then $\wt S=\tH_s^0$ and $n=|Z(G)|$. 
By \eqref{e:psih} $\Lift_G^{\tG}(\pi(\psi))\ne 0$ if and only if
$\pi(\psi)$ is spherical. Thus the lift of the 
spherical principal
series representation with infinitesimal character
$\lambda$ is the sum of the $n$ genuine principal series
representations of $\tG$ with infinitesimal character $\lambda/2$.

On the other hand suppose $C=Z(G)$; this is allowed by Lemma
\ref{zeta2}. Then $\wt S=Z(\tH)$, and the lift consists of a single
principal series.
\end{example}

\begin{example}
We specialize  the preceding example. Suppose $G=\oG=SL(2,\R)$.
Write $H_s\simeq\R^\times$ and suppose $\psi_\nu(x)=|x|^\nu$.
Then $\tH_s\simeq\R^\times\cup i\R^\times$. 
Suppose $\wt\psi_\nu^\pm(x)=|x|^{\nu}$ and $\wt\psi_\nu^{\pm}(i)=\pm i$. 
Then
\begin{equation}
\Lift_G^{\tG}(\pi(\psi_\nu))=\pi(\wt\psi_{\nu/2}^+)+\pi(\wt\psi_{\nu/2}^-).
\end{equation}
Note that $Z(\tG)=\{\pm 1,\pm i\}$ and $\psi_{\nu}^{\pm}$ have
different central characters. Also note that the character of
$\pi(\wt\psi_\nu^+)+\pi(\wt\psi_\nu^-)$ vanishes on $i\R^\times$.
In this case $\phi(H_s)=H_s^0$.

Now take $\oG=PSL(2,\R)\simeq SO(2,1)$.
Then (see Example
\ref{ex:pgl2}) $\phi:\oH_s\rightarrow H_s$ is an isomorphism.
By \eqref{e:psih} the lift of any principal series of $\oG$ is
non-zero. With the obvious notation we have
\begin{equation}
\Lift_{\oG}^{\tG}(\pi(\psi_\nu^\pm))=\pi(\wt\psi_{\nu/2}^{\pm}).
\end{equation}

\end{example}

\sec{Example: Discrete Series on the Compact Cartan}
\label{s:exds}

Suppose $G=\oG$ is connected, semisimple, acceptable, and equal rank. 
We assume also that $\tG$ is connected.   Let $H$ be a compact Cartan subgroup of $G$.
We compute the character formula for the lift of discrete series
representation restricted to $\tH$.

The group $H$ is connected.   
Suppose $\lambda\in \h(\C)^*$ satisfies
$\langle\lambda,\ch\alpha\rangle\in\Z_{\ne 0}$ for all roots
$\alpha$.  Associated to $\lambda$ is  a stable sum of  discrete series,
$\pi^{\st}(\lambda)$ of $G$, whose character $\Theta_{\pi^{\st}(\lambda)}$
satisfies
\begin{equation} 
\label{e:thetadsst}
\Theta_{\pi^{\st}(\lambda)}(h)=
(-1)^q \Delta (\Phi^+,h)^{-1} 
\sum_{w\in W}\epsilon(w)e^{w\lambda}(h)\quad(h \in H').
\end{equation} 
Here $\Phi^+=\Phi^+(\lambda)=\{\alpha\,|\, \langle\lambda,\ch\alpha\rangle>0\}$ and $q=\frac12\dim(G/K)$.

Suppose $(h,\tg)\in X(H,\tH)$.
By Definition \ref{d:transferspecial} (see \eqref{e:formulaforDelta})
the transfer factors are canonical, and given by
\begin{equation}
\Delta(h,\tg)=\frac{\Delta(h,\Phi^+)}{\Delta(\tg,\Phi^+)}\tau(h,\tg).
\end{equation}
By \eqref{e:thetadsst} and \eqref{e:lift1} we have:
\begin{equation}
\Lift_G^{\tG}(\Theta_{\pi^{\st}(\lambda)})(\tg)
=
c\inv (-1)^q\Delta(\tg,\Phi^+)\inv\sum_W\epsilon(w)\sum_{h\in X(H,\tg)}e^{w\lambda}(h)\tau(h,\tg).
\end{equation}

Define $\Lift_H^{\tH}$ using lifting data $(\wt\chi,\wt\chi^2)$ as in
Example \ref{ex:ctori}. 
Then the inner sum is equal to 
$c(H)\Lift_H^{\tH}(e^{w\lambda})(\tg)$. By  Definition
\ref{d:constants} $C(H)=c(H)/c$, so
\begin{equation}
\Lift_G^{\tG}(\Theta_{\pi^{\st}(\lambda)})(\tg)
=C(H)(-1)^q\Delta(\tg,\Phi^+)\inv\sum_W\epsilon(w)\Lift_H^{\tH}(e^{w\lambda})(\tg).
\end{equation}

By Example \ref{ex:ctori} $\Lift_H^{\tH}(e^{w\lambda})\ne 0$ if and
only  if  $e^{w\lambda/2}$ is a  genuine character of $\tH$.    
Let $W_\#$ be the set of $w\in W$ for which this holds. Assume $W_\#\ne \emptyset$.
By replacing $\lambda$ with $w\lambda$ we may assume
$1\in W_\#$. 
Then
\begin{equation}
\label{ex:Wsharpex}
W_\#=\{w\in W\,|\, w\frac\lambda2-\frac\lambda2\in X^*(H)\}
\end{equation}
where $X^*(H)$ is the character lattice of $H$. From Lemma \ref{adams1}
it is clear that 
\begin{equation}
W(G,H)=\{w\in W\,|\, w\frac\lambda2-\frac\lambda2\in R\}\subset W_\#.
\end{equation}
where $R$ is the root lattice.
A short calculation gives

\begin{multline}
\Lift_G^{\tG}(\Theta_{\pi^{\st}(\lambda)})(\tg)
=\\
C(H)\sum_{w\in W(G,H)\bs W_\#}(-1)^q\Delta(\tg,\Phi^+(w\lambda))\inv
\sum_{v\in W(G,H)}\epsilon(v)e^{vw\lambda/2}(\tg).
\end{multline}

For each $w\in W_\#$ the summand is the formula for the character of
$\pi(w\lambda/2)$, the genuine discrete series representation of $\tG$
with Harish-Chandra parameter $w\lambda/2$. Therefore
\begin{equation}
\label{e:oncompact}
\Lift_G^{\tG}(\Theta_{\pi^{\st}(\lambda)})(\tg)
=
C(H)\sum_{w\in W(G,H)\bs W_\#}\Theta_{\pi(w\lambda/2)}(\tg)\quad(\tg\in\tH').
\end{equation}

Since formulas for discrete series characters are more complicated on
noncompact Cartan subgroups, we will not directly compute
$\Lift_G^{\tG}(\Theta_{\pi^{\st}(\lambda)})(\tg)$
for  $\tg$ in an arbitrary
Cartan subgroup.  In Section \ref{s:liftinginveig} we will use results
of Hirai and Harish-Chandra to show that the lift of a stable,
supertempered eigendistribution on $G$ is a supertempered
eigendistribution on $\tG$.  Thus both sides of (\ref{e:oncompact})
are supertempered eigendistributions.  Now a theorem of Harish-Chandra
that says if two supertempered eigendistributions on $\tG$ agree on
the compact Cartan subgroup, they are equal.  This allows us to
conclude that (\ref{e:oncompact}) is valid for all $\tg \in \tG'$.
See Section \ref{s:discreteseries}.

\begin{remark}
\label{r:distinctds}
If $w\in W(G,H)\bs W_\#$ then $w\frac\lambda2-\frac\lambda2\in
X^*(H)/R\simeq \widehat{Z(G(\C))}$. Since $G$ is of equal rank
$Z(G(\C))=Z(G)$. Therefore $\pi(w\lambda/2)$ and $\pi(\lambda/2)$ have
distinct central characters, and the terms in \eqref{e:oncompact} are {\it a fortiori}
not isomorphic. Compare Remark \ref{r:distinctps}.
\end{remark}

\begin{example}
Let $G=Spin(2p,2q)$, and write $\lambda=(a_1,\dots,
a_p;b_1,\dots,b_q)$ in the usual coordinates. Then the lift is
non-zero if and only if all $a_i$ are even, and all $b_i$ are odd, or
vice-versa. If $p\ne q$ then the sum is a singleton, and
$\Lift_G^{\tG}(\pi(\lambda))=C(H)\pi(\lambda/2)$. 
If $p=q$ then the sum consists of two elements, $\lambda/2$ and
$\lambda'/2=\frac12(b_1,\dots, b_p;a_1,\dots, a_p)$.
\end{example}

\begin{example}
Suppose $G$ is compact, so $\tG\simeq G\times\Ztwo$ (which is
admissible, but not connected).    
Then $\pi(\lambda)$ is the finite dimensional
representation with infinitesimal character $\lambda$. 
The calculation of $\Lift_G^{\tG}(\pi(\lambda))$ is the same as above,
except that $\Lift_H^{\tH}(e^{w\lambda})\ne 0$ if and
only  if  $e^{w\lambda/2}$ is a character of $H$, in
which the genuine square root of $e^{w\lambda}$ is
$e^{w\lambda/2}\otimes\sgn$.
In this case $W = W(G,H) = W_\#$ so
the sum on the right hand side of \eqref{e:oncompact} has one element.
Thus
\begin{equation}
\Lift_G^{\tG}(\pi(\lambda))=\pi(\lambda/2)\otimes\sgn
\end{equation}
if $\lambda/2\in X^*(H)$, or is $0$ otherwise.
\end{example}

Also see Example \ref{ex:sl2ds}.

\sec{Invariant Eigendistributions}
\label{s:hirai}

In this section we review results of Hirai and Harish-Chandra on
invariant eigendistributions.  In the next section we apply these
results to study  lifting from $G$ to $\tG$.
 
Througout this section we assume $G$ is a group in 
Harish-Chandra's class (see the Appendix).
This holds for the groups $(\tG,G,\oG)$ in an admissible triple (see
the beginning of Section  \eqref{s:liftinginveig}).
Let $H$ be a Cartan subgroup of $G$.
Recall $H'$ is the set of regular elements of $H$ and
suppose $F:H'\rightarrow \C$ is differentiable.
Corresponding to $X\in\h$ is the differential operator

\begin{equation}\label{defdx} D_X F( h) = \left.\frac{d}{dt}\right|_{t=0} 
F( h\exp(tX))\quad(h \in H').
\end{equation} 
Fix a set of positive roots $\Phi^+$, with corresponding
$\rho=\rho(\Phi^+)$.
Then we also define

\begin{equation}
\label{defdxr}
D_X^{\rho} F( h) = \left.\frac{d}{dt}\right|_{t=0}
e^ {<\rho ,tX>} F( h\exp (tX))\quad(h \in H').
\end{equation} 
The map $D_X
\rightarrow D_X^{\rho} = <\rho ,X> +D_X$ can be extended to an
automorphism $D \rightarrow D^{\rho }$ of $S(\h (\C))$.

Let $\nu $ be a character of the center $\mathfrak Z$
of the universal enveloping algebra of $\g(\C)$.    Using the Harish-Chandra homomorphism
we identify $\nu $ with
a character of $I(\h(\C))=S(\h(\C))^W$.
 
 We say $F:H' \rightarrow \C$ satisfies
condition (C1, $\Phi ^+, \nu$ ) if $F$ is real analytic on $H'$ and

\begin{equation}
\label{e:C1}
\tag{C1,$\Phi^+,\nu$}
D^{\rho}F( h) = \nu (D)F(h)\quad(h \in H', D\in I(\h (\C))).
\end{equation} 
Let $H'(R) = \{ h \in H\,|\,e^{\alpha }(h) \not = 1, \alpha \in \Phi _r\}$.
We say $F$ satisfies condition (C2) if
\begin{equation}
\label{e:C2}
\tag{C2}
F\text{ extends to a real analytic function on }H'(R).
\end{equation}

Let $G'$ be the set of regular semisimple elements of $G$, and
suppose $\Theta $ is a class function on $G'$. With $\Delta^1(\Phi^+,h)$
defined as in \eqref{e:Deltas} let
\begin{equation}
\label{defpsi}
\Psi(H,\Phi^+,h)=\Delta^1(\Phi^+,h)\Theta(h)\quad(h\in H').
\end{equation} 
Assume $\Psi(H,\Phi^+)$ 
satisfies (C1, $\Phi ^+, \nu$ ) and (C2).

Let $\alpha$ be a simple root of $\Phi^+_r$. Let $J$ be the corresponding Cayley
transform of $H$ as in Section \ref{s:cartan}. 
Choose $c\in \Aut(\g (\C))$ 
satisfing \eqref{e:c}(a) and (b), and let
$\beta=c^*(\alpha)$ be the imaginary root of $J$
corresponding to $\alpha$.
Let $\Phi^+_J=c^*\Phi^+$ and $\rho_J=\rho(\Phi^+_J)$.
Let $H(\alpha )$ be the set of $h \in H$
such that $e^{\alpha }(h) = 1$, but $e^{\gamma}(h)\ne1$ for all
$\gamma \ne\pm\alpha$.

Fix $h \in H(\alpha)$.
Then $h\exp(t\ch \alpha )\in H'$ for $0\ne t$ sufficiently small,
and there are well-defined one-sided limits
\begin{equation}\label{deflimpm}
D_{\ch \alpha }^{\rho,\pm} \Psi( H,  \Phi ^+,h) = 
\lim _{t \rightarrow 0^{\pm }} D_{\ch\alpha }^{\rho} \Psi ( H,  \Phi ^+,h  \exp (t\ch \alpha )).
\end{equation}
Note that $i\ch\beta\in\j$ and $h\in J'(R)$.
Thus $D_{\ch \beta }^{\rho _J} = -i
D_{i\ch \beta }^{\rho _J}$ is defined as in (\ref{defdxr}), and
$D_{\ch \beta }^{\rho _J}\Psi (J,\Phi ^+_J,h )$ is defined.  We
say $\Psi (H,\Phi ^+)$ 
satisfies condition (C3) if

\begin{equation}
\tag{C3}
[D_{\ch\alpha}^{\rho,+}-D_{\ch\alpha}^{\rho,-}]\Psi(H,\Phi^+,h) 
= 2D_{\ch\beta }^{\rho _J}\Psi (J, \Phi ^+_J, h)\quad(h \in H(\alpha)).
\end{equation} 

The following theorem was proved by Hirai \cite{hirai_ie_1},
\cite{hirai_ie_2} assuming two extra conditions (see the Appendix)
which may fail to hold in our situation.  In the Appendix we extend
Hirai's theorem to a class of reductive groups that includes
Harish-Chandra's class.

\begin{theorem}[cf.~Theorem \ref{c123a}]
\label{c123}
Let $\Theta $ be a class function on $G'$ and let $\nu$ be a
character of $\mathfrak Z$.  
Then $\Theta $ is an invariant
eigendistribution with infinitesimal character $\nu$ if and only for
every Cartan subgroup $H$ of $G$ and every real root $\alpha$ of $H$, there
is a choice of positive roots $\Phi ^+$ such that $\alpha $ is simple
for $\Phi _r^+$ and the function $\Psi (H,\Phi ^+)$ satisfies
conditions (C1, $\Phi ^+, \nu$), (C2) , and (C3).

Moreover, suppose $\Theta $ is an invariant eigendistribution.
Then $\Psi (H,\Phi^+)$ satisfies (C1, $\Phi ^+, \nu$) and (C2) for
every choice of $\Phi ^+$, and satisfies (C3) for every choice of
$\Phi ^+$ such that $\alpha $ is simple for $\Phi _r^+$.
\end{theorem}

We now review some results of Harish-Chandra.
Write $G = K\exp(\p)$ where $K$ is a maximal compact subgroup and $\g =
\k \oplus \p$ is the Cartan decomposition.  Let $\Vert \cdot \Vert _G
$ be a Euclidean norm on $\p$ and define
\begin{equation} 
\label{defsigma} \tau _G(k \exp X) = \Vert X \Vert ,
\ k \in K ,X \in \p .  \end{equation} Let $\Theta $ be an invariant
eigendistribution on $G$.  Then by \S 12 of \cite{hc_ha_1}, $\Theta $ is
tempered if and only if for each Cartan subgroup $H$ of $G$ there are
numbers $C,r \geq 0$, so that
\begin{equation}\label{hctemp} |D(h)|^{\frac12}|\Theta (h)| \leq
C(1+\tau _G(h))^r, \ h \in H', \end{equation} 
(see \ref{e:Deltas} for the definition of $|D(h)|^{\frac12}$).

Let $\mathfrak z$ be the center of $\g $.  As in \cite{hc_ha_1} we can decompose
$G = {^0G} Z_{\mathfrak p}$ where $^0G$ has Lie algebra $ [\g, \g] \oplus (\z \cap
\k)$ and $Z_{\mathfrak p}= \exp (\z \cap \p)$.  Let $\Theta $ be an invariant
eigendistribution on $G$.  We say $\Theta $ is relatively tempered
(relatively supertempered) on $G$ if its restriction to $^0G$ is
tempered (supertempered) on $^0G$.  By \S 4 of \cite{hc_supertempered}, $\Theta $
is supertempered on $^0G$ if and only if for every Cartan subgroup
$^0H$ of $^0G$ and every $r \geq 0$,

\begin{equation}
\label{hcstemp} \sup _{g \in {^0H'}} |D(h)|^{\frac12}
|\Theta (h)| (1+\tau _G(h))^r < \infty.  \end{equation} We will also
use the following result from \S 4 of \cite{hc_supertempered}.

\begin{theorem}\label{hcsutemp}\cite{hc_supertempered} Assume that $G$ has a
relatively compact Cartan subgroup $B$ and let $\Theta $ be a
relatively supertempered invariant eigendistribution such that $\Theta
(b) = 0, b \in B'$.  Then $\Theta = 0$.
\end{theorem}

\sec{Lifting of Invariant Eigendistributions}
\label{s:liftinginveig}

We now return to the context of lifting. Suppose 
$(\tG, G, \oG )$ is an admissible triple.
We note that each of the groups $\tG,
G$, and $\oG$ are in Harish-Chandra's class (see the Appendix).
For $G$ and $\oG$ this is immediate. For $\tG$ note that $\Cent_{\tG}(\wt
G^0)=\wt{Z(G)}$ by \eqref{e:zgg0}. Therefore
$\Ad(\tG)\simeq\Ad(G)$, so $\tG$ also satisfies
Condition A of the Appendix, and the remaining conditions are easy.
Therefore the results of
Section \ref{s:hirai} apply.

Fix a maximally split Cartan subgroup $H_s$ of $G$ and lifting data
$(\wt\chi_s,\chi_s)\in\caS(H_s)$ as in 
Definition \ref{d:transfer}.
Recall $\mu\in\z(\C)^*$ is given by \eqref{e:mu}.
Let $\Theta $ be a stable invariant eigendistribution on ${\oG}$. 
Define
$\wt\Theta =\Lift_{\oG}^{\tG} \Theta$ as in Definition \ref{d:lift};
$\wt\Theta $ is a genuine class
function on $\tG'$.

\begin{theorem}\label{inveig} Let $\Theta $ be a stable invariant
eigendistribution on ${\oG}$ with infinitesimal character $\nu $.
Then ${\rm Lift}_{\oG}^{\tG} \Theta $ is an invariant
eigendistribution on $\tG$ with infinitesimal character $(\nu
-\mu )/2$.
\end{theorem}

We apply Theorem \ref{c123}. We need to show that conditions (C1-3) hold.

Let $H$ be a Cartan subgroup of $G$, and let $\Phi ^+$ be a
choice of positive roots.  Define
\begin{subequations}
\renewcommand{\theequation}{\theparentequation)(\alph{equation}}  
\begin{equation}
\label{deftpsi0g}
\Psi(\oH,\Phi^+,h)=\Delta^1(\Phi^+,h)\Theta(h)\quad(h\in\oH'),  
\end{equation}
\begin{equation}
\label{deftpsith}
\Psi(\tH,\Phi^+,\tg)
=
\Delta^1(\Phi^+,\tg)\tilde\Theta(\tg)\quad(\tg\in\tH').
\end{equation}
\end{subequations}

The second part of Theorem \ref{c123} says:

\begin{lemma}\label{co123} $\Psi (\oH , \Phi ^+)$ satisfies conditions
(C1, $\Phi ^+,\nu$) and (C2) for every choice of $\Phi ^+$.  Let
$\alpha\in\Phi _r^+$.  Then $\Psi (\oH , \Phi ^+)$ satisfies
condition (C3) for any choice of $\Phi ^+$ such that $\alpha$ is a
simple root for $\Phi _r^+$.
\end{lemma}

We now verify conditions (C1)-(C3) for $\Psi(\tH,\Phi^+)$.
We begin with (C1). 
Note that $\mu$ defines a character of $I(\h(\C))=S(\h(\C))^W$ and
also by the Harish-Chandra homomorphism a character of $\mathfrak Z$.
Assume $\Phi^+$ is a special set of positive roots.

\begin{proposition}
\label{lambda2} For all $D \in I(\h (\C ))$ and $\tg
\in \tH '$,
$$D^{\rho}\Psi(\tH,\Phi^+,\tg)=((\nu -\mu)/2)(D)\Psi(\tH,\Phi^+,\tg).
$$
In other words $\Psi(\tH,\Phi^+)$ satisfies condition (C1,$\Phi^+$,$(\nu-\mu)/2$).

\end{proposition}

Since $\tH,\oH$ and $\Phi^+$ are fixed we write
$\wt\Psi(\tg)= \Psi(\tH,\Phi^+,\tg)$ and $\overline\Psi(h)=\Psi(\oH, \Phi ^+,h)$.
Let $\Gamma=\Gamma(H,\Phi^+)$  (Definition~\ref{d:compatible}). 
Thus for $(h,\tg)\in X(\oH,\tH)$, 
$\Gamma(h,\tg)=\tilde\chi(\tg)\chi^{-1}(h)$
where $(\wt \chi, \chi ) \in S(H,\Phi ^+, \tilde \chi _s, \chi _s)$. 

\begin{lemma}
\label{psi0xg} 
For all $\tg \in \tH , X \in \h
$ such that $\tg \widetilde \exp X \in \tH'$,
$$ 
\wt\Psi(\tg\widetilde\exp X)=c^{-1}
\sum_{h\in X(\oH, \tg)}\Gamma(h,\tg)e^{-<\rho+\mu,X/2>}\overline\Psi(h\overline \exp (X/2) ).
$$
\end{lemma}

\begin{proof} Let $\tg \in \tH'$. 
Recall the definition of $X(\oH,\tg)$ \eqref{e:XoGtg}.
From the definitions it is clear that
$$
\wt\Psi(\tg)=c^{-1}\sum_{h\in X(\oH,\tg)}\Gamma(h,\tg)\overline\Psi(h).
$$

Fix $\tg \in \tH, X \in \h$ such that $\tg \widetilde
\exp X \in \tH'$.  Then $\phi (\overline \exp (X/2)) = \widetilde \exp X$ so $$X( \oH, \tg \widetilde \exp X) = \{
h\overline \exp (X/2): h \in X(\oH,  \tg ) \}.$$ Thus 
$$\wt\Psi(\tg \widetilde \exp X) = \sum _{ h \in X(\oH,  \tg ) }\Gamma(h,
\tg) \Gamma (\overline \exp
(X/2) , \widetilde \exp X)\overline\Psi(h\overline \exp (X/2).
$$ 
But
$$ 
 \Gamma (\overline \exp
(X/2) , \widetilde \exp X) =
  e^{<2d\tilde \chi -  d\chi    ,X/2>} =  e^{-< \rho+\mu , X/2>}
$$ by \eqref{e:mu}.
\end{proof}

\begin{proof}[Proof of Proposition \ref{lambda2}]
  
It follows from Lemmas \ref{co123} and \ref{psi0xg}
that $ \wt\Psi$ is real analytic on $\tH'$.
Suppose $X
\in \h$ and $t \in \R$.  Using Lemma \ref{psi0xg} we have
 $$
e^{<\rho, tX>} \wt\Psi(\tg \widetilde \exp tX ) = c^{-1}
\sum _{ h \in X(\oH, \tg ) }\Gamma(h, \tg)  e^{<\rho -\mu , tX/2  >}  
\overline\Psi(h\overline \exp (tX/2) ).
$$
Thus $$D_X^{\rho} \wt\Psi(\tg ) =c^{-1}
 \sum _{ h \in X(\oH,\tg ) } \Gamma(h, \tg) (D_{X/2}^{\rho} -<\mu ,X/2>) 
\overline\Psi(h) .
$$ 

First assume $X\in \z$. Then 
$D_{X/2} \in I(\h (\C ))$, so
$ D_{X/2}^{\rho} \overline\Psi( h ) =\nu (X/2) \overline\Psi( h )$ 
by Lemma \ref{co123}.  Thus
$$
D_X^{\rho}  \wt\Psi(\tg )  = c^{-1} \sum _{ h \in X( \oH, \tg ) } \Gamma(h, \tg) <\nu - \mu ,X/2> \overline\Psi(   h ) =
((\nu -\mu)/2)(D_X)\wt\Psi(\tg ).
$$ 
By induction we
see that 
\begin{equation}
\label{e:Drhoz}
 D^{\rho} \wt\Psi( \tg ) = ((\nu
-\mu)/2)(D) \wt\Psi( \tg )\quad\text{for all }D \in S(\z (\C )).
\end{equation}

Now suppose $X \in \h _d$. Then $<\mu,X> = 0$,
so
\begin{equation}
\label{e:DxrhoPsi}
D_X^{\rho} \wt\Psi( \tg ) =
\frac1{2c}\sum _{ h \in X(\oH, \tg ) } \Gamma(h, \tg)D_{X}^{\rho} 
\overline\Psi( h).
\end{equation}

Now suppose $Y\in\h _d$, replace $\tg$ with $\tg\wt\exp(tY)$, and
multiply both sides by $e^{\langle\rho,tY\rangle}$. The right hand
side is
\begin{equation}
e^{\langle\rho,tY\rangle}\frac1{2c}\sum_{h\in X(\oH, \tg)}\Gamma(h\overline\exp(tY/2),\tg\wt\exp(tY))D_{X}^{\rho}
\overline\Psi(h\overline\exp(tY/2)).
\end{equation}
As in the proof of Lemma \ref{psi0xg} we can write
$$
\Gamma(h\overline\exp(tY/2),\tg\wt\exp(tY))=\Gamma(h,\tg)e^{\langle-\rho,tY/2\rangle}
$$
and we see
\begin{equation}
e^{\langle\rho,tY\rangle}
D_X^{\rho} \wt\Psi(\tg\wt\exp(tY) ) =
\frac1{2c}\sum_{h\in X(\oH,\tg)}\Gamma(h,\tg)e^{\langle\rho,tY/2\rangle}D_{X}^{\rho}\overline\Psi(h\overline\exp(tY/2)).
\end{equation}
Taking the derivative and evaluating at $t=0$ we see
\begin{equation}
D_Y^{\rho}D_X^{\rho}
 \wt\Psi(\tg ) =
\frac1{4c}\sum_{h\in X(\oH,\tg)}\Gamma(h,\tg)D_Y^\rho D_{X}^{\rho}\overline\Psi(h).
\end{equation}
By induction we conclude that for $D\in I(\h_d(\C))$, homogeneous of
degree $k$, we have
\begin{equation}
\label{e:Drhohd}
\begin{aligned}
D^{\rho}\wt\Psi(\tg)&=
\frac1{2^kc}
\sum_{h\in X(\oH,\tg)}\Gamma(h,\tg)D^{\rho}\overline\Psi(h)\\
&=\frac1{2^kc}
\nu(D)\sum_{h\in X(\oH,\tg)}\Gamma(h,\tg)\overline\Psi(h)\\
&=(\nu/2)(D)\wt\Psi(\tg)\\
&=((\nu-\mu)/2)(D)\wt\Psi(\tg).
\end{aligned}
\end{equation}
The last equality comes from the fact that $\mu(D)=0$.

We have $ I(\h (\C )) = I(\h _d(\C ))S(\z (\C ))$,
and the result follows from \eqref{e:Drhohd} and \eqref{e:Drhoz}.
\end{proof}

We now consider condition (C2).
We continue to write $\wt\Psi(\tg)=\Psi(\tH,\Phi^+,\tg)$ and
$\overline\Psi(h)=\Psi(\oH,\Phi^+,h)$. 
Recall
\begin{equation}
\label{tgdrt} 
\tH '(R) = \{ \tg \in \tH\,|\, e^{\alpha }(\tg) \not = 1, \ \alpha \in \Phi _r
\}. \end{equation}

\begin{proposition}\label{analytic} $ \wt\Psi$ extends to a real
analytic function on $\tH'(R)$, i.e. $\wt\Psi$ satisfies
condition (C2).
\end{proposition}

\begin{proof} Since the restriction of $\tilde \Theta $ to $\tH$
is supported on $Z(\tG)\tH^0$, $\wt\Psi$ is also
supported on $Z(\tG)\tH^0$.  Fix $\tilde t \in Z(\tilde
G)\tilde T^0$.   If $\Phi^+_r$ is a set of positive roots of $\Phi_r$
let

\begin{equation}
\label{e:C}
C(\Phi^+_r)=\{X\in\a\,|\, \alpha(X)>0\text{ for all }
\alpha\in\Phi^+_r\}.
\end{equation}
Then the connected components of 
$\tilde t\wt{H^0}\cap\wt{H'}(R)$ are of the form
$\tilde t\wt{T^0}\wt\exp(C(\Phi^+_r))$.

Fix $h \in X(\oH, \tilde t)$.
Note that for all $\alpha \in \Phi _r$, 
$e^{2\alpha}(h)=e^{\alpha}(\tilde t)=1$, so that $e^{\alpha}(h)=\pm1$.
Let $\Phi_r(h)=\{\alpha\in\Phi_r\,|\,e^\alpha(h)=1\}$.
Then the  connected components of $h\oH ^0 \cap \oH '(R)$ have the
form  $h\oT ^0 \overline \exp(C(\Phi _r^+(h)))$ with $C(\Phi_r^+(h))$
defined as in \eqref{e:C}.

Given $\Phi^+_r$, let $\Phi^+_r(h)=\Phi^+_r\cap\Phi_r(h)$.
If $X\in\t+C(\Phi^+_r)$ then $X/2\in\t+C(\Phi^+_r(h))$.
Thus $X \rightarrow e^{<-\rho -\mu , X/2>} \overline\Psi(h
\overline \exp (X/2))$ extends to a real analytic function on $\t +
C(\Phi _r^+)$.  Therefore 

\begin{equation}
\label{e:extends}
X \rightarrow \sum _{ h\in X(\oH, \tilde t) } \Gamma(h,
\tilde t) e^{<-\rho -\mu , X/2>} \overline\Psi(h\overline \exp (X/2))
\end{equation}
extends to a real analytic function on $\t + C(\Phi _r^+)$.  

By Lemma \ref{psi0xg}, for all $X\in\h$ such that 
$\tilde t\widetilde\exp X\in\tH'$, 
\begin{equation}
\wt\Psi(\tilde t\widetilde\exp X)=c^{-1}\sum_{h\in X(\oH, \tilde t)}\Gamma(h,\tilde t) 
e^{<-\rho-\mu,X/2>}\overline\Psi(h\overline\exp (X/2)).
\end{equation}
By \eqref{e:extends}
there is an open neighborhood $U$ of $\tilde t$ in $\tilde t \tilde T
^0$ such that $ \wt\Psi$ extends to be real analytic on
$U\widetilde \exp ( C(\Phi _r^+))$.
\end{proof}

We now turn to condition (C3).
Fix $\alpha \in \Phi _r$ and let 
$J$ be the Cayley transform of $H$ as in Section \ref{s:cartan}.
Choose $c\in\Aut(\g (\C))$ satisfying \eqref{e:c}(a) and (b), and let
$\beta=c^*(\alpha)$. 
  
\begin{lemma}\label{sstable} Assume that $\Phi ^+$ is a special set of
positive roots, $\alpha\in\Phi ^+$, and that $\Phi _J^+ =c^*\Phi ^+$ is also
special.  Then $\alpha $ is a simple root for $\Phi _r^+$.
\end{lemma}

\begin{proof}
A basic property of Cayley transforms is that for any root $\gamma$ of $H$,
$\sigma(c^*\gamma)=c^*(s_\alpha\sigma\gamma)$. 
Suppose $\gamma\in\Phi_r^+$, $\gamma \not = \alpha$.   Then 
\begin{equation}
\label{e:cbeta}
\sigma(c^*\gamma)=c^*(s_\alpha\gamma).
\end{equation}
If $\langle\gamma,\ch\alpha\rangle=0$
then $s_\alpha(\gamma)=\gamma\in\Phi^+_r$. 
If $\langle\gamma,\ch\alpha\rangle\ne 0$ then
$c^*\gamma\in\Phi^+_J$, and by \eqref{e:cbeta} $c^*\gamma$ is complex. 
Since $\Phi^+_J$ is special  $\sigma(c^*\gamma)\in\Phi^+_J$, and so by \eqref{e:cbeta} again
$s_\alpha\gamma\in\Phi^+_r$. 
Therefore $s_\alpha\gamma\in\Phi^+_r$ for all $\alpha\ne\gamma\in\Phi^+_r$, so
$\alpha$ is simple for $\Phi^+_r$.
\end{proof}

\begin{proposition}
\label{jump} 
Suppose $\Phi^+$ is a special set of positive roots for $H$, such that
$\Phi^+_J$ is a special set of positive roots for $J$. 
Let $\tg \in \tH (\alpha )$.  Then
$$
[D_{\ch\alpha}^{\rho,+}-D_{\ch\alpha}^{\rho,-}]\Psi(\tH,\Phi^+,\tg)=
2D_{\ch\beta}^{\rho_J}\Psi(\wt J,\Phi_J^+,\tg).
$$ 
That is
 $\Psi(\tH,\Phi^+)$ satisfies condition (C3).
\end{proposition}

\begin{proof} 
Define
$\Gamma(H,\Phi^+)$ and $\Gamma(J,\Phi^+_J)$  as in Definition \ref{d:compatible}.

By \eqref{e:DxrhoPsi} we have
$$
[D_{\ch \alpha}^{\rho,+}-D_{\ch \alpha}^{\rho,-}]
\Psi(\tH,\Phi^+,\tg)
=\frac1{2c}\sum_{h\in
X(\oH,\tg)}\Gamma(H,\Phi^+,h,\tg)[D_{\ch\alpha}^{\rho,+}-D_{\ch\alpha}^{\rho,-}]\Psi(\oH,\Phi^+,h).
$$

If $h\in X(\oH,\tg)$ 
then for all 
$\gamma \in \Phi $ by \eqref{e:2alpha} we have $e^{2\gamma }(h ) = e^{\gamma }(\tg)$.
Thus $e^\alpha(h)=\pm1$ and $e^\gamma(h)\ne\pm1$ for all
$\gamma\ne\pm\alpha$.
Thus
$X(\oH,\tg)=X(\oH,\tg,\alpha)_+\cup X(\oH,\tg,\alpha)_-$ where
\begin{equation}
\label{e:Xtgalpha}
X(\oH,\tg, \alpha )_\pm = \{ h \in X(\oH,\tg): e^{\alpha }(h ) =\pm1\}.
\end{equation}

Suppose that $h\in X(\oH,\tg,\alpha)_-$. Then $h$ is regular so
that
\begin{equation}
D_{\ch \alpha }^{\rho,+}\Psi(\oH,\Phi^+,h)=
D_{\ch \alpha }^{\rho,-}\Psi(\oH,\Phi^+,h)\quad(h\in X(\oH,\tg,\alpha)_-).
\end{equation}
On the other hand if $h\in X(\oH,\tg,\alpha)_+$, since 
$\Psi (\oH,\Phi^+)$ satisfies (C3),
\begin{equation}
[D_{\ch \alpha }^{\rho,+}-D_{\ch \alpha }^{\rho,-}]\Psi(\oH,\Phi^+,h) 
= 2 D_{\ch\beta}^{\rho _J} \Psi (\oJ,\Phi^+_J,h)
\end{equation}
where $\beta=c^*\alpha$.
Therefore
\begin{subequations}
\renewcommand{\theequation}{\theparentequation)(\alph{equation}}  
\begin{equation}
\label{e:left}
[D_{\ch \alpha}^{\rho,+}-D_{\ch \alpha}^{\rho,-}]\Psi(\tH,\Phi^+,\tg)  
= c^{-1}\sum _{h \in X(\oH,\tg, \alpha)_+}\Gamma(H,\Phi^+,h, \tg ) D_{\ch\beta}^{\rho _J}
\Psi(\oJ,\Phi^+_J,h).
\end{equation}

On the other hand, by \eqref{e:DxrhoPsi}
applied to $\tilde J$
\begin{equation}
2D_{\ch\beta}^{\rho_J}\Psi(\tilde J,\Phi^+_J,\tg)
=c^{-1}
\sum_{b\in X(\oJ,\tg)}
\Gamma(J,\Phi^+_J,b,\tg)D_{\ch\beta}^{\rho_J}\Psi(\oJ,\Phi^+_J,b).
\end{equation}

Write $X(\oJ,\tg)=X(\oJ,\tg,\beta)_+
\cup X(\oJ,\tg,\beta)_-$ as above, so we have
\begin{equation}
\begin{aligned}
2D_{\ch\beta}^{\rho_J}\Psi(\tilde J,\Phi^+_J,\tg)&=
c^{-1}\sum_{b\in X(\oJ,\tg,\beta)_+}\Gamma(J,\Phi^+_J,b,\tg)D_{\ch\beta}^{\rho_J}\Psi(\oJ,\Phi^+_J,b)\\
&+
c^{-1}\sum_{b\in
X(\oJ,\tg,\beta)_-}\Gamma(J,\Phi^+_J,b,\tg)D_{\ch\beta}^{\rho_J}\Psi(\oJ,\Phi^+_J,b)
\end{aligned}
\end{equation}

Suppose $h\in X(\oJ,\tg,\beta)_+$. Then $h\in\oH$, so $h\in
X(\oH,\tg,\alpha)_+$, and 
$X(\oJ,\tg,\beta)_+=X(\oH,\tg,\alpha)_+$. 
Also by Lemma \ref{l:HHalpha}
$\Gamma(H,\Phi^+,h,\tg)=\Gamma(J,\Phi^+_J,h,\tg)$. 
Therefore

\begin{equation}
\label{e:right}
\begin{aligned}
2D_{\ch\beta}^{\rho_J}\Psi(\tilde J,\Phi^+_J,\tg)&=
c^{-1}\sum_{b\in X(\oH,\tg,\alpha)_+}\Gamma( H,\Phi^+,b,\tg)D_{\ch\beta}^{\rho_J}\Psi(\oJ,\Phi^+_J,b)\\
&+
c^{-1}\sum_{b\in
X(\oJ,\tg,\beta)_-}\Gamma(J,\Phi^+_J,b,\tg)D_{\ch\beta}^{\rho_J}\Psi(\oJ,\Phi^+_J,b)
\end{aligned}
\end{equation}
It is enough to show the second sum on the right hand side is $0$,
 for then the result follows from \eqref{e:left} and \eqref{e:right}.
 \end{subequations}
We restate this as:
\begin{lemma}
\label{l:jump}
Suppose $\alpha$ is an imaginary noncompact root and
$\tg\in\tH(\alpha)$.
Then
\begin{equation}
\sum_{h\in X(\oH,\tg,\alpha)_-}
\Gamma(h,\tg)D_{\ch\alpha}^\rho\Psi(\oH,\Phi^+,h)=0.
\end{equation}
\end{lemma}

\begin{proof}
Since $\oH$ and $\Phi^+$ are fixed let
$\overline\Psi(h)=\Psi(\oH,\Phi^+,h)$.
Let $w$ be the reflection in $\alpha$ and suppose $h\in\oH'$.
Since $\alpha$ is imaginary
$\epsilon_R(\Phi^+,wh)=\epsilon_R(\Phi^+,h)$,
and since $\Theta$ is stable
$\Theta(wh)=\Theta(h)$.
Furthermore since $\epsilon (w) = -1$, 
$$
\Delta^0(\Phi^+,wh)=-e^{\rho-w\rho}(h)\Delta^0(\Phi^+,h).
$$ 
Thus
$$
\overline\Psi(wh)=-e^{\rho-w\rho}(h)\overline\Psi(h).
$$

Now assume $h\in X(\tg,\alpha)_-$, so $e^\alpha(h)=-1$. Then
\begin{equation}
\label{e:wh}
wh=h\overline\exp(\pi i\ch\alpha).
\end{equation}
Therefore $\phi(wh)=p(\tg)$ and $e^\alpha(wh)=-1$, so $wh\in
X(\tg,\alpha)_-$.
Also, for all $t\in\R$,
\begin{equation}
e^{<\rho,it\ch\alpha>}\overline\Psi((wh)\overline\exp(it\ch\alpha))=
-e^{\rho-w\rho}(h)e^{<\rho,-it\ch\alpha>}\overline\Psi(h\overline\exp(-it\ch\alpha)).
\end{equation}
It follows from this and the definitions that
\begin{equation}
D_{\ch\alpha}^\rho\overline\Psi(wh)=
e^{\rho-w\rho}(h) D_{\ch\alpha}^\rho\overline\Psi(h).
\end{equation}
From \eqref{e:wh} we have
$e^{\rho-w\rho}(h)=(-1)^{\langle\rho,\ch\alpha\rangle}$ and
\begin{equation}
\Gamma(wh,\tg)=
\frac{\wt\chi(\tg)}{\chi(wh)}=
\frac{\wt\chi(\tg)}{\chi(h)\chi (\overline\exp(\pi i\ch\alpha))}=
\Gamma(h,\tg)(-1)^{\langle d\chi,\ch\alpha\rangle},
\end{equation}
where $(\tilde \chi, \chi ) \in S(H,\Phi ^+,\tilde \chi _s,\chi _s)$.
Thus
\begin{equation}
\Gamma(wh,\tg)D_{\ch\alpha}^{\rho}\Psi(\oH,wh)=
(-1)^{<d\chi+\rho,\ch\alpha>}\Gamma(h,\tg)D_{\ch\alpha}^{\rho}\Psi(\oH,h).
\end{equation}
By \eqref{e:cond}
$$
\langle d\chi+\rho,\ch\alpha\rangle=
\langle 2d\tilde\chi+2\rho,\ch\alpha\rangle
\equiv\langle 2d\tilde\chi,\ch\alpha\rangle\pmod2.
$$ 
By \cite[Section 7]{dualityonerootlength} 
$\langle 2d\tilde\chi,\ch\alpha\rangle=1\pmod2$.
Therefore
\begin{equation}
\Gamma(wh,\tg)D_{\ch\alpha}^{\rho}\Psi(\oH,wh)=
-\Gamma(h,\tg)D_{\ch\alpha}^{\rho}\Psi(\oH,h).
\end{equation}
Therefore if $wh\ne h$ then the two terms cancel, and if $wh=h$, 
$D_{\ch\alpha}^\rho\Psi(h)=0$.
\end{proof}

This completes the proof of  Proposition  \ref{jump}.
\end{proof}

Theorem \ref{inveig}  follows from Theorem \ref{c123}, Propositions
\ref{lambda2}, \ref{analytic} and \ref{jump}.

\begin{remark}
\label{r:stable}
The only place that stability is used in the proof of Theorem
\ref{inveig} 
is in the proof of Lemma \ref{l:jump}.
\end{remark}

\medskip

It is easy to see that lifting preserves the property of being
relatively tempered and supertempered, and (provided $\mu\in i\z^*$) tempered.

\begin{proposition}
 Let $\Theta $ be a stable
invariant eigendistribution on ${\oG}$.

\begin{enumerate}
\item If $\Theta$ is relatively tempered then ${\rm
Lift}_{\oG}^{\tG}\Theta $ is relatively tempered.
\item If $\Theta$ is relatively supertempered then ${\rm
    Lift}_{\oG}^{\tG}\Theta $ is relatively supertempered.
\item Assume $\mu\in i\z^*$. If $\Theta$ is tempered
then ${\rm Lift}_{\oG}^{\tG}\Theta $ is tempered.
\end{enumerate}
\end{proposition}

\begin{proof}
Since $\oG$ and $\tG$ have the same Lie algebra, we
can use the same Euclidean norm $\Vert \cdot \Vert$ on $\p$ to define
both $\tau _H$ and $\tau _{\tG}$.

Assume $\Theta$ is tempered and $\mu\in i\z^*$.
Let $\h $ be a Cartan subalgebra of $\g $.  Fix $\tg \in \tilde
H'$ and write $\tg = \tilde t \widetilde \exp X$ where $\tilde t
\in \tilde T ', X \in \a $.  Then $\tau _{\tG}(\tg ) = \Vert
X \Vert $.  Further, for any $h \in \oH$ such that $\phi (h)=p(\tilde
g)$, we have $h=t\overline{\exp}X/2$ where $t \in \oT$ with $\phi (t) =
p(\tilde t)$.  Thus $\tau _H(h) = \Vert X/2 \Vert = \tau _{\tilde
G}(\tg )/2$.

Now, using \eqref{e:D}(d) and (\ref{hctemp}) applied to $\oG$, there are $C,r \geq 0$ such that
$$|D(\tg)|^{\frac12}|\tilde \Theta (\tg)| \leq c^{-1}\sum _{ h \in X(\oH \tg)  }
|D(\tg)|^{\frac12} |\Delta_{\oG}^{\tG}(h,\tg)|
|\Theta (h)| $$ $$=c^{-1} \sum _{ h \in X(\oH \tg)  }
|D(h)|^{\frac12}|\Theta (h)| \leq c^{-1}[\oH _1] C(1+\tau _{\tilde
G}(\tg)/2)^r,$$ where $\oH _1 = \{ h \in \oH : \phi (h) = 1 \}$.
Thus there is a constant $C'$ such that
$$|D(\tg)|^{\frac12}  |\tilde \Theta (\tg)| \leq C'(1+\tau _{\tG}(\tg) )^r $$
for all $\tg \in \tH'$.  Now $\tilde \Theta $ is tempered
using (\ref{hctemp}) applied to $\tG$. This proves (3).

Part (1) is similar, 
using
the fact that by \eqref{e:D}(c)
$|D(\tg)|^{\frac12} |\Delta_{\oG}^{\tG}(h,\tg)| = |D(h)|^{\frac12}$
for any choice of $(\tilde \chi _s,\chi _s)$ provided $\tg \in {^0\tG}$.
Part (2) is also similar, using \eqref{hcstemp} in place of
\eqref{hctemp}.
\end{proof}

\sec{Cuspidal Levi Subgroups}
\label{s:levi}

Let $(\tG, G, \oG )$ be an admissible triple, and fix a cuspidal Levi
subgroup $M$ of $G$.
Let $\tM=p\inv(M)$, and $\oM=\overline p(M)$.
These are cuspidal Levi subgroups of $\tG$ and $\oG$ respectively.

\begin{proposition}
\label{cl1} 
$(\tM,M,\oM )$ is an admissible triple (Definition \ref{assumptions}).
\end{proposition}

\begin{proof} 
It is well known that $M$ is the set of real points of
the connected reductive complex group $M(\C )$ and that $M_d(\C )$ is
acceptable, and it is clear that every simple factor of $M(\C)$ is
oddly laced.
Let $\Phi_M$ be the set of roots of $H(\C)$ in $M(\C)$ where $H=TA$ is a relatively compact Cartan subgroup of $M$.
Since $\tG$ is an admissible cover of $G$,
every noncompact root in $\Phi_M$ is metaplectic.  Thus $p: \tM
\rightarrow M$ is also an admissible cover.  Further, $\oM$ 
is the set of real points of $\oM(\C) =M(\C)/C$, where $C
\subset Z_0(G) \subset Z_0(M )$ with $c^2=1$ for all $c \in C$.  
This leaves only condition 3(c) of Definition \ref{assumptions}.

Let
$\Phi ^+$ be a special set of positive roots for $H$ in
$G$ and let $\Phi _M^+ = \Phi ^+ \cap \Phi _M$.  Let $\rho = \rho
(\Phi ^+)$ and $\rho _M = \rho (\Phi _M^+)$.  Then $\rho -\rho _M$ is
zero on $\t  $.  But $C \cap M_d^0 \subset T ^0$, so that $e^{\rho
}(c) = e^{\rho _M}(c)$ for all $c \in C \cap M_d^0$.  Thus if $\tilde
\chi $ is a genuine character of $Z(\tilde M)$, then $\tilde \chi ^2(c)
e^{\rho _M}(c) = \tilde \chi ^2(c) e^{\rho }(c) =1$ for all $c \in C
\cap M_d^0$.
\end{proof}

Suppose $H$ is an arbitrary Cartan subgroup of $M$.
Let $\Phi$ and $\Phi_M$ be the sets of roots of $H(\C)$ in $G(\C)$
and $M(\C)$, respectively.     Let $\Phi^+$ be a special set of positive
roots of $\Phi$, and let $\Phi^+_M=\Phi^+\cap \Phi_M$. This is a special
set of positive roots of $\Phi_M$. Suppose $\chi$ is a character of
$H$.
Define
\begin{equation}
\label{e:chiM}
\chi_M(h)=|e^{\rho_M-\rho}(h)|\chi(h)\quad(h\in \oH).
\end{equation}
To be precise we are writing
$|e^{\rho_M-\rho}(h)|$ for
$|e^{2\rho_M-2\rho}(h)|^{\frac12}$. Note that
$\Phi_i \subset \Phi_M$, so $2\rho_M-2\rho$ is a sum of real
and complex roots; since $\Phi^+$ is special this takes real values on
$\h$. Therefore $\rho_M-\rho$ exponentiates to $\oH^0$, and
\begin{equation}
\label{e:rhorhoMH0}
|e^{\rho_M-\rho}(h)|=
e^{\rho_M-\rho}(h)\quad(h\in \oH^0).
\end{equation}

In particular let $H_s$ be  a maximally split Cartan subgroup of $M$.  
This is also a maximally split Cartan subgroup of $G$.
Let $\Phi _s^+$ be a special set of positive roots for $H_s(\C)$ in
$G(\C)$. Let $\Phi_{M,s}$ be the roots of $H_s(\C)$ in $M(\C)$, and
let
$\Phi _{M,s}^+ =\Phi_s^+ \cap \Phi _{M,s}$.

\begin{lemma}
\label{l:chims}
Fix $(\wt\chi_s,\chi_s)\in\caS(H_s,\Phi_s^+)$.

\smallskip
\noindent (1) Let
$\chi_{_{M,s}}=(\chi_s)_M$ (cf. \ref{e:chiM}). 
Then $(\wt\chi_s,\chi_{_{M,s}})\in \caS(H_s,\Phi^+_{M,s})$.

\smallskip
\noindent (2) Let $H$ be any Cartan subgroup of $M$ and suppose
$(\wt\chi,\chi)\in \caS(H,\Phi^+,\wt\chi_s,\chi_s)$.
Then $(\wt\chi,\chi_M)\in \caS(H,\Phi^+_M,\wt\chi_s,\chi_{M,s})$.

\smallskip
\noindent (3) Let $\mu=\mu(\wt\chi_s,\chi_s)$ and
$\mu_M=\mu(\wt\chi_s,\chi_{M,s})$ (defined with respect to $M$). Then
$\mu=\mu_M$. 
\end{lemma}

\begin{proof} 
Fix $h\in (\oH_s\cap M_d)^0$. 
By \eqref{e:rhorhoMH0} $\chi_{M,s}(h)=\chi_s(h)e^{\rho_M-\rho}(h)$.
On the other hand since 
$(\oH_s\cap M_d)^0\subset (\oH_s\cap G_d)^0$, by \eqref{e:cond}(a) we have
$\chi_s(h)=(\wt\chi^2e^\rho)(h)$. Therefore
\begin{equation}
\chi_{M,s}(h)=\chi_s(h)e^{\rho_M-\rho}(h)=
(\wt\chi^2e^\rho)(h)e^{\rho_M-\rho}(h)=(\wt\chi_s^2e^{\rho_M})(h).
\end{equation}
Therefore \eqref{e:cond}(a) holds.

Next let $\Gamma_r(\oG,\oH)$ (resp. $\Gamma_r(\oM,\oH)$) be the group
$\Gamma_r(\oH)$ with respect to the group $\oG$ (resp. $\oM$) (cf.
\eqref{e:GammaR}).
Then $\Gamma_r(\oM,\oH)\subset \Gamma_r(\oG,\oH)$.
Fix $h\in\Gamma_r(\oM,\oH)$. Then $|e^{\rho_M-\rho}(h)|=1$, and by 
Lemma \ref{l:cayleyH1H2}
\begin{equation}
\chi_{M,s}(h)=\chi_s(h)=\zeta_{cx}(\oG,\oH)(h)=
\zeta_{cx}(\oM,\oH)(h).
\end{equation}
This verifies condition 
\eqref{e:cond}(b), and proves (1).

Now consider  (2).
Let $\chi_0$ be the character of $\oH$ such that
$(\wt\chi,\chi_0)\in\caS(H,\Phi^+_M,\wt\chi_s,\chi_{M,s})$
(cf. Proposition \ref{p:Gammamain} and Definition \ref{d:compatible}).
We need to show $\chi_0=\chi_M$.
We may assume $A\subset A_s$.
Recall $\oH=\Gamma(\oH)Z(\oG)\oH_d^0$. First assume
$h\in\Gamma(\oH)Z(\oG)\subset \Gamma(\oH_s)Z(\oG)$. 
Then $|e^{\rho_M-\rho}(h)|=|e^{\rho_{M,s}-\rho_s}(h)|=1$, so
$\chi_M(h)=\chi(h)$ and $\chi_{M,s}(h)=\chi_s(h)$.
Therefore
\begin{equation}
\begin{aligned}
\chi_0(h)&=
(\wt\chi/\wt\chi_s)(\phi(h))\chi_{M,s}(h)\quad(\text{by }\eqref{e:compatible})\\
&=
(\wt\chi/\wt\chi_{s})(\phi(h))\chi_{s}(h)\\
&=\chi(h)\quad(\text{by }\eqref{e:compatible})\\
&=\chi_M(h).
\end{aligned}
\end{equation}
Now suppose $h=(\oH\cap M_d)^0\subset \oH^0_d$.
Then $\chi(h)=(\wt\chi^2e^\rho)(h)$, and
\begin{equation}
\chi_0(h)=(\wt\chi^2e^{\rho_M})(h)=(\chi e^{-\rho}e^{\rho_M})(h)=\chi_M(h).
\end{equation}

Next let $X \in \h_d\cap \z_M $ where $\z_M$ is the center of $\liem$.    Write $h= \overline{\exp}X$.    

Then:
\begin{equation}
\begin{aligned}
\chi_0(h)&=(\wt\chi/\wt\chi_s)(\phi(h))\chi_{M,s}(h)\\
&=(\wt\chi/\wt\chi_s)(\phi (h))\chi_s(h)e^{\rho_{M,s}-\rho_s}(h) \\
&=(\wt\chi/\wt\chi_s)(\phi (h))(\wt\chi_s^2e^{\rho_{s}})(h)e^{\rho_{M,s}-\rho_s}(h)\\
&=(\wt\chi  /\wt\chi_s )(\exp  2X)(\wt\chi_s^2)(\exp X)e^{<\rho_{s},X>}e^{<\rho_{M,s}-\rho_s,X>}\\
\end{aligned}
\end{equation}
But
$(\wt\chi/\wt\chi_s)(\exp  2X)=(\wt\chi/\wt\chi_s)^2(\exp X)$.     Also
$<\rho_{M,s},X> = <\rho _M,X> = 0$ since $X \in \z_M$.    Therefore
\begin{equation}
\begin{aligned}
\chi_0(h)&=\wt\chi^2(\exp X)\\
&=
\chi(h)e^{-<\rho,X>} \quad(\text{by }\eqref{e:cond}(a))\\
&=\chi(h)e^{<\rho_M-\rho, X>} =\chi_M(h).
\end{aligned}
\end{equation}
The result follows from the fact that $\oH_d^0=\overline{\exp}(\h_d\cap \z_M)(\oH\cap M_d)^0$.

For the final assertion $\mu=d\chi_s-2d\wt\chi_s-\rho$ and
$\mu_M=d\chi_{M,s}-2d\wt\chi_s-\rho_M$, and the result follows
immediately from \eqref{e:chiM}.
\end{proof}

Now fix lifting data $(\wt\chi_s,\chi_s)$ for $G$ (Definition
\ref{d:transfer}.)   By Part (1) of the Lemma $(\wt\chi_s,\chi_{M,s})$ is
lifting data for $(\wt M,M,\oM)$, and we use it to define transfer
factors $\Delta_{\oM}^{\tM}$ and 
$\Lift_{\oM}^{\tM}$.
The remainder of this section is devoted to proving that with these
choices, lifting commutes with parabolic induction.

Suppose $(h,\tg)\in X(\oH,\tH)$.
We need to
compare the transfer factors defined for $G$ and for $M$. Recall
\begin{subequations}
\renewcommand{\theequation}{\theparentequation)(\alph{equation}}  
\begin{align}
\label{transg1} 
\Delta_{\oG} ^{\tG}(h,\tg)&=
\frac{\epsilon_r(h,\Phi^+)\Delta^0(\Phi^+,h)\tilde\chi
(\tg)}{\epsilon_r(\tg,\Phi^+)\Delta^0(\Phi^+,\tg)\chi(h)}\\
\label{transm}
\Delta_{\oM}^{\tM}(h,\tg)&=
\frac{\epsilon_r(h,\Phi_M^+)\Delta^0(\Phi_M^+,h)\tilde\chi
(\tg)}{\epsilon_r(\tg,\Phi_M^+)\Delta^0(\Phi_M^+,\tg)\chi_M(h)}.
\end{align}
\end{subequations}

\begin{lemma}\label{last} Suppose $(h,\tg)\in X'(\oH,\tH)$. 
Then
\begin{equation}
\label{e:last}
\Delta_{\oG}^{\tG}(h,\tg) 
\frac{|D_{\tG}(\tg)|^{\frac12}}{|D_{\oG}(h)|^{\frac12}}
=
\Delta_{\oM }^{\tM }(h,\tg) 
\frac{|D_{\tM}(\tg)|^{\frac12}}{|D_{\oM}(h)|^{\frac12}}.
\end{equation}
\end{lemma}

\begin{proof} 
We first expand the left hand side, using \eqref{transg1} and 
\eqref{e:easy1} for the $|D|^{\frac12}$ terms. 
Let $\Phi ^+$ be a special set of positive roots.
The result is
\begin{equation}
\frac{
\epsilon_r(h,\Phi^+)\Delta^0(\Phi^+,h)|\Delta^0(\Phi^+,\tg)||e^\rho(\tg)|\wt\chi(\tg)
}
{
\epsilon_r(\tg,\Phi^+)\Delta^0(\Phi^+,\tg)|\Delta^0(\Phi^+,h)||e^\rho(h)|\chi(h)
}.
\end{equation}

Write $\Phi^+_r,\Phi^+_i$ and $\Phi^+_{cx}$ for the real, complex and
imaginary roots in $\Phi^+$, respectively. 
Then with the obvious
notation
\begin{equation}
\Delta ^0(h,\Phi ^+) =
\Delta ^0(h,\Phi_r^+)\Delta ^0(h,\Phi_{cx}^+)\Delta ^0(h,\Phi _i^+).
\end{equation}

Then
\begin{equation}
\begin{aligned}
\epsilon_r(h,\Phi^+)\Delta^0(\Phi^+_r,h)&=|\Delta^0(\Phi^+_r,h)|\\
\Delta^0(\Phi^+_{cx},h)&=|\Delta^0(\Phi^+_{cx},h)|\\
\end{aligned}
\end{equation}
Therefore
\begin{equation}
\label{e:fraction}
\frac{
\epsilon_r(h,\Phi^+)\Delta^0(\Phi^+,h)|\Delta^0(\Phi^+,\tg)|}
{
\epsilon_r(\tg,\Phi^+)\Delta^0(\Phi^+,\tg)|\Delta^0(\Phi^+,h)|}
=
\frac{
\Delta^0(\Phi_i^+,h)|\Delta^0(\Phi_i^+,\tg)|}
{\Delta^0(\Phi_i^+,\tg)|\Delta^0(\wt h,\Phi_i^+)|}.
\end{equation}
A similar argument applies to the right hand side of \eqref{e:last}. 
Since $\Phi_i=\Phi_{M,i}$ and the right hand side of
\eqref{e:fraction} only depends on the imaginary roots,
we conclude
\begin{equation}
\frac{
\epsilon_r(h,\Phi^+)\Delta^0(\Phi^+,h)|\Delta^0(\Phi^+,\tg)|}
{
\epsilon_r(\tg,\Phi^+)\Delta^0(\Phi^+,\tg)|\Delta^0(\Phi^+,h)|}
=
\frac{
\epsilon_r(h,\Phi^+_M)\Delta^0(\Phi^+_M,h)|\Delta^0(\Phi^+_M,\tg)|}
{
\epsilon_r(\tg,\Phi^+_M)\Delta^0(\Phi^+_M,\tg)|\Delta^0(\Phi^+_M,h)|}
\end{equation}
It remains to show
\begin{equation}
\frac{|e^{\rho}(\tg)|\wt\chi(\tg)}{|e^{\rho}(h)|\chi(h)}
=
\frac{|e^{\rho_M}(\tg)|\wt\chi(\tg)}{|e^{\rho_M}(h)|\chi_M(h)}.
\end{equation}
Recall
$\chi_M(h)=\chi(h)|e^{\rho_M-\rho}(h)|$.
Inserting this it  suffices to show
\begin{equation}
|e^{\rho-\rho_M}(\tg)|=|e^{\rho-\rho_M}(h)|^2.
\end{equation}
which follows from the fact that
$|e^{2\rho-2\rho_M}(h)|=|e^{\rho-\rho_M}(\phi(h))|=|e^{\rho-\rho_M}(\tg)|$.
\end{proof}

We can rewrite the Lemma in terms of orbits.

\begin{lemma}
\label{lasta} 
Suppose $\wt\O_{\tM}\in\Orb(\tM)$ 
and $\O_{\oM}^{\st}\in\Orbst(\oM)$ are  strongly regular and
semisimple, and satisfy
$\phi(\O_{\oM}^{\st})=p(\wt\O_{\tM})^{\st}$. 

Let $\wt\O_{\tG}$ be the unique $\tG$-orbit containing $\wt\O_{\tM}$
and let $\O^{\st}_{\oG}$ be the unique stable $\oG$-orbit containing
$\O^{\st}_{\oM}$. Then $\phi(\O^{\st}_{\oG})=p(\wt\O_{\tG})^{\st}$, and
$$
\Delta_{\oG}^{\tG}(\O_{\oG}^{\st},O_{\tG})
\frac{|D_{\tG}(\O_{\tG})|^{\frac12}}
{|D_{\oG}(\O_{\oG}^{\st})|^{\frac12}}=
\Delta_{\oM}^{\tM}(\O_{\oM}^{\st},\O_{\tM})
\frac
{|D_{\tM}(\O_{\tM})|^{\frac12}}
{|D_{\oM}(\O_{\oM}^{\st})|^{\frac12}}.
$$
\end{lemma}
 
\begin{lemma}
\label{resind} 
Suppose $\O_{\tG}\in\Orb(\tG)$ is strongly regular
and
$\O^{\st}_{\oM}\in\Orbst(\oM)$.
Let $\O^{\st}_{\oG}$ be the
unique stable orbit of $\oG$ containing $\O^{\st}_{\oM}$. Then
\begin{equation}
 \phi(\O^{\st}_{\oG})  = p(\O_{\tG})^{\st}  
\Leftrightarrow
\phi(\O^{\st}_{\oM}) = p(\O_{\tM})^{\st}
\end{equation}
for some $\tM$-orbit $\O_{\tM}\subset\O_{\tG}$.
Furthermore, if $\O_{\tG}$ is relevant, then  $\O_{\tM}$ is unique.
\end{lemma}

$$
\xymatrix{
\O_{\tG}\ar[d]&\ar@{_{(}-->}[l]\O_{\tM}\ar[d]\\
p(\O_{\tG})^{st}&\ar@{_{(}-->}[l]p(\O_{\wt M})^{st}\\
\O_{\oG}^{st}\ar[u]^\phi&\ar@{_{(}->}[l]\O_{\oM}^{st}\ar@{-->}[u]_\phi
}
$$

\begin{proof}
Choose $\overline h\in \oM$, $\tg\in\tG$ so that
$\O^{\st}_{\oM}=\O^{\st}(\oM,\overline h)$ and
$\O_{\tG}=\O(\tG,\tg)$.
Then $\O^{\st}_{\oG}=\O^{\st}(\oG,\overline h)$.
Let $h=\phi(\overline h)$ and $g=p(\tg)$.
Note that $p(\O_{\tG})=\O(G,g)$ and $\phi(\O_{\oM}^{\st})=\O^{\st}(M,h)$.

Except for the final statement the result does not involve $\tG$ or
$\tM$, and says
\begin{equation}
\label{e:says}
\O^{\st}(G,h) = \O^{\st}(G,g)  \Leftrightarrow
\O^{\st}(M,h) = \O^{\st} (M,y) 
\end{equation}
for some $y\in \O(G,g) \cap M$.
The implication $\Leftarrow$ is obvious: $y\in \O(G,g)$ and $y\in
\O^{\st}(M,h)$ implies $\O^{\st}(G,g)=\O^{\st}(G,y) = \O^{\st}(G,h)$.

Suppose $\O^{\st}(G,h) = \O^{\st}(G,g)$.
By \cite[Section 2]{shelstad_inner} 
\begin{equation}
\label{e:GMstable}
\O^{\st}(G,h)/G\simeq \O^{\st}(M,h)/M.
\end{equation}
More precisely if
$\O^{\st}(M,h)=\cup_i\O(M,h_i)$ with $h_i\in M$, then
$\O^{\st}(G,h)=\cup_i\O(G,h_i)$. Therefore $\O(G,g)=O(G,h_i)$ for some 
$h_i\in M$. Take $y=h_i$. This proves the implication $\Rightarrow$.

Assume that $\tilde g$ is relevant.  
By the preceding paragaraph clearly $\O^{\st}(M,h)$ contains a unique
$M$-orbit of the form $\O(M,y)$ with $y$ $G$-conjugate to $G$.
Choose $\wt y\in p\inv(y)$ so
that $\wt y$ is $\tG$-conjugate to $\tg$, i.e. $\O(\tM,\wt y)\subset
\O_{\tG}$.
It is clear that $\phi(\O^{\st}_{\oM}) = p(\O(\tM,\wt y))^{\st}$, and the
only other choice is 
$\O(\tM,-\wt y)$. By Lemma \ref{l:basicorbits}(3), since $\wt y$ is relevant, $\wt y$ is not
$\tG$-conjugate to $-\wt y$, so $\O(\tM,-\wt y)$ is not contained in $\O_{\tG}$.
\end{proof}

We continue to work with our given Levi subgroups $\tM$ and $\oM$.
Let $\oM \overline N $ and $\tM \tilde N $ be parabolic subgroups
of $\oG$ and $\tG$ with Levi components $\oM $ and $\tM $
respectively.  Let $\Theta_{\oM}$ be a stable character of $\oM $ and let
$\Theta_{\tM}$ be a character of $\tM $.  Then the induced
characters are independent of the choices of $\overline N $ and
$\tilde N$ so by abuse of notation we write
\begin{equation}
\label{induction} 
\Ind_{\oM }^{\oG}(\Theta_{\oM})=
\Ind_{\oM \overline N }^{\oG}(\Theta_{\oM} \otimes 1),
\ \ 
\Ind_{\tM}^{\tG}(\Theta_{\tM})=\Ind_{\tM\wt N}^{\tG}
(\Theta_{\tM}\otimes 1).
\end{equation}
   
\begin{theorem}
\label{liftind}
Let $\Theta_{\oM} $ be a stable character of
$\oM$. 
Then 
\begin{equation}
\Ind_{\tM}^{\tG}(\Lift_{\oM}^{\tM}(\Theta_{\oM}))=
\Lift_{\oG}^{\tG}(\Ind_{\oM}^{\oG}(\Theta_{\oM})).
\end{equation}
\end{theorem}

\begin{proof} 
Suppose $\O_{\tG}\in \Orb(\tG)$ is strongly regular and semisimple.
Let 
$\Res ^{\tG}_{\tM}\O_{\tG}=\{\O_{\tM}\,|\,\O_{\tM}\subset\O_{\tG}\}$.
Then for any character $\Theta _{\tM }$ of $\tM $, the induced
character is given by
(for example see \cite{hecht-schmid}) 
\begin{equation}
\Ind_{\tM}^{\tG}(\Theta_{\tM})(\O_{\tG})
=
|D_{\tG}(\O_{\tG})|^{-\frac12}
\sum_{\O_{\tM}\in \Res^{\tG}_{\tM}(\O_{\tG})}
|D_{\tM}(\O_{\tM})|^{\frac12}
\Theta_{\tM}(\O_{\tM}).
\end{equation}

Therefore by 
\eqref{e:liftorbits}
\begin{equation}
\label{e:sum1}
\begin{aligned}
&\Ind_{\tM}^{\tG}(\Lift_{\oM}^{\tM}(\Theta_{\oM}))(\O_{\tG})=
\sum_{\O_{\tM}}
\frac
{|D_{\tM}(\O_{\tM})|^{\frac12}}
{|D_{\tG}(\O_{\tG})|^{\frac12}}
\Lift_{\oM}^{\tM}(\Theta_{\oM})(\O_{\tM})\\
&=c_{\oM}^{-1}
\sum_{\O_{\tM}}
\sum_{\O_{\oM}^{\st}}
\frac{|D_{\tM}(\O_{\tM})|^{\frac12}}
{|D_{\tG}(\O_{\tG})|^{\frac12}}
\Delta_{\oM}^{\tM}(\O^{\st}_{\oM},\O_{\tM})
\Theta_{\oM}(\O^{\st}_{\oM}).
\end{aligned}
\end{equation}
The double sum is over
\begin{equation}
\{(\O_{\tM},\O^{\st}_{\oM})\,|\,\O_{\tM}\in\Res^{\tG}_{\tM}(\O_{\tG}),
\phi(\O_{\oM}^{\st}) = p(\O_{\tM})^{\st} \}
\end{equation}

Similarly, suppose $\O_{\oG}^{\st}\in\Orbst(\oG)$ is a semisimple orbit
and let
$\Res^{\oG}_{\oM}(\O^{\st}_{\oG})=\{\O^{\st}_{\oM}\,|\,\O^{\st}_{\oM}\subset \O^{\st}_{\oG}\}$.
It is well known that $\Ind_{\oM}^{\oG}(\Theta_{\oM})$ is stable, and
its has character formula
\begin{equation}
\Ind_{\oM}^{\oG}(\Theta_{\oM})(\O_{\oG}^{\st})
=
|D_{\oG}(\O^{\st}_{\oG})|^{-\frac12}
\sum_{\O^{\st}_{\oM}\in \Res^{\oG}_{\oM}(\O^{\st}_{\oG})}
|D_{\oM}(\O^{\st}_{\oM})|^{\frac12}
\Theta_{\oM}(\O^{\st}_{\oM}).
\end{equation}
This follows from the previous induced character formula and \eqref{e:GMstable}.
Therefore
\begin{equation}
\label{e:sum2}
\begin{aligned}
\Lift_{\oG}^{\tG}&(\Ind_{\oM}^{\oG}(\Theta_{\oM})(\O_{\tG})=
c_{\oG}^{-1}
\sum_{\O^{\st}_{\oG}}
\Delta^{\tG}_{\oG}(\O^{\st}_{\oG},\O_{\tG})
\Ind_{\oM}^{\oG}(\Theta_{\oM}))(\O^{\st}_{\oG})\\
&=
c_{\oG}^{-1}
\sum_{\O^{\st}_{\oG}}
\sum_{\O^{\st}_{\oM}}\frac
{|D_{\oM}(\O^{\st}_{\oM})|^{\frac12}}
{|D_{\oG}(\O^{\st}_{\oG})|^{\frac12}}
\Delta^{\tG}_{\oG}(\O^{\st}_{\oG},\O_{\tG})
\Theta_{\oM}(\O^{\st}_{\oM}).
\end{aligned}
\end{equation}
In this case the double sum is over
\begin{equation}
\{(\O^{\st}_{\oG},\O^{\st}_{\oM})\,|\,
\phi(\O^{\st}_{\oG}) = p(\O_{\tG})^{\st} ,
\O^{\st}_{\oM}\in\Res^{\oG}_{\oM}(\O^{\st}_{\oG})
\}
\end{equation}

Suppose that $\O_{\tG}$ is not relevant.  Then \eqref{e:sum2} is zero
since there are no orbits $\O^{\st}_{\oG}$ satisfying
$\phi(\O^{\st}_{\oG}) = p(\O_{\tG})^{\st} $.
Then for all $\O_{\tM}\in\Res^{\tG}_{\tM}(\O_{\tG})$, $\O_{\tM}$ is
not relevant, so \eqref{e:sum1} is also zero.

Assume $\O_{\tG}$ is relevant.  Since a maximally split Cartan
subgroup  $\oH_s$ of $\oM$ is also one for $\oG$, $c_{\oM}=c_{\oG}$.
By Lemma \ref{resind} both sums are over the same set of  orbits
$\O_{\oM}^{\st}$.
For each such orbit, the equality of the corresponding terms is then
given by  Lemma \ref{lasta}.
\end{proof}

By Theorem \ref{liftind}, we see that to understand lifting of
standard representations of $G$, it suffices to understand lifting of
discrete series representations of its cuspidal Levi subgroups.

\sec{Modified Character Data}
\label{s:characterdata}

We describe data which parametrizes L-packets for $\oG$, and also data
for irreducible representations of $\tG$. These are provided by  {\it
character data\/} of Vogan. For $\oG$ see \cite[Definition
6.6.1]{green}, and for $\tG$ see \cite[Section 2]{unitarizability}.
For our purposes it is
convenient to use a modified version of this data.

First let $G$ be the real points of a connected, complex reductive
group. Following \cite[Definition 6.6.1]{green} we define
a {\it regular character} of $G$ to be a triple
$\gamma=(H,\Gamma,\lambda)$ consisting of a $\theta$-stable Cartan subgroup $H$ of $G$, a
character $\Gamma$ of $H$, and an element $\lambda\in\h(\C)^*$.
We assume $\langle\lambda,\ch\alpha\rangle\in\R^\times$ for all
$\alpha\in\Phi_i$. Define $\Phi^+_i=\{\alpha\in\Phi^i\,|\,
\langle\lambda,\ch\alpha\rangle>0\}$ and let
$\Phi_{i,c}^+=\Phi_i\cap \Phi_{i,c}$ (cf.~Section \ref{s:notation}).
Define $\rho_i=\frac12\sum_{\Phi_i^+}\alpha$
and $\rho_{i,c}=\frac12\sum_{\Phi_{i,c}^+}\alpha$.
We assume
\begin{equation}
\label{e:dGamma}
d\Gamma=\lambda+\rho_i-2\rho_{i,c}.
\end{equation}
Write $H=TA$ and let $M=\Cent_G(A)$ as usual. 
Let $P$ be any parabolic subgroup containing $M$.
Associated to $\gamma = (H,\Gamma,\lambda)$ is a relative 
discrete series representation $\pi_M=\pi_M(\gamma)$ of $M$,
and a standard module $\pi_G(\gamma)=\pi(\gamma)=\Ind_P^G(\pi_M)$. 

The central character of $\pi_M(\gamma)$ (respectively $\pi(\gamma)$) 
is $\Gamma|_{Z(M)}$ (resp. $\Gamma|_{Z(G)}$). The Harish-Chandra
parameter of $\pi_M(\gamma)$ is $\lambda$, and the infinitesimal
character of $\pi_M(\gamma)$ and $\pi(\gamma)$ is $\lambda$.

The definition of character data is designed to make the lowest
K-types of $\pi(\gamma)$ evident. For example if $H$ is a compact
Cartan subgroup then $\Gamma$ is the highest weight  of a lowest
K-type of $\pi(\gamma)$. We are interested in character formulas, so
our needs are somewhat different. We modify this data appropriately.

\begin{definition}
\label{d:regchar}
A {\it modified regular character} for $G$ is a triple
$\gamma=(H,\Gamma,\lambda)$ where $H$ is a $\theta$-stable Cartan subgroup of $G$,
$\Gamma$ is a character of $H$, and $\lambda\in\h(\C)^* $. We assume
\begin{subequations}
\renewcommand{\theequation}{\theparentequation)(\alph{equation}}  
\end{subequations}
\begin{equation}
\langle\lambda,\ch\alpha\rangle\in\R^\times\quad\text{for all
}\alpha\in\Phi_i.
\end{equation}
Let $\Phi_i^+(\lambda )=\{\alpha\in\Phi_i\,|\, \langle\lambda,\ch\alpha\rangle>0\}$
and $\rho_i(\lambda)=\frac12\sum_{\alpha\in\Phi_i^+(\lambda)}\alpha$. We assume
\begin{equation}
d\Gamma=\lambda-\rho_i(\lambda ).
\end{equation}

Let $\cd(G)$ be the set of modified regular characters.
If $H$ is fixed let $\cd(G,H)$ be the set of modified
regular characters $(H,\Gamma,\lambda)$. 
If $H$ is given we write $(\Gamma,\lambda)=(H,\Gamma,\lambda)$.
\end{definition}

Write $H=TA$ and let $M=\Cent_G(A)$. Note that $\Phi_i$ is the set of
roots of $\h(\C)$ in $\mathfrak m(\C)$.
Since the definition of modified regular characters only
refers to the imaginary roots, $\cd(G,H)=\cd(M,H)$.

Associated to $\gamma$ is a relative discrete series representation
$\pi_M(\gamma)$, with Harish-Chandra parameter $\lambda$ and central
character $\Gamma|_{Z(M)}$. Let $\Theta_M(\gamma)$ be the character of
$\pi_M(\gamma)$. 
Then for $h\in H'$,
\begin{equation}
\label{e:dsM}
\Theta_M(\gamma)(h)=(-1)^{q_M}\Delta^0(\Phi_i^+(\lambda),h)\inv
\sum_{w\in W(M,H)}\epsilon(w)e^{w\rho_i(\lambda)-\rho_i(\lambda)}(h)\Gamma(w\inv h)
\end{equation}
where $q_M=\frac12\dim(M_d/K\cap M_d)$ and $\Delta^0$ is given by
\eqref{e:Deltas}. 

Let $P=MN$ be any parabolic subgroup containing $M$, and define
$\pi(\gamma)=\pi_G(\gamma)=\Ind_P^G(\pi_M(\gamma))$. Since we are interested only in
characters the choice of $P$ is not important.
Let $\Theta_G(\gamma)$ be the character of $\pi_G(\gamma)$:
\begin{equation}
\label{e:thetaG}
\Theta_G(\gamma)=\Ind_P^G(\Theta_M(\gamma)).
\end{equation}

Suppose $\gamma=(H,\Gamma,\lambda)$ is a modified character. 
Note that (with the obvious notation)
$2\rho_i(\lambda)-2\rho_{i,c}(\lambda)$
is a sum of roots, so $e^{2\rho_i(\lambda)-2\rho_{i,c}(\lambda)}$ is a well defined character
of $H(\C)$, and by restriction of $H$. This character is trivial on $Z(M)$.
It follows easily that $(H,\Gamma e^{2\rho_i(\lambda)-2\rho_{i,c}(\lambda)},\lambda)$ is
a regular character in the sense of Vogan, and 
defines the same relative discrete series representation of $M$ and
standard representation of $G$. This construction is clearly a
bijection between modified character data and character data in the
sense of Vogan.

There is a natural action of $G$ on $\cd(G)$ by conjugation.
The preceding bijection is $G$-equivariant.
By \cite{green} or \cite[Theorem 2.9]{unitarizability} we conclude:

\begin{lemma}
\label{l:cdg}
Suppose $\gamma,\gamma'\in \cd(G)$. Then $\Theta_G(\gamma)=
\Theta_G(\gamma')$ if and only if $\gamma=g\gamma'$ for some $g\in G$.  
\end{lemma}

We don't obtain every irreducible representation of $\oG$ this way,
since we don't allow ``limit'' characters, but this won't matter since
we are only interested in stable virtual characters.
Note that $W_i=W(\Phi_i)$ acts on $H$  and this induces an action on $\cd(G,H)$:
$w(H,\Gamma,\lambda)=(H,w\Gamma,w\lambda)$ where
$w\Gamma(h)=\Gamma(w\inv h)$.  Note that $\pi_G(w\gamma)$ and
$\pi_G(\gamma)$ have the same infinitesimal and central character.

\begin{definition}
\label{d:stable}
Suppose $\gamma=(H,\Gamma,\lambda)$ is a modified regular character.   
Let
\begin{equation}
\pi_M^{\st}(\gamma)=\sum_{w\in W(M,H)\backslash W_i}\pi_M(w\gamma).
\end{equation}
and
\begin{equation}
\pi_G^{\st}(\gamma)=\sum_{w\in W(M,H)\backslash W_i}\pi_G(w\gamma).
\end{equation}
Let $\Theta_M^{\st}(\gamma)$ and $\Theta_G^{\st}(\gamma)$ be the characters\
of $\pi_M^{\st}(\gamma)$ and $\pi_G^{\st}(\gamma)$, respectively.
\end{definition}

The character formula for $\Theta_M^{\st}(\gamma)$ on $H$ is the same as
\eqref{e:dsM} with $W(M,H)$ replaced by $W_i$.

\begin{lemma}[\cite{shelstad_inner}, Lemma 5.2]
\label{l:stable}
Fix $\gamma\in \cd(G,H)$. Then
$\Theta_M^{\st}(\gamma)$ and $\Theta_G^{\st}(\gamma)$ are stable
characters.

Suppose $\Theta$ is a stable virtual character of $G$. Then there
exist $\gamma_1,\dots, \gamma_n\in \cd(G)$ and integers $a_1,\dots,
a_n$ so that $\Theta=\sum a_i\Theta^{G,st}(\gamma_i)$.
\end{lemma}

Stable characters were defined before Definition \ref{d:lift}.
While the second part of the lemma is not stated in
\cite{shelstad_inner} it follows easily from the proof. See
Definition 18.9 and Lemmas 18.10 and 18.11 of \cite{abv}.

Now let $\tG$ be an admissible cover of $G$.  The definition of
regular character extends naturally to $\tG$.  See \cite[Section
27]{hc_ha_1} and \cite[Section 2]{unitarizability}.  In this setting
$\Gamma$ is an irreducible representation of $\tH$, or (by Lemma
\ref{adams3}) a character of $Z(\tH)$.

\begin{definition}
\label{d:regchar2}
A genuine modified regular character of $\tG$   is
a triple $(\tH,\wt\Gamma,\lambda)$.
Here $\tH$ is a Cartan subgroup of $\tG$, $\wt\Gamma$ is a genuine
one-dimensional representation of $Z(\tH)$, 
and $\lambda\in\h(\C)^*$. 
As in Definition \ref{d:regchar}.
we require
$\langle\lambda,\ch\alpha\rangle\in\R^\times$ for all $\alpha\in\Phi_i$,
and
$d\wt\Gamma=\lambda-\rho_i(\lambda)$.
Let $\cdg(\tG)$ be the set of genuine modified regular characters of
$\tG$, and $\cdg(\tG,\tH)$ the subset with given Cartan subgroup $\tH$.
\end{definition}

Let $H=TA$ be the image of $\tH$ in $G$, $M=\Cent_G(A)$, and
$\tM=p\inv(M)$. Associated to 
$\gamma=(\tH,\wt\Gamma,\lambda)$ 
is a relative discrete series representation $\pi_{\tM}(\gamma)$ of $\tM$.
The formula for the character $\Theta_{\tM}(\gamma)$ is
(cf.~\ref{e:dsM})
\begin{equation}
\label{dstM} 
\Theta_{\tM}(\gamma)(h)=
(-1)^{q_M}
\Delta^0(\Phi^+_i(\lambda),h)\inv
\sum_{w\in W(\tM,\tH)}\epsilon(w)e^{w\rho_i(\lambda)-\rho_i(\lambda)}(h)\Tr(\wt\tau(\wt\Gamma)(w^{-1}h))
\end{equation}
for $h\in\tH'$.
Here $\wt\tau(\wt\Gamma)$ is the irreducible representation of $\tH$
associated to $\wt\Gamma$ by Lemma
\ref{adams3}.

As before we also associate to $\gamma$ 
the standard module $\pi_{\tG}(\gamma)$ for $\tG$ induced from $\pi_{\tM}(\gamma)$.
 The central
character of $\pi_{\tM}(\gamma)$ is the restriction of $\wt\Gamma$ to the
center of $\tM$, and similarly for $\pi_{\tG}(\gamma)$.

\begin{lemma}
\label{l:cdgt}
Suppose $\gamma,\gamma'\in\cdg(\tG)$. Then
$\pi_{\tG}(\gamma)\simeq\pi_{\tG}(\gamma')$ if and only if
$\gamma'=g\gamma$ for some $g\in\tG$. Every irreducible genuine
representation of $\tG$ is isomorphic to $\pi_{\tG}(\gamma)$ for
some $\gamma\in\cdg(\tG)$.
\end{lemma}

\begin{proof}

The first statement is a special case of \cite[Theorem 2.9]{unitarizability}.

The second statement amounts to the fact that $\tM$ has no genuine limits of
(relative) discrete series.    This is because every irreducible genuine  representation
of $\tG$ is of the form $\pi_{\tG}(\gamma)$ where $\gamma$ is {\it final}
limit data as in \cite[Definition 2.4]{unitarizability}. For final
limit data we
allow $\langle\lambda,\ch\alpha\rangle=0$ if $\alpha$ is a noncompact
imaginary root.
But if $\wt\Gamma$ is genuine, then by  \cite[Lemma 6.11]{dualityonerootlength}
$\langle d\wt\Gamma,\ch\alpha\rangle\in\Z+\frac12$ for all noncompact
imaginary roots. This also holds for $\lambda$, so
$\langle \lambda,\ch\alpha\rangle\ne 0$ for all noncompact
imaginary roots. Therefore every genuine final limit data for $\tG$ is
in fact regular.
\end{proof}

There is no natural definition of stable distribution for $\tG$, and
we do not define analogues of $\pi^{\st}_G(\gamma)$ and
$\pi^{\st}_G(\gamma)$ for $\tG$.

We will make frequent use of formal sums of modified regular
characters (for example in Lemma \ref{l:stable}).
So if $\gamma_i\in \cd(G)$ and $a_i\in\Z$ ($i\le n$) we define
$\Theta_G(\sum_ia_i\gamma_i)=\sum_ia_i\Theta_G(\gamma_i)$, and similar
notation applies to $\tG$.

\sec{Formal Lifting of Modified Character Data}
\label{s:liftingofdata}

Fix an admissible triple $(\tG, G, \oG )$. Choose a maximally split
Cartan subgroup $H_s$ of $G$, and $(\tilde \chi _s, \chi _s) \in\caS(H_s)$.
This data determines transfer factors for $\oG$ and $\tG$
(Section \ref{s:transfer}) and $\Lift_{\oG}^{\tG}$ is defined (Section
\ref{s:lifting}).
Let $\mu=\mu(\wt\chi_s,\chi_s)\in\mathfrak z^*$ as in \eqref{e:mu}.

Suppose $H$ is a $\theta$-stable Cartan subgroup of $G$ 
and $\gamma=(\oH,\Gamma,\lambda)\in \cd(\oG,\oH)$ (Definition
\ref{d:regchar}). We define the lift of $\gamma$ to $\cd_g(\tG,\tH)$.
We first choose transfer factors for lifting from $\oH$ to
$\tH$. These differ from the obvious choice by a $\rho$-shift.

Write $H=TA$ and let $M=\Cent_G(A)$ as usual.  
Choose an arbitrary genuine character 
$\wt\chi$ of $Z(\tH)$ and a special set $\Phi^+$ of positive roots of
$H$ in $G$.
We assume
\begin{equation}
\label{e:PhiPhi_i}
\Phi^+\cap \Phi_i=\{\alpha\in\Phi_i\,|\,\langle\lambda,\ch\alpha\rangle>0\}.
\end{equation}

Let $\chi_{0}$ be the character of $\oH$ so that
$(\wt\chi,\chi_{0})\in\caS(H,\Phi^+,\wt\chi_s,\chi_s)$
(Definition \ref{d:compatible}).
Define
\begin{equation}
\label{e:shift}
\chi=|e^{\rho_i-\rho }|\chi_{0}.
\end{equation}
Let $\phi_H$ be $\phi$ restricted to $H$, and use it to define
\begin{equation}
\label{e:liftohth}
\Lift_{\oH}^{\tH}(\Gamma)
=
\Lift_{\oH}^{\tH}(\wt\chi,\chi,\Gamma)
\end{equation}
(cf. \ref{e:specify}).

\begin{lemma}
\label{l:chiind}
$\chi$ is independent of the choice of special positive roots $\Phi^+$
satisfying \eqref{e:PhiPhi_i}, and $\Lift_{\oH}^{\tH}$ is independent
of the choice of $\wt\chi$.
\end{lemma}

\begin{proof}
This follows from Lemmas \ref{l:elementary} and \ref{l:GammaPhi}.  
\end{proof}

\begin{definition}
\label{d:liftcharacterdata}
If $\Gamma|_{\Ker(\phi_H)}\ne \chi |_{\Ker(\phi_H)}$ then 
$\Lift_{\oH}^{\tH}(\Gamma)=0$, and we define $\Lift_{\oG}^{\tG}(\gamma)=\emptyset$ 
and $\Theta_{\tG}(\Lift_{\oG}^{\tG}(\gamma))=0$.

Otherwise write
\begin{equation}
\Lift_{\oH}^{\tH}(\Gamma)=\sum_{i=1}^n\wt\tau(\wt\Gamma_i)
\end{equation}
where $\wt\tau(\wt\Gamma_i)$ is given by Lemma \ref{adams3},
$n=|p(Z(\tH))/\phi(\oH)|$, and each $\wt\Gamma_i$ 
is a genuine one-dimensional representation of $Z(\tH)$.
For $1\le i\le n$ let $\wt\gamma_i=(\tH,\wt\Gamma_i,\frac12(\lambda-\mu))$;
this is a genuine regular character of $\tG$.
Define
\begin{equation}
\Lift_{\oG}^{\tG}(\gamma)=\{\wt\gamma_1,\dots,\wt\gamma_n\}.
\end{equation}
and
\begin{equation}
\Theta_{\tG}(\Lift_{\oG}^{\tG}(\gamma))=\sum_{i=1}^n\Theta_{\tG}(\wt\gamma_i).
\end{equation}
\end{definition}

It is important to keep in mind that $\Lift_{\oH}^{\tH}$ depends on
$\chi$, and therefore on $\lambda$ and $\gamma$ (cf.~Lemma \ref{l:chiind}).

The next Lemma follows from the definitions.

\begin{lemma}
\label{l:equiv}
    Let $\tilde g \in \tilde G$ and define $g= \overline p(p(\tilde g)) \in \oG$.    Then
$\Lift_{\oG}^{\tG}(g \gamma) = \tilde g \Lift_{\oG}^{\tG}( \gamma)$.
\end{lemma}

The following Lemma is an easy consequence of Lemma \ref{l:chiind} and Corollary \ref{gammaexists}.

\begin{lemma}
\label{l:dataexists}
Fix $\wt\gamma =(\tH,\wt\Gamma ,\lambda )\in \cdg(\tG)$.  Then there
is a unique $\gamma \in \cd(\oG)$ such that $\wt\gamma$ occurs in
$\Lift_{\oG}^{\tG}(\gamma)$.  It is given by $\gamma = (H, \Gamma,
2\lambda +\mu)$ where $\Gamma $ is the unique character of $\oH$ such
that $\wt\tau(\wt\Gamma)$ occurs in $\Lift_{\oH}^{\tH}(\Gamma)$ (cf.~Corollary \ref{gammaexists}).
\end{lemma}

The reason for the shift \eqref{e:shift} in the definition of
$\Lift_{\oH}^{\tH}$ is that it makes the next Lemma hold.
Fix a cuspidal Levi subgroup $M$ of $G$. Recall (Section 
\ref{s:levi})  $(\tM,M,\oM)$ is an admissible triple, where
$\tM=p\inv(M)$ and $\oM=\overline p(M)$, and $\Lift_{\oM}^{\tM}$ was
defined in Section \ref{s:levi} (following Lemma \ref{l:chims}).

\begin{lemma}
\label{l:indliftmliftg}
\begin{equation}
\Ind_{\tM}^{\tG}(\Theta_{\tM}(\Lift_{\oM}^{\tM}(\gamma)))
=
\Theta_{\tG}(\Lift_{\oG}^{\tG}(\gamma)).
\end{equation}
\end{lemma}

\begin{proof}
Fix a genuine character $\wt\chi$ of $Z(\tH)$.

Let $\chi_0$ be the unique character of $\oH$ so that
$(\wt\chi,\chi_0)\in\caS(H,\Phi^+,\wt\chi_s,\chi_s)$ and set
$\chi_1=\chi_0|e^{\rho_i-\rho }|$. 
Let 
$\mu=\mu(\tilde \chi _s,\chi_s) = d\chi_0-2d\wt\chi -\rho(\Phi^+)$.
Then the right hand side is the sum of terms
$\Ind_{\tM}^{\tG}(\Theta_{\tM}(\tH,\wt\Gamma,\frac12(\lambda-\mu))$
where $\wt\Gamma\in\Lift_{\oH}^{\tH}(\wt\chi,\chi_1,\Gamma)$.

On the other hand let $\chi_2$ be the unique character of $\oH$ so that
$(\wt\chi,\chi_2)\in\caS(H,\Phi_M^+,\wt\chi_s,\chi_{M,s})$.
Let $\mu_M=\mu(\tilde \chi _s,\chi_{M,s})=d\chi_2-2d\wt\chi-\rho(\Phi^+_i)$.
Then (cf. Lemma \ref{l:chims}) the left hand side is the sum of terms
$\Ind_{\tM}^{\tG}(\Theta_{\tM}(\tH,\wt\Gamma,\frac12(\lambda-\mu_M)))$
where $\wt\Gamma\in\Lift_{\oH}^{\tH}(\wt\chi,\chi_2,\Gamma)$.

It is enough to show $\chi_1=\chi_2$, i.e.
$(\wt\chi,\chi_0)\in\caS(H,\Phi^+,\wt\chi_s,\chi_s)$ implies
$(\wt\chi,\chi_0|e^{\rho_i-\rho }|)\in\caS(H,\Phi^+_M,\wt\chi_s,\chi_{M,s})$.
This is Lemma \ref{l:chims}(2).
\end{proof}
  
\begin{corollary}
Suppose $\wt\gamma=(\tH,\wt\Gamma,\wt\lambda)\in\cd_g(\tG,\tH)$.
Then there is a character $\Gamma$ of $\oH$ such that
$\gamma=(\oH,\Gamma,2\wt\lambda+\mu)\in\cd(\oG,\oH)$, and 
$\wt\gamma$ is a summand of $\Lift_{\oG}^{\tG}(\gamma)$.
\end{corollary}

\begin{proof}
This follows easily from Corollary \ref{gammaexists}: take
$\Gamma(h)=\chi(h)(\wt\Gamma/\wt\chi\inv)(\phi(h))$. We leave the
details to the reader.
\end{proof}

Given $\gamma=(\oH,\Gamma,\lambda)\in\cd(\oG,\oH)$ it is 
helpful to know when $\Lift_{\oG}^{\tG}(\gamma)$ is non-zero for some
choice of lifting data.

\begin{lemma}
We can choose lifting data so that $\Lift_{\oG}^{\tG}(\gamma)\ne 0$ if
and only if
\begin{equation}
\label{e:liftnonzero}
\begin{aligned}
\Gamma(h)=
\begin{cases}
(\wt\chi^2e^\rho)(h)&h\in \Ker(\phi_H)\cap \oH_d^0,\\  
\zeta_{cx}(h)&h\in \Gamma_r(\oH).
\end{cases}
\end{aligned}
\end{equation}
Here, as usual,  $\wt\chi$ is any genuine character of $\tH$, and
$\rho$ is defined with respect to a special set of positive roots
(Definition \ref{d:special}). 
\end{lemma}

\begin{proof}
Fix lifting data $(\chi_s,\tilde \chi_s)$, and 
let $\chi_0$ be the character of $\oH$ such that
$(\tilde\chi,\chi_0)\in S(H,\Phi^+,\tilde\chi_s,\chi_s)$. 
Then the lift is non-zero if $\Gamma (h) = |e^{\rho _i - \rho}(h)|
\chi_0(h)$ for all $h \in Ker(\phi_H)$.
It is easy to see $|e^{\rho_i-\rho}(h)|=1$ for all $h\in
\Ker(\phi_H)$, so the condition is
\begin{subequations}
\renewcommand{\theequation}{\theparentequation)(\alph{equation}}  
\begin{equation}
\Gamma(h)=\chi_0(h)\quad(h\in\Ker(\phi_H)).
\end{equation}
Using \eqref{e:cond} the condition in
\eqref{e:liftnonzero} is equivalent to
\begin{equation}
\label{e:liftnonzero2}
\Gamma(h)=\chi_0(h)\quad (h\in \oG_d\cap\Ker(\phi_H)).
\end{equation}
This condition is obviously necessary; we need to show it is
sufficient.

Let $\lambda$ be a character of $\oH$ satisfying
\begin{align}
\lambda(h)&=1\quad (h\in\oH\cap\oG_d)\\
\lambda(h)&=\Gamma(h)\chi_0\inv(h)\quad(h\in\Ker(\phi_H)).  
\end{align}
This is possible since
$\oH\cap\oG_d\cap\Ker(\phi_H)=\oG_d\cap\Ker(\phi_H)$, and
$\Gamma\chi_0\inv=1$ on this group by 
(b).
\end{subequations}
Then $\lambda$ extends uniquely to a one dimensional representation of
$\oH\oG_d$, and then (possibly not uniquely) to a one dimensional
representation $\psi$ of $\oG$.

The result follows from Lemma \ref{choicedata}:
the lift is non-zero if we replace
$\Delta_{\oG}^{\tG}(\chi_s,\wt\chi_s)$ with 
$\psi\Delta_{\oG}^{\tG}(\chi_s,\wt\chi_s)$.
\end{proof}

\sec{Lifting Stable Discrete Series}
\label{s:discreteseries}

Assume that $(\tilde G, G, \oG )$ is an admissible triple.  Fix
lifting data $(\tilde \chi _s, \chi _s)$ (Definition \ref{d:transfer}).
This data determines choices
of transfer factors for $\oG$ and $\tilde G$ (Section 6) and $\Lift
_{\oG}^{\tG}$ is defined (Section 7).  We assume for this section only
that $G$ has a relatively compact Cartan subgroup $H$, i.e. such that
$H\cap G_d$ is compact, 
and hence $\tilde
G, G$, and $\oG$ have relative discrete series representations.
We consider lifting of the  relative discrete series from $\oG$
to $\tG$. See Section \ref{s:exds} for the case when
$G=\oG$ is connected and semisimple.

Fix a relatively compact Cartan subgroup $H$ of $G$ with roots $\Phi
$.  Since $H$ is relatively compact $\Phi = \Phi _i$.
Fix $\gamma=(\oH,\Gamma,\lambda)\in\cd(\oG, \oH)$, and let 
$\Theta_{\oG}^{\st}(\gamma)$ be the associated stable 
discrete series character as in Section \ref{s:characterdata}.
We want to compute $\Lift_{\oG}^{\tG}(\Theta_{\oG}^{\st}(\gamma))$.

Define $\Lift_{\oH}^{\tH}(\Gamma)$ and
$\Lift_{\oG}^{\tG}(\gamma)=(\tH,\Lift_{\oH}^{\tH}(\Gamma),\frac12(\lambda-\mu))$ 
as in Definition
\ref{d:liftcharacterdata}. In this case by \eqref{e:liftohth}  
\begin{equation}
\Lift_{\oH}^{\tH}(\Gamma)=
\Lift_{\oH}^{\tH}(\wt\chi,\chi,\Gamma)
\end{equation}
where $(\wt\chi,\chi)\in\caS(H,\Phi^+(\lambda),\wt\chi_s,\chi_s)$
(since $\Phi=\Phi_i$, $\Phi^+$ is uniquely determined and
$\chi=\chi_0$).  Fix $(\wt\chi,\chi)\in
\caS(H,\Phi^+(\lambda),\wt\chi_s,\chi_s)$. 

Let $W=W(\Phi)=W_i$. 
For $w\in W$ let $w\gamma=(\oH,w\Gamma,w\lambda)$.
Then $\Lift_{\oH}^{\tH}(w\gamma)=(\tH,\Lift_{\oH}^{\tH}(w\Gamma),
\frac12w(\lambda-\mu))$ is defined.
It is important to keep in mind
that the lifting used in defining $\Lift_{\oH}^{\tH}(w\Gamma)$ depends
on $w$: by the preceding discussion
\begin{equation}
\Lift_{\oH}^{\tH}(w\Gamma)=\Lift_{\oH}^{\tH}(\wt\chi_w,\chi_w,w\Gamma)
\end{equation}
where $(\wt\chi_w,\chi_w)\in\caS(H,\Phi^+(w\lambda),\wt\chi_s,\chi_s)$.
One possibility is to take
\begin{equation}
\label{e:take}
(\wt\chi_w,\chi_w)=(\wt\chi,\chi e^{w\rho-\rho}).
\end{equation}
To see this, write $h=\gamma h_0$ where $\gamma\in \Gamma(\oH)$ and $h_0\in
\oH^0$. 
By Lemma \ref{l:GammaPhi} we can take
$\wt\chi_w=\wt\chi$ and
$\chi_w(h)=\chi(h)e^{w\rho-\rho}(h_0)$.
Since $\oH$ is relatively compact $\Gamma(\oH)\subset Z(\oG)$, so 
$e^{w\rho-\rho}(h_0)=e^{w\rho-\rho}(h)$. 

We will also use the fact that if $x\in W(\tG,\tH),y\in W$ then
\begin{equation}
\label{e:differentlift1}
\Lift_{\oH}^{\tH}(xy\Gamma)(\tilde h)=
\Lift_{\oH}^{\tH}(y\Gamma)(x\inv\tilde h),
\end{equation}
or more explicitly
\begin{equation}
\label{e:differentlift2}
\Lift_{\oH}^{\tH}(\wt\chi_{xy},\chi_{xy},xy\Gamma)(\tilde h)=
\Lift_{\oH}^{\tH}(\wt\chi_y,\chi_y,y\Gamma)(x\inv\tilde h).
\end{equation}
This follows  from the fact we can take
$(\wt\chi_y,\chi_y)$ as in \eqref{e:take}
and  $(\wt\chi_{xy},\chi_{xy})=(x\wt\chi_y,x\chi_y)$.
Then
\begin{equation}
\label{e:differentlift3}
\Lift_{\oH}^{\tH}(\wt\chi_{xy},\chi_{xy},xy\Gamma)(\tilde h)=
\Lift_{\oH}^{\tH}(x\wt\chi_y,x\chi_y,xy\Gamma)(\tilde h)\\
\end{equation}
and it follows immediately from the definitions that this equals
\begin{equation}
\Lift_{\oH}^{\tH}(\wt\chi_y,\chi_y,y\Gamma)(x\inv \tilde h).
\end{equation}

\begin{proposition}
\label{p:liftds}
Let $\gamma \in\cd(\oG, \oH)$.  Then
\begin{equation}
\label{e:propliftds}
\Lift _{\oG}^{\tG}(\Theta_{\oG}^{\st}(\gamma )) =
C_{\oG}(H) \sum _{w \in W(\oG, \oH)\backslash W} \Theta_{\tG}(\Lift _{\oG}^{\tG} (w\gamma )).  \end{equation}
\end{proposition}

\begin{proof}
Let $\overline\Theta=\Theta_{\oG}^{\st}(\gamma)$.
By \eqref{e:dsM}
\begin{equation}
\overline\Theta(h)=
(-1)^q\Delta^0( \Phi^+,h)^{-1}\sum_{w\in W}\epsilon(w)e^{w\rho-\rho}(h)\Gamma(w^{-1}h)\quad(h\in H')
\end{equation}
where $q=q_G$.
Fix a strongly regular element $\tilde h\in\tH$.
Since $H$ has no real roots we have
$$
\Delta_{\oG}^{\tG}(h,\tilde h) =\frac{\Delta ^0(\Phi ^+, h) \tilde
\chi (\tilde h)}{\Delta ^0(\Phi ^+, \tilde h) \chi ( h)}\quad(h\in X(\oH,\tilde h)).
$$
By definition of lifting we have 
$\Lift _{\oG}^{\tilde G}(\overline\Theta)(\tilde h)=$

\begin{equation}
c\inv(-1)^q 
\Delta ^0(\Phi ^+, \wt h)^{-1} 
\sum_{w \in W} \epsilon (w)
\tilde \chi (\tilde h)
\sum _{h \in X(\oH,\tilde h)} \chi (h)^{-1} e^{w\rho - \rho}(h) w\Gamma
(h).
\end{equation}

If $h\in X(\oH,\tilde h)$ then $e^{w\rho-\rho}(\tilde
h)=e^{2w\rho-2\rho}(h)$, so
$e^{w\rho - \rho}(h) = e^{w\rho - \rho}(\tilde h)e^{\rho - w\rho
}(h)$.
Therefore

\begin{equation}
\begin{aligned}
\Lift_{\oG}^{\tilde G}(\overline\Theta)(\tilde h)&=
c\inv(-1)^q 
\Delta ^0(\Phi ^+,\tilde h)^{-1}\sum_{w \in W} \epsilon
(w)e^{w\rho-\rho}(\tilde h)\times
\\
&\quad\quad\tilde\chi(\tilde h)
\sum_{h\in X(\oH, \tilde h)}\chi(h)^{-1}e^{\rho-w\rho}(h)w\Gamma(h).
\end{aligned}
\end{equation}
By \ref{abelian1}
\begin{equation}
\wt\chi(\tilde h)\sum_{h\in X(\oH,\tilde h)}\chi(h)^{-1}e^{\rho-w\rho}(h)w\Gamma(h)
=c(H)\Lift_{\oH}^{\tH}(\wt\chi,\chi e^{w\rho-\rho},w\Gamma)(\tilde h).
\end{equation}
By \eqref{e:take} 
$\Lift_{\oH}^{\tH}(\wt\chi,\chi e^{w\rho-\rho},w\Gamma)=\Lift_{\oH}^{\tH}(w\Gamma)$.
Recalling $C_{\oG}(H) = C(H)=c(H)/c$ this gives
\begin{equation}
\Lift_{\oG}^{\tG}(\overline\Theta)(\tilde h)=
C(H)(-1)^q\Delta^0(\wt h,\Phi^+)^{-1}
\sum_{w \in W} \epsilon (w)e^{w\rho-\rho}(\tilde h)
\Lift_{\oH}^{\tH}(w\Gamma)(\tilde h)\\
\end{equation}

Write the sum as
\begin{equation}
\sum_{y\in W(\tG,\tH)\backslash W}\epsilon (y)
e^{y\rho-\rho}(\tilde h)
\sum_{x\in W(\tG,\tH)}\epsilon (x)
e^{xy\rho-y\rho}(\tilde h)
\Lift_{\oH}^{\tH}(xy\Gamma)(\tilde h).
\end{equation}
By \eqref{e:differentlift2} 
\begin{equation}
\Lift_{\oH}^{\tH}(xy\Gamma)(\tilde
h)=\Lift_{\oH}^{\tH}(y\Gamma)(x\inv\tilde h).
\end{equation}
Also, for all $y\in W$ we have
\begin{equation}
\Delta^0(\wt h,\Phi^+)\inv\epsilon(y)e^{y\rho-\rho}=
\Delta^0(\wt h,y\Phi^+)\inv.
\end{equation}
This gives
\begin{equation}
\begin{aligned}
&\Lift _{\oG}^{\tilde G}(\overline\Theta)(\tilde h)=
C(H)(-1)^q\sum_{y\in W(\tG,\tH)\backslash W}
\Delta^0(y\Phi^+,\tilde h)^{-1}\times\\
&\quad\quad\sum_{x\in W(\tG,\tH)}\epsilon (x)
e^{xy\rho-y\rho}(\tilde h)
\Lift_{\oH}^{\tH}(y\Gamma)(x\inv\tilde h).
\end{aligned}
\end{equation}

By Proposition \ref{p:tori} 
$\Lift_{\oH}^{\tH}(y\Gamma)=0$ unless $y\in W(\gamma)$ where
\begin{equation}
W(\gamma)=\{w\in W\,|\, \chi(h)=e^{\rho-w\rho}(h)w\Gamma(h)\,\text{
  for all }h\in \Ker(\phi_H)\}.
\end{equation}
Fix $y\in W(\gamma)$ and define
\begin{equation} 
\wt\Gamma _y(\tilde h)=\tilde\chi(\tilde h)
\chi^{-1}(h)e^{\rho - y\rho}(h) y\Gamma (h) \quad (h \in X(\oH, \tilde h)).
\end{equation} 
Then by Proposition \ref{p:tori}
\begin{equation}
\Lift_{\oH}^{\tH}(y\Gamma)(\tilde h)=\sum_{\tilde\Gamma\in X_g(\tilde
H,\tilde\Gamma_y)}\Tr(\wt\tau(\tilde\Gamma)(\tilde h))\quad(\tilde h\in\tH).
\end{equation}

Thus
\begin{equation}
\begin{aligned}
&\sum _{x \in W(\tG, \tH)}
\epsilon (x) e^{xy\rho - y\rho}(\tilde h)\Lift _{\oH}^{\tH} ( y\Gamma
)(x^{-1}\tilde h)\\
&=
\sum _{\tilde \Gamma \in  X_g( \tH, \tilde \Gamma _y)}
\sum _{x \in W(\tG, \tH)} \epsilon (x) e^{xy\rho -
  y\rho}(\tilde h)\Tr(\wt\tau(\tilde \Gamma)(x^{-1}\tilde h)).
\end{aligned}
\end{equation}
But for each $\tilde \Gamma \in X_g( \tH, \tilde \Gamma _y)$ by
\eqref{dstM} we have
\begin{equation}
(-1)^q\Delta ^0(y\Phi ^+, \tilde h)^{-1}\sum _{x \in W(\tG, \tH)}
\epsilon (x) e^{xy\rho - y\rho}(\tilde h)\Tr(\wt\tau(\tilde \Gamma)(
x^{-1}\tilde h))= \Theta_{\tG}(\tilde \gamma )(\tilde h)\end{equation}
where $\tilde \gamma = (\tH, \tilde \Gamma ,(y\lambda - \mu )/2)$
is the corresponding element of $\Lift _{\oG}^{\tG}(y\gamma )$.  
Thus
\begin{equation}
Lift _{\oG}^{\tilde G}(\overline\Theta)(\tilde h)  =
C(H)\sum _{y \in W(\tilde G, \tH)\backslash W(\gamma) } 
\sum _{\tilde \gamma \in\Lift _{\oG}^{\tG}(y\gamma )}\Theta_{\tG}(\tilde \gamma )(\tilde h).
\end{equation}
In other words  for
every strongly regular element $\tilde h \in \tH$ we have
\begin{equation}
\label{e:liftds}
\Lift _{\oG}^{\tilde G}(\overline\Theta)(\tilde h) = 
C(H)\sum _{w \in W(\oG, \oH)\backslash W} 
\Theta_{\tG}(\Lift _{\oG}^{\tG} (w\gamma ))(\tilde h).
\end{equation} 

Since
$\overline\Theta $ is a stable discrete series character of $\oG$ with
Harish-Chandra parameter $\lambda $, it is a stable invariant
eigendistribution with infinitesimal character $\lambda $ and is
relatively supertempered.  Thus by Theorems 10.32 and 10.35 $\Lift
_{\oG}^{\tilde G}(\overline\Theta)$ is an invariant eigendistribution on $\tilde
G$ with infinitesimal character $(\lambda -\mu)/2$ and is relatively
supertempered.  Now each $\Theta_{\tilde G}(\Lift _{\oG}^{\tG}
(w\gamma )), w \in W$, is either zero or is a sum of discrete series
characters which all have infinitesimal character $(\lambda -\mu)/2$,
and hence is also a relatively supertempered invariant
eigendistribution on $\tilde G$ with infinitesimal character $(\lambda
-\mu)/2$.  Thus (\ref{e:liftds}) holds for all $\tilde g \in \tilde G'$ by
Harish-Chandra's Theorem \ref{hcsutemp}.
\end{proof}

Not all of the terms in \eqref{e:propliftds}   are non-zero.
For a precise description of which ones are non-zero see Proposition
\ref{p:Wsharp}.

\begin{example}
\label{ex:sl2ds}
Let $G=SL(2,\R)$. A discrete series representation of $\tG$ or $G$ is
determined by its Harish-Chandra parameter $\lambda\in\t^*(\C)$ where $T$
is a compact Cartan subgroup of $G$. 
Write $\pi_G^{\st}(\lambda)$ for the corresponding stable sum of discrete
series of $G$, and $\pi_{\tG}(\lambda)$ for a discrete series
representation of $\tG$. 
Let $P=\Z\langle\rho\rangle$ be the weight lattice.
We compute $C(T)=2$, and
then
\begin{equation}
\Lift_G^{\tG}(\pi_G^{\st}(\lambda))=
\begin{cases}
2\pi_{\tG}(\lambda/2)+2\pi_{\tG}(-\lambda/2)&\lambda\in2P+\rho\\
0&  \text{otherwise}
\end{cases}
\end{equation}
In the usual coordinates the lift is non-zero if and only if $\lambda$
is an odd integer, in which case $\pm\lambda/2\in \Z+\frac12$.

If $\oG=PSL(2,\R)\simeq SO(2,1)$ the result is the similar, except
that $C(T)=1$.
\end{example}

\sec {Lifting Standard Representations}
\label{s:standard}

Fix an admissible triple $(\tG,G,\oG )$ and lifting data
$(\wt\chi_s,\chi_s)$ for $\tG$.

Suppose $\Theta_{\oG}^{\st}$ is a stable virtual character of
$\oG$. Recall (Lemma \ref{l:stable}) that $\Theta_{\oG}^{\st}$ is a sum of standard
characters. Therefore, assuming we can write $\Theta_{\oG}^{\st}$ as 
a sum of standard modules, it is enough to compute the lifting of an
arbitrary standard character. This is the main result of this section
(Theorem \ref{t:liftstandard}).

We recall some constructions from Section \ref{s:characterdata}.
Suppose $\gamma$ is modified character data for $\oG$ (Definition
\ref{d:regchar}), with associated stable standard character
$\Theta_{\oG}^{\st}(\gamma)$.
Write $\gamma=(\oH,\Gamma,\lambda)$ and let $\oM$ be the corresponding
Levi factor of $\oG$. Recall $\gamma$ is also character data for
$\oM$, so the stable relative discrete series  character $\Theta_{\oM}^{\st}(\gamma)$ is
defined, and 
$\Theta_{\oG}^{\st}(\gamma)=\Ind_{\oM}^{\oG}(\Theta_{\oM}^{\st}(\gamma))$.
(As in the previous section we write $\Ind_{\oM}^{\oG}$ instead of
$\Ind_{\oM\overline N}^{\oG}$). 

Recall (Section \ref{s:characterdata}) $\Lift_{\oG}^{\tG}(\gamma)$
is a set
of modified character data
for $\tG$. Recall $W_i=W(\Phi _i)$ acts on $\cd(\oG,\oH)$.

\begin{theorem}
\label{t:liftstandard}
\begin{equation}
\label{e:liftstandard}
\Lift_{\oG}^{\tG}(\Theta_{\oG}^{\st}(\gamma))
=
C_{\oG}(H)\sum_{w\in W(\oM,\oH)\bs W_i}
\Theta_{\tG}(\Lift_{\oG}^{\tG}(w\gamma)).
\end{equation}
Explicitly for  $w\in W_i$ define 
$\Lift_{\oG}^{\tG}(w\gamma)$ by Definition \ref{d:liftcharacterdata}, 
and define
$\wt\gamma(w,i)\in \cdg(\tG)$ by:
\begin{equation}
\Lift_{\oG}^{\tG}(w\gamma)=\sum_{i=1}^{n_w}\wt\gamma(w,i)
\end{equation}
(if the lift is empty take $n_w=0$). Then
\begin{equation}
\label{e:liftwi}
\Lift_{\oG}^{\tG}(\Theta_{\oG}^{\st}(\gamma))
=
C_{\oG}(H)\sum_{w\in W(\oM,\oH)\bs W_i}\sum_{i=1}^{n_w}
\Theta_{\tG}(\wt\gamma(w,i)).
\end{equation}
\end{theorem}

\begin{proof}
This is merely a question of assembling the pieces.
Write $\sum_w$ for the sum over  $W(\oM,\oH)\bs W_i$. 
Then
\begin{equation}
\begin{aligned}
\Lift_{\oG}^{\tG}(\Theta_{\oG}^{\st}(\gamma))&=
\Lift_{\oG}^{\tG}(\Ind_{\oM}^{\oG}(\Theta_{\oM}^{\st}(\gamma)))
\quad(\eqref{e:thetaG} \text{ and Definition }\eqref{d:stable}\\
&=\Ind_{\tM}^{\tG}(\Lift_{\oM}^{\tM}(\Theta_{\oM}^{\st}(\gamma)))
\quad(\text{Theorem }\ref{liftind})
\\
&=C_{\oG}(H)\Ind_{\tM}^{\tG}(\sum_w
\Theta_{\tM}(\Lift_{\oM}^{\tM}(w\gamma))
\quad(\text{Proposition }\ref{p:liftds})
\\
&=C_{\oG}(H)
\sum_w\Theta_{\tG}(\Lift_{\oG}^{\tG}(w\gamma))\quad(\text{Lemma \ref{l:indliftmliftg}}).
\end{aligned}
\end{equation}
\end{proof}

It is important to know that the terms in the  sum
\eqref{e:liftwi} are distinct.

\begin{proposition}
In the sum \eqref{e:liftwi}, $\pi_{\tG}(\wt\gamma(w,i))\simeq
\pi_{\tG}(\wt\gamma(w',i'))$ if and only if
$W(\oM,\oH)w= W(\oM,\oH)w'$ and $i=i'$.
\end{proposition}

\begin{proof}
Suppose for $i=1,2$, $w_i\in W_i$, and  $\wt\gamma_i$ occurs in
$\Lift_{\oG}^{\tG}(w_i\gamma)$.
We may as well assume 
$w_1=1$, and write
$\wt\gamma_1=(\tH,\wt\Gamma_1,\frac12(\lambda-\mu))\in\Lift_{\oG}^{\tG}(\gamma)$, 
$\wt\gamma_2=(\tH,\wt\Gamma_2,\frac12(w\lambda-\mu))\in\Lift_{\oG}^{\tG}(w\gamma)$
for some
$w\in W_i$. We want to show
\begin{equation}
(\tH,\wt\Gamma_2,\frac12(w\lambda-\mu)=
g(\tH,\wt\Gamma_1,\frac12(\lambda-\mu))
\end{equation}
for some $g\in \tG$ implies $w=1$ and $\wt\Gamma_1=\wt\Gamma_2$.
Since $g$ normalizes $\tH$ 
let $u$ be the image of $\overline p(p(g))$
in $W(\oG,\oH)$ and set $v=w\inv u$.
Both $w$ and $u$ are in $W^\theta$,
the elements of $W$ fixed by $\theta$, so $v\in W^\theta$.
Also $u\lambda=w\lambda$, so $v\lambda=\lambda$ and therefore
$v\Phi^+_i(\lambda)=\Phi^+_i(\lambda)$. 

We now apply \cite[Propositions 3.12 and 4.16]{ic4}. 
By  Proposition 3.12(c)
$v\in W^\theta, v\Phi^+_i(\lambda)=\Phi^+_i(\lambda)$ implies 
$v\in W^q$, the Weyl group of the roots perpendicular to $\rho_i$.
Then by Proposition 4.16(a) $v\in W^q \subset W(\oG,\oH)$. Therefore 
$w=uv\inv\in W(\oG,\oH)$. But $w\in W_i$ so $w\in W(\oG,\oH)\cap
W_i=W(\oM,\oH)$.  

Since we are summing over cosets of $W(\oM,\oH)$ we may as well assume
$w=1$. Therefore 
$\wt\gamma_1=(\tH,\wt\Gamma_1,\frac12(\lambda-\mu))$, 
$\wt\gamma_2=(\tH,\wt\Gamma_2,\frac12(\lambda-\mu))$. Since these are both
contained in 
$\Lift_{\oG}^{\tG}(\gamma )$, $\wt\Gamma_1$ and $\wt\Gamma_2$
agree on $\tH^0$. Since they are conjugate by $g\in\tG$ they agree on
$Z(\tG)$. Since $\tH=Z(\tG)\tH^0$ (Proposition \ref{p:ZH}) this proves they are equal, proving
the Lemma.
\end{proof}

\begin{corollary}
In the setting of the Theorem, let $\{\wt\gamma_1,\dots, \wt\gamma_n\}$ be the set of
constituents of $\Lift_{\oG}^{\tG}(w\gamma)$ as $w$ runs over $W_i$,
considered without multiplicity.
Then
\begin{equation}
\label{e:liftwi2}
\Lift_{\oG}^{\tG}(\pi_{\oG}^{\st}(\gamma))
=
C_{\oG}(H)\sum_{i=1}^n \pi_{\tG}(\wt\gamma_i).
\end{equation}
\end{corollary}

\begin{remark}
It is a remarkable fact that not only are the constituents of
\eqref{e:liftwi} distinct, they have distinct central characters. 
See Remarks \ref{r:distinctps} and \ref{r:distinctds}.
As in \cite{adams:jussieu} one
can use this, together with Fourier inversion on $Z(\tG)$ to obtain a
character formula for $\Theta^{\tG}(\wt\gamma(w,i))$. Since this uses
a number of structural results not needed elsewhere we omit the proof.
\end{remark}

\begin{corollary}
\label{c:main}
Suppose $\pi$ is an admissible virtual representation of $\oG$ and
$\Theta_\pi$ is stable. Then
$\Lift_{\oG}^{\tG}(\Theta_\pi)$ is the character of a genuine virtual
representation of $\tG$, or $0$.
If $\pi$ is tempered and 
$\mu(\wt\chi,\chi_s)\in i\mathfrak z$ then so is 
$\Lift_{\oG}^{\tG}(\Theta_\pi)$.
\end{corollary}

Every genuine standard module for $\tG$ occurs in some Lift.

\begin{lemma}
\label{thetastexists}
Fix $\tilde\gamma \in \cdg(\tG)$.
Define $\gamma \in \cd(\oG)$ as
in Lemma \ref{l:dataexists}, so that $\wt\gamma$ occurs in
$\Lift_{\oG}^{\tG}(\gamma)$.  Then $\Theta _{\tG}(\tilde\gamma)$
occurs in $\Lift_{\oG}^{\tG}(\Theta_{\oG}^{\st}(\gamma))$.
Conversely, suppose  $\Theta_{\tG}(\tilde\gamma)$ occurs in
$\Lift_{\oG}^{\tG}(\Theta_{\oG}^{\st}(\gamma '))$ 
for some $\gamma'\in\cd(\oG)$.
Then $\Theta_{\oG}^{\st}(\gamma) = \Theta_{\oG}^{\st}(\gamma')$.
\end{lemma}

\begin{proof}  
The first statement is an immediate consequence of Theorem
\ref{t:liftstandard}.  Now suppose that $\gamma ' \in \cd(\oG)$ such
that $\Theta _{\tG}(\tilde\gamma)$ occurs in
$\Lift_{\oG}^{\tG}(\Theta_{\oG}^{\st}(\gamma '))$.  Then $\Theta
_{\tG}(\tilde\gamma) = \Theta _{\tG}(\tilde\gamma ')$ where $\tilde
\gamma '$ occurs in $\Lift_{\oG}^{\tG}(w\gamma ')$ for some $w \in
W_i'$, the imaginary Weyl group for the Cartan subgroup associated
to $\gamma '$.  Since $\Theta _{\tG}(\tilde\gamma) = \Theta
_{\tG}(\tilde\gamma ')$, by Lemma \ref{l:cdgt} there is $\tilde g
\in \tilde G$ such that $\tilde \gamma = \tilde g\tilde \gamma '$.
Define $g = \overline p(p(\tilde g))$.  Then by Lemma \ref{l:equiv}
$\tilde \gamma = \tilde g\tilde \gamma '$ occurs in $\tilde g
\Lift_{\oG}^{\tG}(w\gamma ') = \Lift_{\oG}^{\tG}(gw\gamma ')$.  But
by Lemma \ref{l:dataexists}, this implies that $\gamma = gw\gamma
'$.  Hence $\Theta_{\oG}^{\st}(\gamma) =
\Theta_{\oG}^{\st}(\gamma')$.
\end{proof}

We now make the sum \eqref{e:liftwi} more explicit.
Assume $\Lift_{\oG}^{\tG}(\Theta_{\oG}^{\st}(\gamma))\ne0$. Then 
$\Lift_{\oG}^{\tG}(w\gamma)\ne\emptyset$ for some $w\in W_i$.
and after replacing $\gamma$ with $w\gamma$ we assume $w=1$. 
We may then describe the set of $w\in W_i$ such that
$\Lift_{\oG}^{\tG}(w\gamma)$ is non-empty.
See \eqref{ex:Wsharpex} for a special case.

Fix $\gamma=(\oH,\Gamma,\lambda)\in\cd(\oG)$ and suppose
$\Lift_{\oG}^{\tG}(\pi_{\oG}^{\st}(\gamma))\ne 0$. 
Without loss of generality we may assume $\Lift_{\oG}^{\tG}(\gamma)\ne
0$. 
Then each component of $\Lift_{\oG}^{\tG}(\gamma)$ is of the form
$(\tH,\wt\Gamma,\frac12(\lambda-\mu))$. In particular
$\frac12(\lambda-\mu) - \rho _i(\lambda)$ is the differential of a genuine character of $\tH$.

Let
\begin{subequations}
\label{e:Wsharp}
\renewcommand{\theequation}{\theparentequation)(\alph{equation}} 
\begin{equation}
L=\{X\in \h\,|\, \exp(X)\in \Gamma(H)\cap C\cap H^0\}.
\end{equation}
Note that $\Gamma(H)\cap C\cap H^0$ is a finite central subgroup, 
and $L$ is a lattice.
Let
\begin{equation}
\Wsharp=\{w\in W_i\,|\,\exp((w\lambda/2-\lambda/2)(X))=1\text{ for all }X\in L\}.
\end{equation}

Here is an alternative description of  $\Wsharp$ in terms of
$M=\Cent_G(A)$. 
Let $\ch P_M$ be the weight lattice of the derived group of $M$, and
\begin{equation}
L'=\{\ch\gamma\in \ch P_M\,|\, \exp(2\pi i\ch\gamma)\in \Gamma(H)\cap
C\cap H^0\}.
\end{equation}
Note that
$X_{M_*}\subset L'\subset \ch P_M$ where $X_{M_*}=X_*(T(\C)\cap M_d(\C))$.
Then
\begin{equation}
\Wsharp=\{w\in W_i\,|\, \langle
w\lambda/2-\lambda/2,\ch\gamma\rangle\in \Z\text{ for all
}\ch\gamma\in L'\}.
\end{equation}
\end{subequations}

Suppose $\alpha\in\Phi_i$.
By \cite[Lemma 6.11]{dualityonerootlength})
and the fact that $\frac12(\lambda-\mu)- \rho _i(\lambda )$ is the differential of a
genuine character of $\tH$, 
we conclude
$\langle\lambda/2,\ch\alpha\rangle\in \Z$ or $\Z+\frac12$ depending on
whether $\alpha$ is compact or noncompact. It follows that
$\Wsharp\subset\Norm_{W_i}(\Phi_{i,c})$, and  
it is easy to see  $W(M,H)\subset\Wsharp$, so we have:
\begin{equation}
W(\Phi_{i,c})\subset W(M,H)\subset\Wsharp\subset \Norm_{W_i}(\Phi_{i,c})
\end{equation}
Compare \cite[Proposition 4.16(d)]{ic4}.

Frequently $\Gamma(H)\cap C\cap H^0=1$,
so $L$ is the kernel of $\exp$ restricted to $\h$ and $L'=X_{M_*}$.

\begin{proposition}
\label{p:Wsharp}
Fix $\gamma\in\cd(\oG)$ and suppose $\Lift_{\oG}^{\tG}(\gamma)\ne 0$. 
With $\Wsharp$ as in \eqref{e:Wsharp} we have
\begin{equation}
\label{e:Wsharpsum}
\Lift_{\oG}^{\tG}(\Theta_{\oG}^{\st}(\gamma))
=C_{\oG}(H)
\sum_{w\in W(M,H)\bs W_{\#}}
\sum_{\wt\gamma\in\Lift_{\oG}^{\tG}(w\gamma)}\Theta_{\tG}(\wt\gamma).
\end{equation}
There are $|W(M,H)\bs W_{\#}||Z_0(H)/\phi(\oH)|$ terms
in the sum.    
\end{proposition}

For example suppose $\gamma$ is a discrete series parameter for 
$G=SL(2,\R)$, and $\oG=SL(2,\R)$ or $PSL(2,\R)$. Then
$\phi(\oH)=Z_0(H)=H$, $W(G,H)=1$ and $W_\#=\Ztwo$. See Example
\ref{ex:sl2ds}. 

\begin{proof} 
Fix $w\in W_i$.
It is enough to show 
$\Lift_{\oG}^{\tG}(w\gamma)\ne 0$ if and only if $w\in\Wsharp$, and by 
Lemma \ref{l:indliftmliftg}
this is equivalent to $\Lift_{\oM}^{\tM}(w\gamma)\ne 0$.

By assumption $\Lift_{\oH}^{\tH}(\wt\chi,\chi,\Gamma)\ne\emptyset$ where
$(\wt\chi,\chi) \in S(H,\Phi ^+,\tilde \chi _s,\chi _s)$.
By Proposition \ref{p:tori} $\Gamma(h)=\chi(h)$ for all
$h\in\Ker(\phi_H)$ where $\phi_H$ is the restriction of $\phi$ to $\oH$.
Then
$w\gamma=(\oH,w\Gamma,w\lambda)$, and  (cf. \ref{e:take})
$\Lift_{\oH}^{\tH}(w\Gamma)=\Lift(\wt\chi,\chi_w,w\Gamma)$ where
$\chi_w(h)=e^{w\rho-\rho}(h)\chi (h)$.
Thus $\Lift_{\oG}^{\tG}(w\gamma)\ne 0$ if and only if
$w\Gamma(h)=e^{w\rho-\rho}(h)\chi(h)$ for all
$h\in\Ker(\phi_H)$. But $w\Gamma(h)=\Gamma(w\inv h)=\chi(w\inv h)$, so
the condition is
\begin{subequations}
\renewcommand{\theequation}{\theparentequation)(\alph{equation}}  
\begin{equation}
\chi((w\inv h)h\inv)=e^{w\rho-\rho}(h)\quad\text{for all }h\in \Ker(\phi_H).
\end{equation}
Write $\oH=Z(\oM)(\oH\cap\oM_d)^0$, and $h=zh_0$ accordingly.
Since $z\in Z(\oM)$ and $w\in W_i=W(M(\C),H(\C))$,
$(w\inv h)h\inv\in (\oH\cap \oM_d)^0$.
By \eqref{e:cond1} we can write condition (a) as 
\begin{equation}
(\tilde\chi^2e^\rho)((w\inv h)h\inv)=e^{w\rho-\rho}(h)\quad\text{for all }h\in \Ker(\phi_H).
\end{equation}
Let $H(\C)_2=\{h\in H(\C)\,|\, h^2=1\}$. 
Then $\Ker(\phi_H)=\overline p(H(\C)_2)\cap \oH$, so we can 
we can replace $\Ker(\phi_H)$ with $H(\C)_2\cap\overline
p\inv(\oH)$
on the right hand side of (b), which then becomes equivalent to
\begin{equation}
\tilde\chi^2((w\inv h)h\inv)=1\quad\text{for
 all }h\in H(\C)_2\cap\overline p\inv(\oH).
\end{equation}

We claim this is equivalent to
\begin{equation}
\tilde\chi^2((w\inv h)h\inv)=1\quad\text{for
 all }h\in H^0, h^2\in \Gamma(H)\cap C.
\end{equation}
\end{subequations}
If $h$ is in the set in (c), write $h=ta$ with $t\in\exp(\mathfrak
t(\C))$ and $a\in\exp(\mathfrak a(\C))=A(\C)$. Then $t$ is in the set in
(d), and $t=ha\inv$ with $a\in A(\C)$. Conversely for $h$ as in (d),
suppose $h^2=\exp(iX)\in\Gamma(H)$ with
$X\in \a$ (cf.~\ref{defzah}) and let $a=\exp(iX/2)\in A(\C)$. Then $x=ha$
satisfies the condition in (c).  The claim follows since $W_i$ acts
trivially on $A(\C)$.

For $h$ in the set in  (d) write $h=\exp(X/2)$ for $X\in \h$.
As in the discussion preceding the Proposition choose $\wt\chi$ to be
a genuine character of $Z(\tH)$ with differential
$\frac12(\lambda-\mu) - \rho _i$.   
Then  since $\mu\in\z(\C)^*$ 
\begin{equation}
\begin{aligned}
\wt\chi^2((w\inv h)h\inv)&=\exp((\lambda-\mu - 2\rho _i)(w\inv X/2-X/2))\\
&=\exp(w\lambda/2-\lambda/2)(X)\exp((2\rho_i-w2\rho_i)(X/2))\\
&=\exp(w\lambda/2-\lambda/2)(X)
\end{aligned}
\end{equation}
The last equality follows from 
$\exp((2\rho_i-w\rho_i)(X/2))=\exp((\rho_i-w\rho_i)(X))=1$ since $\exp(X)=h^2\in
\Gamma(H)\cap C\subset Z(G)$.
This gives \eqref{e:Wsharp}(b).
The alternative description of $\Wsharp$ is fairly standard. The
kernel of the exponential map restricted to $\h$ is contained in $\t$,
and now everything is taking place in $M$, and in fact in $M_d$.
\end{proof}

At least if  $\phi(\oH)=Z_0(H)$ the right hand side of 
\eqref{e:Wsharpsum}, which is a sum over $W(M,H)\backslash W_\sharp$, is
a reasonable candidate for an {\it L-packet} for $\tG$. 
Unlike the linear case the terms in this sum have different central
characters; this also happens for $Mp(2n,\R)$ \cite{metalift}.
However $W_\#$ may depend on  $\oG$, so this sum
is not canonical.
We also note that there is no obvious notion of stability for genuine
virtual characters of $\tG$. 
Therefore the precise definition  of L-packet for $\tG$
remains to be determined.

\sec{Appendix}
\label{app}

In this section we extend Hirai's Theorem on invariant
eigendistributions
to a class of groups
containing all reductive groups of Harish-Chandra's class.  We first
give some definitions from \cite{hirai_ie_2}.  Let $G$ be a reductive Lie
group with real Lie algebra $\g$. Let $\g_d(\C)$ be the derived
algebra of the complex Lie algebra $\g(\C)$.
Let $G^*(\C)$ be the connected complex adjoint group of $\g_d(\C)$. 
 
\begin{definition}\label{CA} $G$ satisfies condition A if 
the image of $G$ under the adjoint map $\Ad$ is contained in $G^*(\C)$.
\end{definition}

Recall \cite{green} $G$ is of Harish-Chandra's class if $G$ has
finitely many connected components, $Z(G_d)$ is finite, and condition
A holds.

The kernel of $\Ad$ is $\Cent_G(G^0)$, so let $\Ad(G)=G/\Cent_G(G^0)$.
Since $Z(G)\subset \Cent_G(G^0)$ there is a natural map
$p:G/Z(G)\rightarrow\Ad(G)$.

\begin{definition}\label{CB} $G$ satisfies condition B if it satisfies
condition A and there is connected, complex group $G^1(\C)$ with
adjoint group $G^*(\C)$ and an injective homomorphism $\phi$ making the following
diagram commute:

$$
\xymatrix{
G/Z(G)\ar[r]^\phi\ar[d]_p&G^1(\C)\ar[d]^\Ad\\
\Ad(G)\ar[r]^\Ad&G^*(\C)
}
$$ 
\end{definition}

For example $G$ satisfies condition B
if it satisfies condition A and $Z(G)=\Cent_{G}(G^0)$ (take
$G^1(\C)=G^*(\C)$). Note that  $GL(2,\R)$ satisfies
condition B, but any admissible two-fold cover of $GL(2,\R)$
(cf. Section \ref{s:admissible}) satisfies condition A (and is of
Harish-Chandra's class) but not condition B. 
The failure of nonlinear groups to satisfy condition B requires us to
extend Hirai's results.

Recall a connected complex Lie group is {\it acceptable} if $\rho$
(one-half the sum of the positive roots) exponentiates to a character of
a Cartan subgroup (cf.~Section \ref{s:notation}).

\begin{definition}\label{acceptable} 
$G$ is  acceptable if it satisfies condition A 
and there is an acceptable, connected, complex group $G^1(\C)$ with
adjoint group $G^*(\C)$ and a homomorphism $\phi$ making the following
diagram commute:
$$
\xymatrix{
G\ar[r]^\phi\ar[d]_p&G^1(\C)\ar[d]^\Ad\\
\Ad(G)\ar[r]^\Ad&G^*(\C)
}
$$ 
\end{definition}

This definition of acceptability for disconnected groups is very
strong.
For example, $GL(2, \R)$ is not an acceptable group (there is
no nontrivial homomorphism from $GL(2,\R)$ to $SL(2,\C))$, even 
though its identity component is acceptable.

Hirai's Theorem gives necessary and sufficient conditions for a class
function on $G'$ to be an invariant eigendistribution when $G$
satisfies condition B and is acceptable.  That the conditions are
necessary requires only condition A and acceptability, and is a fairly
straightforward extension of results of Harish-Chandra \cite{harish_inveigssgroup}.
Hirai's main contribution is to show that the conditions are
sufficient, and for this he requires condition B.

Let $G$ be a reductive Lie group satisfying condition A.  Then the
algebra of all left and right invariant differential operators on $G$
is equal to the center $\mathfrak Z$ of the universal envelopping
algebra.  In \cite{hirai_supplements} Hirai showed how to extend the results of
Harish-Chandra to show that when $\Theta $ is a $G^0$-invariant
eigendistribution for the action of $\mathfrak Z$, then $\Theta $ is
given by integration against a locally integrable function on $G$
which we also call $\Theta $.  Further, $\Theta $ is analytic on $G'$,
the set of regular semisimple elements, and is a $G^0$-invariant
function, that is $\Theta (xgx^{-1}) = \Theta (g), x \in G^0, g \in
G'$.
 
 We assume for the remainder of this appendix that $G$ is a reductive
Lie group satisfying condition A and use the notation of Section
\ref{s:hirai}.

 \begin{definition}\label{cthh} Let $\nu $ be a character of
$\mathfrak Z$ and let $\Theta $ be a $G^0$-invariant function on $G'$.
We say $\Theta \in IE(\nu )$ if $\Theta $ is the analytic function on
$G'$ corresponding to a $G^0$-invariant eigendistribution with
infinitesimal character $\nu $.
\end{definition}

Let $\Theta $ be a $G^0$-invariant function on $G'$.  Let $H$ be a
Cartan subgroup of $G$, $\Phi ^+$ a choice of positive roots, and use
the notation of Section \ref{s:lifting}.  Define
\begin{equation}\label{defpsi0} \Psi ^0(H ,\Phi ^+, h) = \Delta
^0(\Phi ^+,h) \Theta (h), \ \ \Psi (H ,\Phi ^+, h) = \epsilon _r(\Phi
^+, h) \Psi ^0(H ,\Phi ^+, h).\end{equation}
 
  Suppose that $G$ is acceptable with $\phi :G \rightarrow G^1(\C )$
satisfying the conditions of Definition \ref{acceptable}.  Let $\rho =
\rho (\Phi ^+)$ and write $\xi _{\rho}(h) = e^{\rho}(\phi(h)), h \in H$.
For $h \in H'$ define
    \begin{equation}\label{defpsirho} \Psi _{\rho}^0 (H ,\Phi ^+, h) =
\xi _{\rho}(h) \Psi ^0(H ,\Phi ^+, h), \ \ \Psi _{\rho} (H ,\Phi ^+,
h) = \xi _{\rho}(h)\Psi (H ,\Phi ^+, h).\end{equation} The functions
$\Psi _{\rho}^0 (H ,\Phi ^+)$ and $\Psi _{\rho}(H ,\Phi ^+)$ are the
ones used by Hirai to state his conditions.  They can only be defined
for acceptable groups.  We first restate Hirai's conditions in terms
of the functions $\Psi ^0 (H ,\Phi ^+)$ and $\Psi (H ,\Phi ^+)$ which
can be defined without the assumption of acceptability.
 
  For $X \in \h $, define $D_X$ and $D_X^{\rho}$ as in (\ref{defdx})
and (\ref{defdxr}).  The following lemma follows easily from the
definitions.
 
\begin{lemma}\label{lemdrho} Assume that $G$ is acceptable.  Let $F
\in C^{\infty}(H')$ and define $F_{\rho}(h) = \xi _{\rho}(h)F(h), h
\in H'$.  Then
$$DF_{\rho}(h) =  \xi _{\rho}(h) D^{\rho}F(h), h \in H', D \in S(\h (\C)).$$
\end{lemma}

Let $F:H' \rightarrow \C$, and let $\nu $ be a character of $\mathfrak
Z$.  Define conditions (C1, $\Phi ^+, \nu$) and (C2) as in Section
\ref{s:hirai}.  In the case that $G$ is acceptable, we say $F$
satisfies condition (CA1, $\nu$) if it satisfies (C1, $\Phi ^+, \nu$)
with $D^{\rho}$ replaced by $D$.
 
 \begin{lemma}\label{lemc2} \begin{enumerate}
\item If $\Psi (H ,\Phi ^+)$ satisfies condition (C1, $\Phi ^+, \nu$)
for one choice of positive roots, then $\Psi (H ,\Phi ^+)$ satisfies
condition (C1, $\Phi ^+, \nu$) for any choice of positive roots.

\item Assume that $G$ is acceptable.  Then $ \Psi _{\rho} (H ,\Phi
^+)$ satisfies condition (CA1, $\nu$) if and only if $\Psi (H ,\Phi
^+)$ satisfies condition (C1, $\Phi ^+, \nu$).
 \end{enumerate}
 \end{lemma}
 
 \begin{proof} (1) Let $\Phi ^+$ be any choice of positive roots.
Since $\epsilon _r(\Phi ^+,h)$ is locally constant on $H'$, $ \Psi (H
,\Phi ^+)$ satisfies condition (C1, $\Phi ^+, \nu$) if and only if $
\Psi ^0(H ,\Phi ^+)$ satisfies condition (C1, $\Phi ^+, \nu$). Let $w
\in W(\Phi )$ so that $w\Phi ^+$ is another choice of positive roots
for $\Phi $ with $\rho (w\Phi ^+) = w\rho $.  Then $\Delta ^0(w\Phi
^+,h) = \epsilon (w)e^{\rho - w\rho }(h) \Delta ^0(\Phi ^+,h), h \in
H$, where $\epsilon (w) = \pm 1$ is the determinant of $w$.  Thus
\begin{equation}\label{wphiplus} \Psi ^0(H ,w\Phi ^+, h) = \epsilon
(w)e^{\rho - w\rho }(h) \Psi ^0 (H , \Phi ^+, h), h \in
H'.\end{equation} Thus for any $X \in \h , h \in H'$,
$$D_X^{w\rho }\Psi  ^0(H ,w\Phi ^+, h) = \epsilon (w)e^{\rho  - w\rho }(h) D_X^{\rho} \Psi ^0 (H , \Phi ^+, h).$$
Thus $ \Psi ^0(H ,w\Phi ^+)$ satisfies condition (C1, $w\Phi ^+, \nu$)
if and only if $\Psi ^0(H ,\Phi ^+)$ satisfies condition (C1, $\Phi
^+, \nu$).
 
 (2) follows from Lemma \ref{lemdrho} and the fact that $\xi _{\rho
}(h)$ and $\xi _{\rho }^{-1} (h)$ are both real analytic on $H$.
\end{proof}

The proof of the following lemma is similar to that of Lemma
\ref{lemc2} since $\epsilon _r(\Phi ^+,h)$ is constant on each
connected component of $H'(R)$.

\begin{lemma}\label{lemc1} \begin{enumerate}
\item If $ \Psi (H ,\Phi ^+)$ satisfies condition (C2) for one choice
of positive roots, then $ \Psi (H ,\Phi ^+)$ satisfies condition (C2)
for every choice of positive roots.

\item Assume that $G$ is acceptable.  Then $ \Psi _{\rho} (H ,\Phi
^+)$ satisfies condition (C2) if and only if $ \Psi (H ,\Phi ^+)$
satisfies condition (C2).
 \end{enumerate}
\end{lemma}

Let $\alpha \in \Phi _r^+$ and define $J,\Phi _J^+ = c\Phi^+,H(\alpha
)$, and $\beta = c^*(\alpha) $ as in Section \ref{s:hirai}.  Let $h\in
H(\alpha )$ and assume that $ \Psi (H ,\Phi ^+)$ satisfies condition
(C1, $\Phi ^+, \nu$) and (C2).  Then $D_{\check \alpha }^{\rho, \pm
}\Psi (H, \Phi ^+,h)$ are defined as in (\ref{deflimpm}).  There are
also well-defined limits
\begin{equation}\label{dxb} D_{\check \alpha }^{\rho} \Psi ^0(H, \Phi
^+, h) = \lim _{s \rightarrow 0} D_{\check \alpha }^{\rho} \Psi ^0 ( H
, \Phi ^+, h \exp (s\check \alpha )).\end{equation} and when $G$ is
acceptable
\begin{equation}\label{dxa} D_{\check \alpha }\Psi _{\rho}^0(H, \Phi
^+, h) = \lim _{s \rightarrow 0} D_{\check \alpha }\Psi _{\rho}^0 ( H
, h \exp (s\check \alpha )).\end{equation} Further, $ h \in J'(R)$, so
that $D_{\check \beta }^{\rho _J}\Psi ^0( J, \Phi _J ^+, h)$ and
$D_{\check \beta }^{\rho _J} \Psi ( J, \Phi _J ^+, h)$ are defined,
and when $G$ is acceptable, $D_{\check \beta }\Psi _{\rho _J}^0( J,
\Phi _J ^+, h)$ is defined.

Define
\begin{equation}\label{epsilonb} \epsilon _r^{\alpha }(\Phi ^+, h) =
sign \prod _{\gamma \in \Phi _r^+, \gamma \not = \alpha } (1 -
e^{-\gamma}(h)), h \in H(\alpha ).
\end{equation}
    
\begin{lemma}\label{c3ao} For all $h \in H(\alpha )$,
$$D_{\check \alpha }^{\rho, \pm} \Psi (H, \Phi ^+, h) = \pm \epsilon _r^{\alpha }(\Phi ^+, h)D_{\check \alpha }^{\rho} \Psi  ^0(H, \Phi ^+, h);$$
$$D_{\check \beta }^{\rho _J} \Psi ( J, \Phi _J^+,  h) = \epsilon _r(\Phi _J^+, h)D_{ \check \beta }^{\rho _J} \Psi  ^0( J,  \Phi _J^+, h).$$
\end{lemma}

\begin{proof} Fix $h \in H(\alpha )$.  The first equation holds
because for all small $s \not = 0$,
$$\epsilon _r(\Phi ^+, h \exp (s\check \alpha )) = \epsilon _r^{\alpha }(\Phi ^+, h) sign (s).$$
The second equation holds because $h \in J'(R)$ so that $\epsilon
_r(\Phi _J^+, h)$ is constant in a neighborhood of $h$.
\end{proof}

\begin{lemma}\label{goodplus} Assume that $\Phi ^+$ is a choice of
positive roots such that $\alpha $ is a simple root for $\Phi _r^+$.
Then $\epsilon _r^{\alpha }(\Phi ^+, h) = \epsilon _r(\Phi _J^+, h)$
for all $h \in H(\alpha )$.
\end{lemma}

\begin{proof} Fix $h \in H(\alpha )$.  Write $\Phi _r^+ \backslash \{
\alpha \} = \Phi _1 \cup \Phi _2$ where
$$\Phi _1  = \{ \gamma \in \Phi _r^+: <\gamma , \check \alpha > = 0 \}$$ and $\Phi _2$ is its complement in
$\Phi _r^+ \backslash \{ \alpha \}$.  Then the positive real roots in
$\Phi _J$ are given by $\Phi _{r,J}^+ = c^*\Phi _1$.  Thus
$$\epsilon _r(\Phi _J^+, h) = sign \prod _{\gamma  \in \Phi _1} (1-e^{-c^*\gamma }(h)) =sign \prod _{\gamma  \in \Phi _1} (1-e^{- \gamma }(h))$$
since $e^{ c^*\gamma }(h) =e^{ \gamma }(h)$ for all $\gamma \in \Phi, h
\in H(\alpha )$.  Thus $$\epsilon _r^{\alpha }(\Phi ^+, h) = \epsilon
_r(\Phi _J^+, h) sign \prod _{\gamma \in \Phi _2} (1-e^{-\gamma
}(h)).$$ Let $\gamma \in \Phi _2$.  Then since $\alpha $ is simple for
$\Phi _r^+$, $\gamma \not = \alpha $, and $<\gamma , \check \alpha >
\not = 0$, $s_{\alpha }\gamma \in \Phi _2, s_{\alpha }\gamma \not =
\gamma $.  But $e^{-s_{\alpha }\gamma}(h) = e^{- \gamma}(h)$.  Thus
$$sign \prod _{\gamma \in \Phi _2} (1-e^{-\gamma }(h)) = 1.$$
\end{proof}

Consider the following conditions.
 \begin{equation}\label{dc3a0} D_{\check \alpha }\Psi _{\rho}^0(H,
\Phi ^+, h) = D_{\check \beta }\Psi ^0_{\rho _J}( J, \Phi _J^+, h), h
\in H(\alpha ).\end{equation}
 \begin{equation}\label{dc30} D_{\check \alpha }^{\rho}\Psi ^0(H, \Phi
^+, h) = D_{\check \beta }^{\rho _J}\Psi ^0 ( J, \Phi _J^+, h), h \in
H(\alpha ).\end{equation} Condition (\ref{dc3a0}) is the condition
used by Hirai in \cite{hirai_supplements} for acceptable groups, and (\ref{dc30})
is the analogous condition for not necessarily acceptable groups.  It
is useful for our applications in Section \ref{s:liftinginveig} to
replace (\ref{dc30}) by condition (C3) of Section \ref{s:hirai}.  The
following lemma clarifies the relationship between these conditions.

 \begin{lemma}\label{lemc3} \begin{enumerate}
 \item If $\Psi ^0(H,\Phi ^+)$ satisfies (\ref{dc30}) for one choice
of positive roots, then $\Psi ^0(H,\Phi ^+)$ satisfies (\ref{dc30})
for every choice of positive roots.
 
 \item Let $\Phi ^+$ be a choice of positive roots such that $\alpha $
is simple for $\Phi _r^+$.  Then $\Psi ^0(H,\Phi ^+)$ satisfies
(\ref{dc30}) if and only if $\Psi (H,\Phi ^+)$ satisfies condition
(C3).
  
\item Assume that $G$ is acceptable.  Then $\Psi _{\rho}^0(H,\Phi ^+)$
satisfies (\ref{dc3a0}) if and only if $\Psi ^0(H,\Phi ^+)$ satisfies
(\ref{dc30}).
  \end{enumerate}
   \end{lemma}

\begin{proof} Let $h \in H(\alpha )$ and let $\Phi ^+$ be any choice
of positive roots.

(1) Let $w \in W(\Phi )$.  Then $w_J = c^*w(c^*)^{-1} \in W(\Phi _J)$,
$c^*w\Phi ^+ = w_J\Phi _J^+$ and $w_J\rho _J = \rho (w_J\Phi _J^+)$.
Further,
$$D_{\check \alpha }^{w\rho }\Psi  ^0(H ,w\Phi ^+, h) = \epsilon (w)e^{\rho  - w\rho }(h) D_{\check \alpha }^{\rho} \Psi  ^0(H, \Phi ^+, h);$$
$$D_{\check \beta }^{w_J \rho _J }\Psi  ^0(J ,w_J\Phi _J^+, h) = \epsilon (w_J)e^{\rho _J  - w_J\rho _J}(h) D_{\check \beta }^{\rho _J} \Psi  ^0(J, \Phi _J^+, h).$$
But $\epsilon (w) = \epsilon (w_J)$ and $e^{\rho _J - w_J\rho _J}(h) =
e^{c^*(\rho -w\rho )}(h) =e^{ \rho -w\rho }(h)$.  Thus $\Psi ^0(H,w\Phi
^+)$ satisfies (\ref{dc30}) when $\Psi ^0(H, \Phi ^+)$ does.

(2) follows from combining Lemmas \ref{c3ao} and \ref{goodplus}.

(3) Using Lemma \ref{lemdrho} we have
$$D_{\check \alpha }\Psi _{\rho}^0(H, \Phi ^+, h) = \xi _{\rho}(h)D_{\check \alpha }^{\rho} \Psi  ^0(H, \Phi ^+, h);$$
$$D_{\check \beta } \Psi _{\rho _J} ^0 ( J,  \Phi _J^+, h) = \xi _{\rho _J}(h) D_{\check \beta }^{\rho _J}\Psi ^0 ( J,  \Phi _J^+, h).$$
Recall that $\xi _{\rho}, \xi _{\rho _J}$ are defined using a
homomorphism $\phi :G \rightarrow G ^1(\C )$.  Let $H^1(\C )$ and $J^1(\C
)$ be the Cartan subgroups of $G^1(\C )$ with Lie algebras $\h (\C )$
and $\j (\C )$.  Then we can pick a representative of the Cayley
transform $c \in G^1(\C )$ with $cH^1(\C )c^{-1} = J^1(\C )$ and
$chc^{-1} = h$ for all $h \in H^1(\C ) \cap J^1(\C )$.  Then for $h
\in H(\alpha )$,
 $$\xi _{\rho _J }( h ) = e^{\rho _J}(\phi ( h )) = e^{c^*\rho }(\phi( h )) = e^{\rho }(c\phi ( h )c^{-1}) = e^{\rho  }(\phi( h )) = \xi _{\rho }( h ), $$ 
 since $\phi ( h ) \in H^1(\C ) \cap J^1(\C )$.  Thus $\Psi
_{\rho}^0(H,\Phi ^+)$ satisfies (\ref{dc3a0}) if and only if $\Psi
^0(H,\Phi ^+)$ satisfies (\ref{dc30}).
 \end{proof}
 
\begin{definition}\label{cth} $\Theta $ satisfies condition $C(\nu )$
if for every Cartan subgroup $H$ of $G$ and choice of positive roots
$\Phi ^+$, $\Psi (H,\Phi ^+)$ satisfies conditions (C1,$\Phi ^+, \nu$)
and (C2), and $\Psi ^0 (H,\Phi ^+)$ satisfies (\ref{dc30}) for all
$\alpha \in \Phi _r^+$.
 \end{definition}
 
 Hirai's necessary and sufficient conditions are (CA1, $\nu$), (C2),
and (\ref{dc3a0}).  Using Lemmas \ref{lemc2}, \ref{lemc1}, and
\ref{lemc3}, when $G$ is acceptable, $\Theta $ satisfies $C(\nu )$ if
and only if $\Theta $ satisfies Hirai's conditions.  Thus we can state
Hirai's Theorem as follows.

\begin{theorem}\label{c123H} (Hirai, \cite{hirai_ie_2}) Let $G$ be a
reductive Lie group that satisfies condition B and is acceptable.  Let
$\Theta $ be a $G^0$-invariant function on $G'$ and let $\nu $ be a
character of $\mathfrak Z$.  Then $\Theta \in IE(\nu )$ if and only if
$\Theta $ satisfies $C(\nu )$.
 \end{theorem}
 
 \begin{remark} Suppose that $\Theta $ is a $G$-invariant function on
$G'$.  If $\Theta \in IE(\nu )$, then the corresponding
eigendistribution is $G$-invariant.  Thus $\Theta $ corresponds to a
$G$-invariant eigendistribution if and only if $\Theta $ satisfies
$C(\nu )$.
 \end{remark}
 
 \begin{remark} In Hirai's statement of Theorem \ref{c123H} in \S 11 of
\cite{hirai_ie_2}, he starts with the functions $\Psi _{\rho}(H,\Phi ^+)$
(which he calls $\kappa ^j$) on each of a set of representatives $H^j$
of $G^0$-conjugacy classes of Cartan subgroups.  He then gives an
extra condition that he calls $\epsilon$-symmetric which guarantees
that they can be patched together to give a $G^0$-invariant function
on $G'$ which is our $\Theta $.  We don't use this condition since we
assume from the beginning that we have a $G^0$-invariant function on
$G'$.
 \end{remark}
 
We will prove the following extension of Hirai's Theorem.
 
\begin{theorem}\label{c123a} Let $G$ be a reductive Lie group which
satisfies condition A.  Let $\Theta $ be a $G^0$-invariant function on
$G'$ and let $\nu $ be a character of $\mathfrak Z$.  Then $\Theta \in
IE(\nu )$ if and only if $\Theta $ satisfies $C(\nu )$.
 \end{theorem}
 
We will prove Theorem \ref{c123a} by making a number of reductions and
then applying Hirai's theorem.  A number of routine lemmas are stated
without proof.

Let $G$ be a reductive Lie group which satisfies condition A.  Let
$\Theta $ be a $G^0$-invariant function on $G'$ and let $\nu $ be a
character of $\mathfrak Z$.  Since $G$ is a Lie group, $G$ has at most
countably many connected components.  Write $G = \cup _i x_iG^0$, and
let $\Theta _i = \Theta \chi _i$ where $\chi _i$ is the characteristic
function of $x_iG^0$.  Then each $\Theta _i$ is a $G^0$-invariant
function on $G'$ which is supported on $x_iG^0$.
 
\begin{lemma}\label{reduction1} $\Theta \in IE(\nu )$ if and only if
$\Theta _i \in IE(\nu )$ for all $i$.  Further, $\Theta $ satisfies
$C(\nu )$ if and only if $\Theta _i$ satisfies $C(\nu )$ for all $i$.
\end{lemma}

\begin{lemma}\label{gammao} Each connected component of $G$ has a
representative $x$ such that $x^2 \in C_G(G^0)$ and $<x> \cap
C_G(G^0)G^0 \subset C_G(G^0)$ where $<x>$ denotes the cyclic subgroup
generated by $x$.  \end{lemma}

\begin{proof} Let $G^*(\C )$ denote the connected complex adjoint
group with Lie algebra $\g _d(\C )$.  Since $G$ satisfies condition A
there is a homomorphism $p:G \rightarrow G^*(\C)$ with kernel
$C_G(G^0)$ such that $Ad \ p(g) = Ad \ g, g \in G$.  In fact, $p(G)
\subset G^*(\R )$, the real points of $G^*(\C)$.  Fix $x \in G$.
Then $p(xG^0)$ is a connected component of $G^*(\R)$.  If $p(xG^0)
= (G^*)^0$, then $xG^0 \subset C_G(G^0)G^0$ so we can pick our
representative $x \in C_G(G^0)$.  Otherwise, $p(xG^0) = t(G^*)^0$
where $t ^2 = 1$ and $t \not \in (G^*)^0$.  We can pick our
representative $x$ so that $p(x) = t$.  Then $x \not \in C_G(G^0)G^0$,
but $p(x^2) = 1$ so $x^2 \in C_G(G^0)$.  Thus $<x> \cap C_G(G^0)G^0 =
<x^2> \subset C_G(G^0)$.
\end{proof}

\noindent {\bf First Reduction.}  Because of Lemmas \ref{reduction1}
and \ref{gammao} we may as well assume that $\Theta $ is supported on
$xG^0$ where $x \in G$ such that $x^2 \in C_G(G^0)$ and $<x> \cap
C_G(G^0)G^0 \subset C_G(G^0)$.

Fix $x \in G$ such that $x^2 \in C_G(G^0)$ and $<x> \cap C_G(G^0)G^0
\subset C_G(G^0)$.  Define $\tau :G^0 \rightarrow G^0$ by $\tau (g) =
xgx^{-1}, g \in G^0$.  Since $x^2 \in C_G(G^0)$ we have $\tau ^2 = 1$.
Let $G_{\tau } = <\tau > \ltimes G^0$ be the semidirect product of
$<\tau > = \{ 1, \tau \}$ and $G^0$.  That is, if $x \in C_G(G^0)$,
then $\tau = 1$ and $G_{\tau } = G^0$.  Otherwise $G_{\tau } = G^0
\cup \tau G^0$ has two connected components.

Define $\phi :xG^0 \rightarrow \tau G^0$ by $\phi (xg) = \tau g, g \in
G^0$.  Then $\phi $ is a diffeomorphism with $Ad(xg) = Ad(\phi (xg))$
for all $g \in G^0$ and $\phi (hxgh^{-1}) = h\phi (xg)h^{-1}$ for all
$h,g \in G^0$.  Thus there is a bijection between $G^0$-invariant
functions $\Theta $ on $G'$ which are supported on $xG^0$ and
$G^0$-invariant functions $\Theta _{\tau }$ on $G_{\tau}'$ which are
supported on $\tau G^0$ given by $\Theta (xg) = \Theta _{\tau }(\phi
(xg)), g \in G^0$.

\begin{lemma}\label{reducetotau} Let $\Theta $ be a $G^0$-invariant
function on $G'$ which is supported on $xG^0$, and let $\Theta _{\tau
}$ be the corresponding $G^0$-invariant function on $G_{\tau }'$ which
is supported on $\tau G^0$.  Then $\Theta \in IE(\nu )$ as a function
on $G'$ if and only if $\Theta _{\tau } \in IE(\nu )$ as a function on
$G_{\tau }'$.  Further, $\Theta $ satisfies $C(\nu )$ if and only if
$\Theta _{\tau }$ satisfies $C(\nu )$.
\end{lemma}

\begin{lemma}\label{gtau} $C_{G_{\tau }}(G^0) = Z(G^0)$.  \end{lemma}

\begin{proof} Clearly $C_{G_{\tau }}(G^0) \cap G^0 = Z(G^0)$.  Suppose
there is $g \in G^0$ such that $\tau g \in C_{G_{\tau }}(G^0)$.  Then
$\tau \in C_{G_{\tau}}(G^0)G^0$ so that $x \in <x> \cap C_G(G^0)G^0
\subset C_G(G^0)$ and $\tau = 1$.  Thus $\tau G^0 \cap C_{G_{\tau
}}(G^0) = \emptyset $ unless $\tau G^0 = G^0$.  In any case we have
$C_{G_{\tau }}(G^0) = C_{G_{\tau }}(G^0) \cap G^0 = Z(G^0)$.
\end{proof}

\noindent {\bf Second Reduction.}  Because of Lemma \ref{reducetotau}
we may as well replace $G$ by $G_{\tau }$.  Thus we assume there is $x
\in G$ with $x^2 = 1$ and $G = <x > \ltimes G^0$.  Moreover, because
of Lemma \ref{gtau} we have $C_G(G^0) = Z(G^0)$.  We continue to write
$\tau (g) = xgx^{-1}, g \in G^0$.

Let $C$ be a closed discrete subgroup of $Z(G^0)$ such that $\tau (C)
= C$.  Then $C$ is a normal subgroup of $G$.  Define $G_1 = G/C$ and
$x_1 = xC$.  Then $G_1^0 = G^0/C, x_1^2=1$, and $G_1 = <x_1> \ltimes
G_1^0$.

 Suppose that $\Theta _1$ is a $G_1^0$-invariant function on $G_1$.
Then $\Theta _1$ lifts to a $G^0$-invariant function $\Theta $ on $G'$

\begin{lemma}\label{cover1} $\Theta \in IE(\nu )$ if and only if
$\Theta _1 \in IE(\nu )$.  Further, $\Theta $ satisfies $C(\nu )$ if
and only if $\Theta _1$ satisfies $C(\nu )$.  \end{lemma}

\begin{lemma}\label{defc} Let $C = \{ \tau (z) z^{-1}: z \in
Z(G^0)\}$.  Then $C$ is a closed discrete subgroup of $Z(G^0)$ with
$\tau (C) = C$.  Define $G_1 = G/C$.  Then $Z(G_1) = C_{G_1}(G_1^0) =
Z(G_1^0)$.
\end{lemma}

\begin{proof} Define $\psi :Z(G^0) \rightarrow Z(G^0)$ by $\psi (z) =
\tau (z)z^{-1}$.  Then $\psi $ is a homomorphism so that $C = \phi
(Z(G^0))$ is a subgroup of $Z(G^0)$.  Write $\g = \g _d \oplus \z$
where $\g _d = [\g , \g]$ and $\z$ is the center of $\g$.  Then
$Z(G^0) = Z(G_d^0)Z$ where $Z = \exp \z$.  Now $Z(G_d^0)$ is a closed
discrete subgroup of $Z(G^0)$ and $\tau (z)=z$ for all $z \in Z$.
Thus $C = \{ \tau (z)z^{-1}: z \in Z(G_d^0)\} \subset Z(G_d^0)$ and so
is a closed discrete subgroup of $Z(G^0)$.  Further, since $\tau ^2 =
1$, we have $\tau (\psi (z)) = \psi (z^{-1}) \in C$ for all $z \in
Z(G^0)$.  Thus $\tau (C) = C$.

Write $p:G \rightarrow G_1 = G/C$ for the projection.  Then $Z(G_1^0)
= p(Z(G^0))$.  Let $z \in Z(G^0)$.  Then $\tau (z) z^{-1} \in C$ so
that $p(x)p(z)p(x)^{-1} = p(\tau (z)) = p(z)$.  Thus $p(z) \in
Z(G_1)$.  Thus $Z(G_1^0) \subset Z(G_1) \subset C_{G_1}(G_1^0)$.
Conversely, let $g \in G$ and suppose that $p(g) \in C_{G_1}(G_1^0)$.
Then $Ad \ g = Ad\ p(g) = 1$ so that $g \in C_G(G^0) = Z(G^0)$ by
Lemma \ref{gtau}.  Thus $p(g) \in Z(G_1^0)$.  Thus $C_{G_1}(G_1^0)
\subset Z(G_1^0)$.
\end{proof}

\begin{lemma}\label{ongx} Define $C$ as in Lemma \ref{defc}.  Let
$\Theta $ be a $G^0$-invariant function on $G'$ which is supported on
$xG^0$.  Then $$\Theta (gc) = \Theta (g), g \in G' , c \in C.$$ That
is, $\Theta $ factors to a $G_1^0$-invariant function $\Theta _1$ on
$G_1'$ supported on $x_1G_1^0$.
 \end{lemma}
 
\begin{proof} Since $\Theta $ is supported on $xG^0$ and $C \subset
G^0$, the result is true if $g \not \in xG^0$.  Let $g \in G^0$ with
$xg \in G'$ and let $c \in C$.  Then there is $z \in Z(G^0)$ such that
$c=\tau (z)z^{-1} = x zx^{-1}z^{-1} = x^{-1} zxz^{-1}$.  Since $\Theta
$ is $G^0$-invariant and $x^{-1} zx \in Z(G^0)$ we have
$$\Theta (xg) = \Theta (zxgz^{-1}) = \Theta (x x^{-1} zx gz^{-1}) = \Theta (xg x^{-1} zxz^{-1}) = \Theta (xgc).$$
 \end{proof}
 
 \noindent {\bf Third Reduction.}  By Lemmas \ref{cover1} and
\ref{ongx} we may as well assume that $G=G_1$.  But by Lemma
\ref{defc} we have $Z(G_1) = C_{G_1}(G_1^0) = Z(G_1^0)$.  That is, we
can assume that $Z(G) = C_G(G^0) = Z(G^0)$ and there is $x \in G$ with
$x^2 = 1$ and $G = <x > \ltimes G^0$.
  
The following Lemma completes the proof of Theorem \ref{c123a}.

\begin{lemma}\label{cover2} Let $G$ be a reductive Lie group
satisfying condition A.  Assume that $Z(G) = C_G(G^0) = Z(G^0)$ and
there is $x \in G$ with $x^2 = 1$ and $G = <x > \ltimes G^0$.  Let
$\Theta $ be a $G^0$-invariant function on $G'$ and let $\nu $ be a
character of $\mathfrak Z$.  Then $\Theta \in IE(\nu )$ if and only if
$\Theta $ satisfies $C(\nu )$.
\end{lemma}

\begin{proof} Write $\g = \g _d \oplus \z$ as in Lemma \ref{defc}.
Then $G^0 = G_d^0Z$ where $G_d^0$ and $Z$ are the connected subgroups
of $G$ corresponding to $\g _d$ and $\z$ respectively.  Since $G_d^0$
is a connected semisimple Lie group there is a connected finite
central extension $p_1:G_1 \rightarrow G_d^0$ such that $G_1$ is
acceptable.  That is, there is a connected acceptable complex Lie
group $G_2(\C)$ with Lie algebra $\g _d(\C )$ such that the inclusion
map of $\g _d \subset \g _d(\C )$ lifts to a homomorphism $\phi_1:G_1
\rightarrow G_2(\C )$.  Let $K = ker \ p_1 \cap ker\  \phi_1 \subset
Z(G_1)$ and define $\oG_1 = G_1/K$.  Then $p_1$ and $\phi_1$ both
factor to $\oG_1$, so that we have $\overline p_1:\overline
G_1 \rightarrow G_d^0$ and $\overline\phi _1:\oG_1 \rightarrow
G_2(\C )$.  Thus $\oG_1$ is also an acceptable cover of
$G_d^0$.  Thus we may as well assume that $ker \ p_1 \cap ker \ \phi_1 =
\{ 1 \}$.

Since $G$ satisfies condition A and $G_2(\C )$ is a complex Lie group
with Lie algebra $\g _d( \C)$ there is $x_2 \in G_2(\C )$ such that
$Ad \ x = Ad \ x_2$.  Then $Ad \ x_2^2 = Ad\ x^2 = 1$ so that $x_2^2
\in Z(G_2(\C ))$.  Thus $x _2$ has finite order.  Let $p_2:\tilde G_1
\rightarrow G_1$ be the simply connected cover of $G_1$.  Then there
is a unique automorphism $\tau $ of $\tilde G_1$ such that
 $$\tau (\widetilde \exp X) = \widetilde  \exp (Ad\  xX) =
\widetilde \exp (Ad \ x_2 X) , X \in \g .$$ It satisfies $\tau ^2 =
1$, and
 $$p_1p_2\tau (\tilde g) = x p_1p_2(\tilde g)x^{-1}, \ \  \phi_1p_2\tau (\tilde g) = x_2  \phi_1p_2(\tilde g)x_2^{-1}, \ \tilde g \in \tilde G_1.$$
Suppose that $\tilde g \in ker \ p_2$.  Then
$$p_1p_2\tau (\tilde g) = x p_1p_2(\tilde g)x^{-1} =1, \ \  \phi_1p_2\tau (\tilde g) = x_2 \phi_1p_2(\tilde g)x_2^{-1} =1.$$
Thus $p_2\tau (\tilde g) \in ker \ p_1 \cap ker\  \phi_1 = \{1 \}$, so
that $\tau (\tilde g) \in ker \ p_2$.  Thus $\tau $ descends to give a
well-defined automorphism $\tau :G_1 \rightarrow G_1$ satisfying $\tau
^2 = 1$ and
 \begin{equation}\label{p1b} p_1(\tau (g_1)) = x p_1(g_1)x^{-1}, \ \ \
 \phi_1(\tau (g_1)) = x_2 \phi_1(g_1)x_2^{-1}, g_1 \in G_1.  \end{equation}
Assume that $z_1 \in Z(G_1)$.  Then $p_1(z_1) \in Z(G_d^0) \subset
Z(G^0) = Z(G)$ by assumption.  Further, $Ad(\phi_1(z_1)) = Ad\ z_1 = 1$.
This implies that $\phi_1(z_1) \in Z(G_2(\C ))$ since $G_2(\C )$ is
connected.  Thus $p_1(\tau (z_1)) = xp_1(z_1)x^{-1} = p_1(z_1)$ and
$\phi_1(\tau (z_1)) = x_2\phi_1(z_1)x_2^{-1} = \phi_1(z_1)$.  Thus $\tau
(z_1)z_1^{-1} \in ker \ p_1 \cap ker \  \phi_1 = \{ 1 \}$.  That is, $\tau
(z_1) = z_1$ for all $z_1 \in Z(G_1)$.
 
Let $m = 2m_2$ where $m_2$ is the order of $x_2$ and write $\Z _m$ for
the additive group of integers mod $m$.  Let $\tilde G = \{( z, k, g):
z \in Z, k \in \Z _m , g \in G_1\}$.  For $z_1,z_2 \in Z, k,q \in
\Z_m, g_1,g_2 \in G_1$, define the product
\begin{equation}\label{mult} (z_1, k, g_1)(z_2, q, g_2) = (z_1z_2,
k+q, \tau ^{-q}(g_1)g_2). \end{equation} Thus $\tilde G$ is the direct
product of $Z$ with the semidirect product of $\Z _m$ and $G_1$ where
$k \in \Z _m$ acts on $G_1$ by $\tau ^k$.  $\tilde G^0$ is isomorphic
to the direct product of $Z$ and $G_1$ and so $\tilde G$ is a
reductive Lie group with Lie algebra $\g $.  We calculate that
\begin{equation}\label{p2b} (z , k, g )(z_2,q, g_2)(z , k, g )^{-1} =
(z_2, q, \tau ^k(\tau ^{-q}(g )g_2g ^{-1})).\end{equation} Thus $Ad(z,
k, g ) = Ad(x^kp_1(g))$ so that $\tilde G$ satisfies condition A.  We
will show that
\begin{equation}\label{centerg} Z(\tilde G) = C_{\tilde G}(\tilde G^0)
= \{ (z,k,g): z \in Z, x^k =1, g \in Z(G_1) \}.\end{equation} Using
(\ref{p2b}) and the fact that $\tau (g) = g$ for all $g \in Z(G_1)$,
it is clear that $ \{ (z,k,g): z \in Z, x^k =1, g \in Z(G_1) \}
\subset Z(\tilde G) \subset C_{\tilde G}(\tilde G^0)$.  Suppose that
$(z,k,g) \in C_{\tilde G}(\tilde G^0)$.  Then $Ad(x^kp_1(g)) =Ad(z, k,
g ) = 1$, so that $x^kp_1(g) \in C_G(G^0) = Z(G) = Z(G^0)$ by
assumption.  Thus $x^k \in <x> \cap G^0 = \{ 1\}$.  Thus $x^k=1$ and
$p_1(g) \in Z(G) \cap G_d^0 = Z(G_d^0)$ so that $g \in Z(G_1)$.  Thus
$C_{\tilde G}(\tilde G^0) \subset \{ (z,k,g): z \in Z, x^k =1, g \in
Z(G_1) \} \subset Z(\tilde G)$.  This shows that $\tilde G$ satisfies
condition B.

Define $\phi :\tilde G \rightarrow G_2(\C )$ by $\phi (z, k, g ) = x_2^k
\phi_1(g)$.  It is well-defined because $m$ is a multiple of the order of
$x_2$.  It is easy to show it is a homomorphism using (\ref{p1b}).  It
is a continuous homomorphism which induces canonically the natural
injection of Ad($\tilde G$) into Ad ($G_2(\C )$).  Thus $\tilde G$ is
acceptable.

Define $p :\tilde G \rightarrow G$ by $p (z, k,g) = zx^kp_1(g)$.  The
mapping is well-defined since $m$ is even, and it satisfies $Ad\ p (z,
k,g) = Ad (z, k,g)$.  It is easy to show that it is homomorphism using
(\ref{p1b}), and it is clearly surjective.  Let $(z,k,g) \in ker \ p$
so that $zx^kp_1(g) =1$.  Then $Ad(z,k,g) = Ad p (z, k,g) = 1$ so that
$(z,k,g) \in C_{\tilde G}(\tilde G^0) = Z(\tilde G)$.  Thus $ker \ p
\subset Z(\tilde G)$.  Further, using (\ref{centerg}) we have
$Z(\tilde G )^0 = \{ (z,0,1): z \in Z\}$.  Thus the identity component
of $ker \ p$ is contained in $Z(\tilde G )^0 \cap ker \ p = \{ (1,0,1)
\}$ so that $ker \ p$ is a discrete subgroup of $Z(\tilde G)$.
   
Let $\Theta $ be a $G^0$-invariant function on $G'$ and let $\nu $ be
a character of $\mathfrak Z$.  Lift $\Theta $ to a $\tilde
G^0$-invariant function $\tilde \Theta $ on $\tilde G'$.  Then $\Theta
\in IE(\nu )$ if and only if $\tilde \Theta \in IE(\nu )$ and $\Theta
$ satisfies $C(\nu )$ if and only if $\tilde \Theta $ satisfies $C(\nu
)$.  But $\tilde G$ is acceptable and satisfies condition B.  Thus
$\tilde \Theta \in IE(\nu )$ if and only if $\tilde \Theta $ satisfies
$C(\nu )$ by Theorem \ref{c123H}.
\end{proof}

\bibliographystyle{plain}
\def\cprime{$'$} \def\cftil#1{\ifmmode\setbox7\hbox{$\accent"5E#1$}\else
  \setbox7\hbox{\accent"5E#1}\penalty 10000\relax\fi\raise 1\ht7
  \hbox{\lower1.15ex\hbox to 1\wd7{\hss\accent"7E\hss}}\penalty 10000
  \hskip-1\wd7\penalty 10000\box7}

\end{document}